\documentclass[12pt]{amsart}

\usepackage{amsfonts, amssymb, latexsym}
\usepackage{amsmath, amsthm}

\usepackage{color}
\definecolor{Gray}{cmyk}{0,0,0,0.50}
\definecolor{Black}{cmyk}{0,0,0,1}

\usepackage{geometry}
\input prepictex
\input pictexwd
\input postpictex

\newtheorem{theorem}{Theorem}[section]
\newtheorem{lemma}[theorem]{Lemma}
\newtheorem{prop}[theorem]{Proposition}
\newtheorem{corollary}[theorem]{Corollary}
\newtheorem{remark}[theorem]{Remark}
\newtheorem{example}[theorem]{Example}
\newtheorem{defn}[theorem]{Definition}
\numberwithin{equation}{subsection}

\newcommand\bbC{\mathbb{C}}

\newcommand\bbK{\mathbb{K}}
\newcommand\bbN{\mathbb{N}}
\newcommand\bbQ{\mathbb{Q}}
\newcommand\bbR{\mathbb{R}}

\newcommand\bbZ{\mathbb{Z}}

\newcommand\calA{\mathcal{A}}
\newcommand\calB{\mathcal{B}}

\newcommand\calE{\mathcal{E}}

\newcommand\calK{\mathcal{K}}
\newcommand\calL{\mathcal{L}}

\newcommand\fh{\mathfrak{h}}

\newcommand\fsl{\mathfrak{sl}}

\newcommand\fS{\mathfrak{S}}

\newcommand\GL{\mathrm{GL}}
\newcommand\hgt{\mathrm{ht}}

\newcommand\shf{\mathrm{sh}}

\newcommand\rmasc{\mathrm{asc}}
\newcommand\rmdes{\mathrm{des}}

\newcommand\ttA{\mathtt{A}}
\newcommand\ttC{\mathtt{C}}
\newcommand\ttb{\mathtt{b}}

\newcommand\ttd{\mathtt{d}}
\newcommand\tte{\mathtt{e}}
\newcommand\ttf{\mathtt{f}}
\newcommand\ttH{\mathtt{H}}
\newcommand\tti{\mathtt{i}}

\newcommand\ttwt{\mathtt{wt}}

\title{A Littlewood-Richardson rule for Macdonald polynomials}
\author{Martha Yip \\
Department of Mathematics\\
University of Wisconsin\\
Madison, WI 53706}
\date{\today}

\begin{document}

\begin{abstract}
Macdonald polynomials are orthogonal polynomials associated to root systems, and in the type A case, the symmetric Macdonald polynomials are a common generalization of Schur functions, Macdonald spherical functions, and Jack polynomials.
We use the combinatorics of alcove walks to calculate products of
monomials and intertwining operators of the double affine Hecke algebra.  From this, we obtain a product formula for Macdonald polynomials.
\end{abstract}
\maketitle \tableofcontents
\parskip=9pt

\section{Introduction}

In~\cite{M88}, Macdonald introduced a remarkable family of orthogonal polynomials $P_\lambda(q,t)$ associated with root systems.  For special values of $q$ and $t$, they specialize to various well-known functions, including Weyl characters and spherical functions for $p$-adic groups.  These polynomials are a basis for symmetric functions, and are a common generalization of Schur functions $s_\lambda$, monomial symmetric functions, Hall-Littlewood polynomials, and the symmetric Jack polynomials.  The symmetric Macdonald polynomials are indexed by dominant weights of the weight lattice $P$.

Classically, the Littlewood-Richardson coefficients $c_{\lambda\mu}^\nu$ are the structure constants of the ring of symmetric functions with respect to the Schur basis:
$$s_\lambda s_\mu = \sum_{\nu} c_{\lambda\mu}^\nu s_\nu.$$
In the representation theory of the general linear group $\GL_n(\bbC)$, the Littlewood-Richardson coefficients also give the multiplicity of the irreducible highest weight module $V(\nu)$ in $V(\lambda)\otimes V(\mu)$.  The coefficient $c_{\lambda\mu}^\nu$ is given combinatorially as the number of Young tableaux of shape $\nu\backslash \lambda$ admitting a Littlewood-Richardson filling of type $\mu$.

Littelmann introduced the path model in~\cite{Li94} as a tool for calculating formulas for characters of complex symmetrizable Kac-Moody algebras, and showed that it can also be used to compute Littlewood-Richardson coefficients.  Instead of a sum over tableaux, his formula for $c_{\lambda\mu}^\nu$ is a sum over certain paths in the vector space $P\otimes_\bbZ \bbR$, where the endpoint (weight) of a path takes the place of the filling of a tableau.  
Several variations of the Littelmann path model were introduced to obtain character formulas, including the gallery model of Gaussent-Littelmann~\cite{GL02}, and the model of Lenart-Postnikov~\cite{LP08} based on $\lambda$-chains.  In~\cite{R06}, Ram developed the alcove walk model for working in the affine Hecke algebra, and the paper~\cite{RY} showed that alcove walks are a useful tool for expanding products of intertwining operators of the double affine Hecke algebra. 

Cherednik developed the theory of double affine Hecke algebras, using it to solve Macdonald's constant term conjectures~\cite{C95a}, and in~\cite{C95b}, he showed that products of intertwining operators of the double affine Hecke algebra generate the nonsymmetric Macdonald polynomials $E_\lambda(q,t)$, which are a family of orthogonal polynomials indexed by points of the weight lattice.  These polynomials were first introduced by Opdam~\cite{O95} in the case $q\rightarrow 1$ (see~\cite[p.147]{M03}).  By applying a symmetrizing operator $\mathbf{1}_0$ to $E_\lambda$, one can obtain the symmetric polynomials $P_\lambda$.

In this paper, we use the alcove walk model to calculate products of monomials and intertwining operators of the double affine Hecke algebra (Theorem~\ref{thm.X-tau}), and give a product formula for two symmetric Macdonald polynomials (Theorem~\ref{thm.PP-P}).  This is a  generalization of the classical formula for products of Weyl characters, where the generalized Littlewood-Richardson coefficents are rational functions in $q$ and $t$.  In particular at $q=0$, Theorem~\ref{thm.PP-P} reduces to Schwer's product formula~\cite[Theorem 1.3]{Sc06} for Hall-Littlewood polynomials in terms of positively folded galleries, and at $q=t$, the formula reduces to the product formula of Littelmann~\cite{Li94} for Weyl characters phrased in terms of Littelmann paths.

Section~\ref{sec.daha} of this paper introduces the basic definitions and properties of double affine Hecke algebras for reduced root systems.  The alcove walk statistics needed in the later sections are addressed here.  Section~\ref{sec.poly} discusses how alcove walks can be used to calculate the coefficients of products of monomials and intertwining operators in the double affine Hecke algebra.
In Section~\ref{sec.lrrule}, we state and prove the main results:

\noindent{\bf Theorem~\ref{thm.EP-E}}
Let $E_\mu$ be the nonsymmetric Macdonald polynomial indexed by the weight $\mu$, and let $P_\lambda$ be the symmetric Macdonald polynomial indexed by the dominant weight $\lambda$.  Let $m_{\varpi(h)}^{-1}$ be the alcove where the walk $h$ ends.  Then
$$E_\mu P_\lambda = \sum_{h} a_h(q,t) E_{\varpi(h)},$$
where the sum is over alcove walks of type determined by $\mu$ and contained in the dominant chamber, and the coefficients $a_h(q,t)$ are certain rational functions in $q$ and $t$.

\noindent{\bf Theorem~\ref{thm.PP-P}} 
Let $P_\lambda$ be the symmetric Macdonald polynomial indexed by the dominant weight $\lambda$.  Let $\ttwt(h)$ be the weight of the path $h$, and let $w_0$ be the longest element of the Weyl group.  Then
$$P_\mu P_\lambda = \sum_{h} c_h(q,t) P_{-w_0\ttwt(h)},$$
where the sum is over alcove walks of type determined by $\mu$ and contained in the dominant chamber, and the coefficients $c_h(q,t)$ are rational functions in $q$ and $t$.  

This section concludes by explaining how special cases of Theorem~\ref{thm.PP-P} relate to Macdonald's Pieri formula for symmetric Macdonald polynomials in terms of tableaux, and also explains the connection to Schwer's formula for Hall-Littlewood polynomials in terms of positively folded galleries.  

The final Section~\ref{sec.eg} contains many examples and illustrations.   A number of calculations can be made completely explicit in the case of the reduced rank-one root system of type $A_1$.

{\bf Acknowledgements.} This research was supported by the National Sciences and Engineering Research Council of Canada.  The author would like to thank A. Ram for his guidance and insight, and also J. Haglund and C. Lenart for helpful conversations.

\section{Alcoves, Weyl groups, and double affine Hecke algebras}\label{sec.daha}
\subsection{Root systems and Weyl groups}

Let $(\fh_\bbZ^*, R, \fh_\bbZ, R^\vee)$ be a reduced root datum with a pairing 
$$\langle \cdot , \cdot \rangle : \fh_\bbZ^* \times \fh_\bbZ \rightarrow \bbZ.$$  
That is, $\fh_\bbZ^*$ and $\fh_\bbZ$ are lattices of finite rank, containing the finite subsets $R$ and $R^\vee$ respectively, and there is a bijection $R \rightarrow R^\vee: \alpha \mapsto \alpha^\vee$ such that $\langle \alpha, \alpha^\vee \rangle =2$.

Let $\alpha_1, \ldots, \alpha_n \in R$ be the simple roots, and $\alpha_1^\vee, \ldots, \alpha_n^\vee \in R^\vee$ be the simple coroots.  The fundamental weights $\{\omega_1,\ldots, \omega_n\}$ are defined by $\langle \omega_i, \alpha_j^\vee \rangle = \delta_{ij}$, and the fundamental coweights $\{\omega_1^\vee, \ldots, \omega_n^\vee \}$ are defined by $\langle \alpha_i, \omega_j^\vee \rangle = \delta_{ij}$.
Let 
$$Q = \sum_{i=1}^n \bbZ \alpha_i,\quad 
Q^\vee = \sum_{i=1}^n \bbZ \alpha_i^\vee,\quad 
P = \sum_{i=1}^n \bbZ \omega_i,\quad 
P^\vee = \sum_{i=1}^n \bbZ \omega_i^\vee$$ 
be the root, coroot, weight, and coweight lattice respectively.  Then $Q \subseteq \fh_\bbZ^* \subseteq P$ and $Q^\vee \subseteq \fh_\bbZ \subseteq P^\vee$ as lattices.

Let $\fh_\bbR^* = \fh_\bbZ^* \otimes \bbR$ and $\fh_\bbR = \fh_\bbZ \otimes \bbR$.
For $\alpha \in R$, the map
$$s_\alpha : \fh_\bbR^* \rightarrow \fh_\bbR^*:x\mapsto x-\langle x, \alpha^\vee\rangle \alpha$$
acts on the lattice $\fh_\bbZ^*$ and is a reflection in the hyperplane
$$\ttH_{\alpha^\vee} = \{ x\in \fh_\bbR^* \mid \langle x, \alpha^\vee \rangle =0 \},$$
and sends $\alpha$ to $-\alpha$.
For $1\leq i \leq n$, let
$$s_i = s_{\alpha_i}.$$
The {\em Weyl group} $W_0$ is generated by $s_1, \ldots, s_n$ subject to the relations
$$s_i^2 =1, \qquad \textrm{and}\qquad
s_is_js_i \cdots = s_js_is_j \cdots \quad \textrm{$(m_{ij}$ factors each side),}$$
where $\pi/m_{ij}$ is the angle between $\ttH_{\alpha_i^\vee}$ and $\ttH_{\alpha_j^\vee}$.  See~\cite[p.69]{K01}.

\subsubsection{Double affine Weyl groups}
Let $e$ be the smallest positive integer which satisfies $\langle \fh_\bbZ^*, \fh_\bbZ \rangle \subseteq \hbox{$\frac{1}{e}$}\bbZ$.
Let 
$X=\{x^\mu \mid \mu \in \fh_\bbZ^*\}$ and $Y=\{y^{\lambda^\vee} \mid \lambda^\vee \in \fh_\bbZ\}$ be abelian groups isomorphic to $\fh_\bbZ^*$ and $\fh_\bbZ$ respectively, with multiplication 
\begin{equation}\label{eqn.xyabeliangroups}
x^\mu x^\lambda = x^{\mu+\lambda}, \quad \textrm{and} \quad
y^{\lambda^\vee}y^{\mu^\vee} = y^{\lambda^\vee+\mu^\vee}.
\end{equation}  
The {\em double affine Weyl group} $\widetilde W$ is
$$\left\{q^k x^\mu w y^{\lambda^\vee} \mid 
k \in \hbox{$\frac{1}{e}$}\bbZ, \mu \in \fh_\bbZ^*, w\in W_0, \lambda^\vee \in \fh_\bbZ\right\},$$
subject to the relations \eqref{eqn.xyabeliangroups}
and
$$wx^\mu = x^{w\mu} w, \quad
wy^{\lambda^\vee} = y^{w\lambda^\vee}w,\quad
x^\mu y^{\lambda^\vee} = q^{\langle \mu, \lambda^\vee \rangle} y^{\lambda^\vee}x^\mu, \quad
q^{1/e}\in Z(\widetilde W).$$
See~\cite[Corollary 4.6]{H06}.

The {\em extended affine Weyl groups}
\begin{equation}\label{eqn.awg1}
W = \{w y^{\lambda^\vee}\mid w\in W_0, \lambda^\vee\in \fh_\bbZ \} = W_0 \ltimes Y,
\end{equation}
\begin{equation}\label{eqn.awg2}
W^\vee = \{x^\mu w \mid \mu \in \fh_\bbZ^*, w\in W_0\} = X \rtimes W_0
\end{equation}
are subgroups of $\widetilde W$, and
$W$ acts by conjugation on $\{q^kx^\mu \mid k\in\hbox{$\frac{1}{e}$}\bbZ, \mu \in \fh_\bbZ^* \}$.
Define
\begin{equation}
x^{\mu + k\delta} = q^kx^\mu \quad\textrm{and}\quad
y^{\lambda^\vee + kd} = q^{-k}y^{\lambda^\vee}.
\end{equation}
Then $W$ acts on the lattice $\fh_\bbZ^*\oplus \bbZ\delta$, where for $w\in W$ and $\nu = \mu+k\delta\in \fh_\bbZ^*\oplus \bbZ\delta$, $w\nu$ is defined by 
\begin{equation}\label{eqn.conjaction}
x^{w\nu} = wx^\nu w^{-1} \textrm{ in }\widetilde W.
\end{equation}

For $\alpha \in R$ and $j\in \bbN$, the map
$$x^{j\alpha}s_\alpha:\fh_\bbR^* \rightarrow \fh_\bbR^* : x \mapsto s_\alpha x + j\alpha$$
acts on the lattice $\fh_\bbZ^*$ and is a reflection in the affine hyperplane
$$\ttH_{-\alpha^\vee+jd} = \{x\in \fh_\bbR^* \mid \langle x, \alpha^\vee \rangle = j\}.$$
Let $\varphi^\vee \in R^\vee$ be the maximal coroot, and $\varphi\in R$ be the maximal (short) root.  Define
$$\alpha_0 = -\varphi + \delta,\quad
\alpha_0^\vee= -\varphi^\vee + d,\quad
s_0 = y^{\varphi^\vee}s_\varphi \in W, \quad \textrm{and}\quad
s_0^\vee = x^{\varphi}s_\varphi \in W^\vee.$$
Then $s_0^\vee$ is a reflection in the hyperplane $\ttH_{\alpha_0^\vee} = \ttH_{-\varphi^\vee+d}$.

The affine Weyl group $W_a = W_0 \ltimes Q^\vee$ is generated by $s_0, s_1,\ldots, s_n$ subject to the relations
$$s_i^2 =1, \qquad \textrm{and}\qquad
s_is_js_i \cdots = s_js_is_j \cdots \quad \textrm{$(m_{ij}$ factors each side),}$$
where $\pi/m_{ij}$ is the angle between $\ttH_{\alpha_i^\vee}$ and $\ttH_{\alpha_j^\vee}$.  See~\cite[p.123]{K01}.

The extended affine Weyl group $W$ has an alternate presentation~\cite[p.132]{K01}
$$W = W_a \rtimes \Pi,$$
where $\Pi \cong \fh_\bbZ/Q^\vee$.

The dual version of the above statements for $W$ holds for $W^\vee$ as well.  That is, 
$$W^\vee = Q \rtimes W_0 = \Pi^\vee \ltimes W_a^\vee,$$ 
$\Pi^\vee \cong \fh_\bbZ^*/Q$, and $W_a^\vee$ is the group generated by $s_0^\vee, s_1, \ldots, s_n$.  For notational convenience, we sometimes write $s_i^\vee = s_i$ for $i=1,\ldots, n$.

\subsection{The alcove picture}\label{sec.alcovepicture}
See for example, a picture for the $\fsl_2$ root system in Section~\ref{sec.eg}.

Denote the positive roots and coroots by $R_+$ and $R_+^\vee$. The {\em positive affine coroots} are
$$S_+^\vee = \left\{\alpha^\vee+jd\mid  \alpha^\vee\in R_+^\vee, j \in \bbZ_{\geq0} \right\} \cup \left\{-\alpha^\vee+jd\mid  \alpha^\vee\in R_+^\vee, j \in \bbZ_{\geq1} \right\}.$$
The {\em chambers} of $W_0$ are the connected components of $\fh_\bbR^* \backslash \cup_{\alpha \in R_+}\ttH_{\alpha^\vee}$, and 
the {\em alcoves} of $W_a^\vee$ are the connected components of $\fh_\bbR^* \backslash \cup_{a\in S_+} \ttH_{a^\vee}$.  

The {\em fundamental chamber} or {\em dominant chamber} is the region
$$\ttC 
=\left \{x\in \fh_\bbR^* \mid 0 < \langle x, \alpha^\vee \rangle \textrm{ for } \alpha \in R_+\right\}
= \bigcap_{i=1}^n \{x\in \fh_\bbR^* \mid 0<\langle x, \alpha_i^\vee \rangle\},$$
whose {\em walls} (the hyperplanes which have nonempty intersection with the closure of $\ttC$) are~  $\ttH_{\alpha_1^\vee}, \ldots, \ttH_{\alpha_n^\vee}$.
The {\em fundamental alcove} is the region
$$\ttA 
= \{x\in \fh_\bbR^* \mid 0 < \langle x, \alpha^\vee \rangle <1 \textrm{ for } \alpha \in R^+\}
= \ttC \cap \{x\in \fh_\bbR^* \mid \langle x, \varphi^\vee \rangle <1\},$$
and its walls are the hyperplanes $\ttH_{\alpha_0^\vee}, \ldots, \ttH_{\alpha_n^\vee}$.

By Proposition 4-6 and Proposition 11-5~\cite{K01}, $W_0$ acts freely transitively on the chambers, and $W_a^\vee$ acts freely transitively on the alcoves so that there is a bijection
\begin{eqnarray*}
W_a^\vee &\longleftrightarrow & \{\textrm{alcoves}\} \\
w & \leftrightarrow & w\ttA.
\end{eqnarray*}
In the above correspondence, the elements of $\Pi^\vee \subseteq W^\vee = \Pi^\vee \rtimes W_a^\vee$ fix the fundamental alcove $\ttA$.
Since $|P/Q| = \det[\langle \alpha_i, \alpha_j^\vee\rangle]_{1\leq i, j \leq n}$ is finite, then $\Pi^\vee\cong \fh_\bbZ^*/Q \subseteq P/Q$ is a finite abelian group.
The extended affine Weyl group $W^\vee$ acts freely transitively on $|\Pi^\vee|$ copies (sheets) of alcoves so that there is a bijection
\begin{eqnarray*}
W^\vee &\longleftrightarrow & \{\textrm{alcoves}\} \times |\Pi^\vee|\\
w & \leftrightarrow & w\ttA,
\end{eqnarray*}
where elements in $W_a^\vee$ permute alcoves in the base sheet, and elements $\pi_j^\vee \in \Pi^\vee$ send the fundamental alcove to the copy of the fundamental alcove on the $j$th sheet.  This correspondence will be used frequently, and we will often use the shorthand $w=w\ttA$.

The {\em periodic orientation} is the orientation of the hyperplanes $\left\{\ttH_{a^\vee} \mid a^\vee \in S_+^\vee\right\}$ such that

\begin{enumerate}
\item $\ttA$ is on the positive side of $\ttH_{\alpha^\vee}$ for $\alpha^\vee \in R_+^\vee$,
\item $\ttH_{\alpha^\vee+jd}$ and $\ttH_{\alpha^\vee}$ have parallel orientations.
\end{enumerate}

The figure in Section~\ref{sec.eg} illustrates the alcove picture of the extended affine Weyl group $W^\vee$ for the $\fsl_2$ root system, showing the periodic orientation of the hyperplanes.

\subsubsection{The length function}
Given $w\in W^\vee$ with a reduced expression $w = \pi_j^\vee s_{i_1}^\vee \cdots s_{i_r}^\vee$, the set of positive coroots
\begin{equation}\label{calLw}
\calL(w) = \left\{\pi_j^\vee\alpha_{i_1}^\vee,\quad 
\pi_j^\vee s_{i_1}^\vee \alpha_{i_2}^\vee,\quad
\pi_j^\vee s_{i_1}^\vee s_{i_2}^\vee \alpha_{i_3}^\vee, \quad
\ldots,\quad 
\pi_j^\vee s_{i_1}^\vee \cdots s_{i_{r-1}}^\vee \alpha_{i_r}^\vee 
\right\}
\end{equation}
indexes the hyperplanes that separate the fundamental alcove $\ttA$ and the alcove $w\ttA$.  In Macdonald's notation~\cite[(2.2.1), (2.2.9)]{M03}, 
$$\calL(w) = \{b^\vee \in S_+^\vee \mid w^{-1}b^\vee \in S_-^\vee\} =S(w^{-1}).$$
The {\em length} of $w$ is
$$\ell(w) = |\calL(w)|,$$
the number of hyperplanes that separate $\ttA$ and $w\ttA$.  In particular, $\ell(\pi^\vee) =0 $ for all $\pi^\vee \in \Pi^\vee$.

More generally, for $v, w\in W^\vee$, the set of positive coroots
\begin{equation}\label{calLvw}
\calL(v,w) = \left(\calL(v) \cup \calL(w)\right) \backslash (\calL(v) \cap \calL(w))
\end{equation}
indexes the set of hyperplanes that separate the alcoves $v\ttA$ and $w\ttA$.  If $v\leq w$, we may write $\calL(v,w) = v\calL(1,v^{-1}w)$, where $\calL(1,w) = \calL(w)$.



\subsubsection{The group $\Pi^\vee$}
Let $w_0\in W_0$ be the unique longest element in the finite Weyl group $W_0$.
For $\mu \in \fh_\bbZ^*$, let $\mu_+$ denote the unique dominant weight in the orbit, so that $\mu_- = w_0\mu_+$ is the unique antidominant weight.
Let $v_\mu$ be the shortest element of $W_0$ such that $v_\mu\mu = \mu_-$.  Define
\begin{equation}\label{eqn.mlcr}
m_\mu = x^\mu v_\mu^{-1} \in W^\vee.
\end{equation}
By~\cite[(2.4.5)]{M03}, $m_\mu$ is the unique shortest element in the coset $x^\mu W_0$.

The {\em minuscule weights} are the fundamental weights $\omega_j$ that satisfy 
$\langle \omega_j, \alpha^\vee \rangle \leq 1$ for $\alpha^\vee \in R_+^\vee$.  In other words, these are the fundamental weights which are contained in the closure of the fundamental alcove.
Let 
$$J = \{j \mid \omega_j \in \fh_\bbZ^* \textrm{ is a minuscule weight} \} \cup \{0\}.$$
Let $\pi_0^\vee = 1$, and 
for $j \in J\backslash \{0\}$, let 
\begin{equation}
\pi_j^\vee = m_{\omega_j} = x^{\omega_j}v_{\omega_j}^{-1}.
\end{equation}  
By~\cite[(2.5.4)]{M03}, the subgroup of length zero elements in $W^\vee$ is 
$$\Pi^\vee = \left\{\pi_j^\vee \mid j\in J\right\}.$$  

\subsection{Braid groups and Hecke algebras}
\subsubsection{Double affine braid groups}
The relevant facts about braid groups from~\cite{M03} are stated here in our notation.

Let $\left\{q^k X^\mu \mid k\in\hbox{$\frac{1}{e}$}\bbZ, \mu \in \fh_\bbZ^*\right\}$
be the multiplicative group isomorphic to $\fh_\bbZ^* \oplus \bbZ\delta$, and write $X^{\mu+k\delta} = q^kX^\mu$.
By equation~\eqref{eqn.conjaction}, the conjugation action of $W = \langle s_0, \ldots, s_n \rangle \rtimes \Pi$ on 
$\left\{q^kX^\mu \mid k\in \hbox{$\frac{1}{e}$}\bbZ, \mu\in\fh_\bbZ^*\right\}$ is given by
$$w(\mu+k\delta)= w\mu + k\delta, \quad
y^{\lambda^\vee}(\mu+k\delta)=\mu -\langle \mu, \lambda^\vee \rangle \delta +k\delta,\quad
\textrm{ for } w\in W_0, \lambda^\vee \in \fh_\bbZ.$$

The {\em double affine braid group} $\widetilde B$ is generated by 
$\left\{q^k X^\mu \mid k\in\bbZ, \mu \in \fh_\bbZ^*\right\}$, $T_0, T_1, \ldots, T_n$, and $\Pi$, subject to the relations
\begin{enumerate}
\item $T_iT_jT_i \cdots = T_jT_iT_j \cdots$ \quad 
	($m_{ij}$ factors each side),
\item $T_i X^\mu = X^{\mu}T_i$ \quad 
	if $\langle \mu,\alpha_i^\vee\rangle=0$ for $0\leq i \leq n$,
\item $T_i X^\mu T_i = X^{s_i\mu}$ \quad 
	if $\langle \mu,\alpha_i^\vee\rangle=1$ for $0\leq i \leq n$,
\item $\pi T_i \pi^{-1} = T_j$ \quad 
	if  $\pi\alpha_i = \alpha_j$ for $\pi \in \Pi$,
\item $\pi X^\mu \pi^{-1} = X^{\pi\mu}$ \quad 
	for $\pi \in \Pi$,
\end{enumerate}
where $\langle \mu, \alpha_0^\vee\rangle = \langle \mu, -\varphi^\vee\rangle$.
Since $W= W_0 \rtimes Y$ fixes $\delta$, then $q^{1/e}$ is a central element of $\widetilde B$.  See~\cite[Sec 3.4]{M03}. 

For $w\in W$ with a reduced expression $w=\pi_j s_{i_1}\cdots s_{i_r}$ where $\pi\in \Pi$ and $s_{i_j}\in \langle s_0,\ldots, s_n\rangle$, define
\begin{equation}
T_w = \pi_j T_{i_1}\cdots T_{i_r}.
\end{equation}
By~\cite[(3.1.1)]{M03}, the element $T_w$ is independent of the choice of a reduced word for $w$.

Identify the reduced expression $w=\pi_j s_{i_1}\cdots s_{i_r}$
with the minimal path $p$ from $\ttA$ to $w\ttA$ via the sequence of alcoves $\pi_j\ttA, \ \pi_js_{i_1}\ttA, \ \ldots,\ \pi_js_{i_1}\cdots s_{i_r}\ttA$,  
and define
\begin{equation}\label{eqn.Y}
Y^w = \pi_j T_{i_1}^{\epsilon_1}\cdots T_{i_r}^{\epsilon_r}, \quad
\textrm{ where }
\epsilon_k = \begin{cases}
+1, & \textrm{if the $k$th step of $p$ is 
\beginpicture
\setcoordinatesystem units <1cm,1cm>                    
\setplotarea x from -0.8 to 0.8, y from -0.5 to 0.5     
\put{$\scriptstyle{-}$}[b] at -0.3 0.25 
\put{$\scriptstyle{+}$}[b] at 0.3 0.25
\plot  0 -0.4  0 0.5 /
\arrow <5pt> [.2,.67] from -0.5 0 to 0.5 0   %
\endpicture
,}\\
-1, & \textrm{if the $k$th step of $p$ is 
\beginpicture
\setcoordinatesystem units <1cm,1cm>                    
\setplotarea x from -0.8 to 0.8, y from -0.5 to 0.5     
\put{$\scriptstyle{-}$}[b] at -0.3 0.25 
\put{$\scriptstyle{+}$}[b] at 0.3 0.25
\plot  0 -0.4  0 0.5 /
\arrow <5pt> [.2,.67] from 0.5 0 to -0.5 0   %
\endpicture
,} 
\end{cases}
\end{equation}
with respect to the periodic orientation (section~\ref{sec.alcovepicture}) of the hyperplanes.

For simplicity, write $Y^{y^{\lambda^\vee}} = Y^{\lambda^\vee}$.
Then $Y^{\lambda^\vee}Y^{\mu^\vee} = Y^{\lambda^\vee+\mu^\vee}$ for  
$\lambda^\vee, \mu^\vee \in \fh_\bbZ,$
and
$$Y=\{Y^{\lambda^\vee} \mid \lambda^\vee \in \fh_\bbZ\}$$
is a multiplicative group isomorphic to $\fh_\bbZ$.
The elements $Y^{\lambda^\vee}$ satisfy the relations
$$\begin{array}{ll}
T_i^{-1} Y^{\lambda^\vee} = Y^{\lambda^\vee}T_i^{-1},
	& \hbox{if $\langle \alpha_i, \lambda^\vee\rangle=0$ 
	for $1\leq i \leq n$,}\\
T_i^{-1} Y^{\lambda^\vee} T_i^{-1} = Y^{s_i\lambda^\vee},
	& \hbox{if $\langle \alpha_i, \lambda^\vee\rangle=1$ 
	for $1\leq i \leq n$.}
\end{array}
$$
See~\cite[(3.2.4)]{M03} for details.

Define
\begin{equation}
T_0^\vee = (X^\varphi T_{s_\varphi})^{-1} \quad
\textrm{and}\quad
Y^{-\alpha_0^\vee} = qY^{\varphi^\vee}.
\end{equation}
For notational convenience, we sometimes write $T_i^\vee=T_i$ for $i=1,\ldots, n$.  

Then, identifying the reduced expression $z = s_{i_r}^\vee \cdots s_{i_1}^\vee (\pi_j^\vee)^{-1} \in W^\vee$ with the minimal path $b$ from $z\ttA$ to $\ttA$ 
via the sequence of alcoves
$$z\ttA,\quad 
z\pi_j^\vee s_{i_1}^\vee \ttA,\quad
z\pi_j^\vee s_{i_1}^\vee s_{i_2}^\vee \ttA,\quad
\ldots,\quad 
(z\pi_j^\vee s_{i_1}^\vee \cdots s_{i_r}^\vee) \ttA= \ttA,$$
let 
\begin{equation}\label{eqn.epsilonk}
X^z = (T_{i_r}^\vee)^{\epsilon_r} \cdots (T_{i_1}^\vee)^{\epsilon_1} (\pi_j^\vee)^{-1}, \quad
\hbox{where }
\epsilon_k = \begin{cases}
+1, &\hbox{if the $k$th step of $b$ is 
	\beginpicture
	\setcoordinatesystem units <1cm,1cm>                    
	\setplotarea x from -0.7 to 0.7, y from -0.5 to 0.5     
	\put{$\scriptstyle{-}$}[b] at -0.3 0.25 
	\put{$\scriptstyle{+}$}[b] at 0.3 0.25
	\plot  0 -0.4  0 0.5 /
	\arrow <5pt> [.2,.67] from -0.5 0 to 0.5 0   %
	\endpicture
,}\\
-1, &\hbox{if the $k$th step of $b$ is
		\beginpicture
	\setcoordinatesystem units <1cm,1cm>                    
	\setplotarea x from -0.7 to 0.7, y from -0.5 to 0.5     
	\put{$\scriptstyle{-}$}[b] at -0.3 0.25 
	\put{$\scriptstyle{+}$}[b] at 0.3 0.25
	\plot  0 -0.4  0 0.5 /
	\arrow <5pt> [.2,.67] from 0.5 0 to -0.5 0   %
	\endpicture
,}
\end{cases}
\end{equation} 
with respect to the periodic orientation of the hyperplanes (section~\ref{sec.alcovepicture}).  Note that for $\mu\in\fh_\bbZ^*\subseteq W^\vee$, $X^\mu=X^{x^{\mu}}$.

%

\subsubsection{Double affine Hecke algebras}
Let $\bbK$ be a field.  Fix $t_0, t_1, \ldots, t_n \in \bbK$ such that $t_i =t_j$ if $s_i$ and $s_j$ are conjugate in $W$.  For $\alpha\in R$ and $k\in \frac{1}{e}\bbZ$, define $t_{\alpha+k\delta} = t_i$  if $\alpha = w\alpha_i$ for some $w\in W$.

The {\em double affine Hecke algebra} $\widetilde H$ is the quotient of the group algebra $\bbK\widetilde B$ of the double affine braid group by the relations
$$T_i^2 = (t_i^{1/2}-t_i^{-1/2})T_i +1, \quad \textrm{for } 0 \leq i \leq n.$$

The {\em intertwining operators} (also called {\em creation operators} in~\cite[Sec 5.10]{M03}) are
\begin{align*}
\tau_i^\vee 
&= T_i^\vee + \frac{t_i^{-1/2}-t_i^{1/2}}{1-Y^{-\alpha_i^\vee}} 
= (T_i^\vee)^{-1} + \frac{(t_i^{-1/2}-t_i^{1/2})Y^{-\alpha_i^\vee}}{1-Y^{-\alpha_i^\vee}}, \qquad\textrm{for } 0\leq i \leq n,\\
\pi_j^\vee 
&= X^{\omega_j}T_{v_{\omega_j}^{-1}}, \qquad\textrm{for } j\in J.
\end{align*}
For $w\in W^\vee$ with a reduced expression $w = \pi_j^\vee s_{i_1}^\vee \cdots s_{i_r}^\vee$, define
\begin{equation}
\tau_w^\vee = \pi_j^\vee \tau_{i_1}^\vee \cdots \tau_{i_r}^\vee.
\end{equation}
By~\cite[(5.10.13)]{M03},
$\tau_w^\vee$ is independent of the choice of a reduced word for $w$.  Moreover, the intertwiners satisfy
\begin{equation}\label{eqn.intertwiningrelns}
\tau_w^\vee Y^{\lambda^\vee} = Y^{w\lambda^\vee}\tau_w^\vee, \qquad \textrm{for } w\in W^\vee, \lambda^\vee \in \fh_\bbZ.
\end{equation}
For $0\leq i \leq n$ and $w\in W^\vee$,
\begin{equation} \label{eqn.tauisquared}
\tau_i^\vee \tau_w^\vee =
\begin{cases}
\tau_{s_i^\vee w}^\vee, & \textrm{if } s_i^\vee w > w,\\
(\tau_i^\vee)^2 \tau_{s_i^\vee w}^\vee, & \textrm{if } s_i^\vee w < w,
\end{cases}
\hbox{  where  }
(\tau_i^\vee)^2 = 
\left(\frac{1-t_i^{-1}Y^{-\alpha_i^\vee}}{1-Y^{-\alpha_i^\vee}}\right) 
\!\!\!\left(\frac{1-t_iY^{-\alpha_i^\vee}}{1-Y^{-\alpha_i^\vee}}\right).
\end{equation}

\subsubsection{Polynomial representation}
The {\em affine Hecke algebra} $H$ is the subalgebra of the double affine Hecke algebra $\widetilde H$ generated by $T_0, \ldots, T_n$ and $\Pi$.  A basis for $H$ is $\{T_w Y^{\lambda^\vee} \mid w\in W_0, \lambda^\vee \in \fh_\bbZ\}$.
Let $\bbK\mathbf{1}$ be the $H$-module given by
$$\pi \mathbf{1} = \mathbf{1}, \quad
T_i\mathbf{1} = t_i^{1/2}\mathbf{1},\quad
\textrm{for $\pi\in \Pi$ and $0\leq i \leq n$.}$$
The {\em polynomial representation} of $\widetilde H$ is 
$$\bbK[X]\mathbf{1} = \mathrm{Ind}_H^{\widetilde H}\mathbf{1},$$
with basis $\{X^\mu\mathbf{1} \mid \mu\in \fh_\bbZ^*\}$.

For $0\leq i \leq n$, the operators $T_i$ act on $\bbK[X]\mathbf{1}$ by
\begin{equation}
T_i X^\mu \mathbf{1} = t_i^{1/2} X^{s_i\mu}\mathbf{1} + \left(t_i^{1/2}-t_i^{-1/2}\right) \frac{X^\mu - X^{s_i\mu}}{1-X^{\alpha_i}}\mathbf{1}.
\end{equation}

For an affine coroot $\beta^\vee + jd$, define {\em shift} and {\em height}
\begin{equation}\label{eqn.shiftheight}
q^{\shf(\beta^\vee + jd)} = q^{-j}, \qquad \textrm{and} \qquad 
t^{\hgt(\beta^\vee + jd)} = \prod_{\alpha \in R_+} 
t_\alpha^{\frac12\langle \alpha, \beta^\vee \rangle},
\end{equation}
so that
\begin{equation}
Y^{\beta^\vee+jd}\mathbf{1} = q^{-j}Y^{\beta^\vee}\mathbf{1} = q^{-j} \prod_{\alpha \in R_+} t_{\alpha}^{\frac12 \langle \alpha, \beta^\vee \rangle}\mathbf{1} = q^{\shf(\beta^\vee + jd)}t^{\hgt(\beta^\vee + jd)}\mathbf{1}. 
\end{equation}
If $t_\alpha=t$ for all $\alpha \in R_+$, then $t^{\hgt(\beta^\vee + jd)} = t^{\langle \rho, \beta^\vee \rangle}$ for $\rho = \frac12\sum_{\alpha\in R_+} \alpha^\vee$.

\subsection{Alcove walks}\label{sec.alcovewalks}
Fix a reduced factorization of $w=\pi_j^\vee s_{i_1}^\vee\cdots s_{i_r}^\vee \in W^\vee$.  
An {\em alcove walk} of type $\vec{w} = (i_1, \ldots, i_r)$ beginning at $z$ is a sequence of steps in the alcove picture, where for $i=0,\ldots, n$, a step of type $i$ is one of the following:
\begin{equation}
    \beginpicture
	\setcoordinatesystem units <0.75cm,0.75cm>         
	\setplotarea x from -1 to 1, y from -1 to 1 
    	\plot 0 -0.9 0  1  / 
	\setplotsymbol(.)
	\arrow <5pt> [.2,.67] from -0.7 0.3 to 0.7 0.3
	\put{\footnotesize $v$} at -0.7 -0.5 
	\put{\footnotesize $vs_i^\vee$} at 0.7 -0.5
	\put{$i$-crossing,} at 0 -1.4
	\endpicture
\qquad\qquad\qquad
    \beginpicture
	\setcoordinatesystem units <0.75cm,0.75cm>         
	\setplotarea x from -1 to 1, y from -1 to 1 
    	\plot 0 -0.9 0  1  /
	\setplotsymbol(.) 
	\plot -0.75 0.15 0 0.15 / \plot -0.02 0.15 -0.02 0.35 /
	\arrow <5pt> [.2,.67] from 0 0.35 to -0.75 0.35	
	\put{\footnotesize $v$} at -0.7 -0.5 
	\put{\footnotesize $vs_i^\vee$} at 0.7 -0.5
	\put{$i$-folding.} at 0 -1.4
    \endpicture
\end{equation} 
In addition, a `step' of type $\vec{\pi}^\vee$ for $\pi^\vee\in \Pi^\vee$ can be thought of as a change of sheets from the alcove $v$ to the alcove $v\pi_j^\vee$.  See the figure in section~\ref{sec.eg}. 

Let $\Gamma(\vec{w},z)$ be the set of alcove walks of type $\vec{w}$ beginning in $z$. There are $2^r$ walks in $\Gamma(\vec{w},z)$, since each step can be either a crossing or a folding.  For a walk $h\in\Gamma(\vec{w},z)$, let 
\begin{equation}\label{eqn.pkvee}
h_k^\vee \hbox{ be the positive coroot such that $\ttH_{h_k^\vee}$ separates the alcoves $v\ttA$ and $vs_{i_k}^\vee\ttA$,}
\end{equation}
in which $v\ttA$ is the alcove where $k$th step of $h$ begins.  We call $\ttH_{h_k^\vee}$ the $k$th {\em active hyperplane} of the walk $h$.
For $k=1,\ldots, r$, let 
\begin{equation}\label{eqn.bkvee}
b_k^\vee = s_{i_r}^\vee s_{i_{r-1}}^\vee \cdots s_{i_{k+1}}^\vee\alpha_{i_k}^\vee, \quad\hbox{and}\quad
t_{b_k^\vee} = t_{i_k},
\end{equation}
so that $b_r^\vee, \ldots, b_1^\vee$ is the sequence of labels of hyperplanes crossed by the walk of type $\vec{w}^{-1}=(i_r,\ldots, i_1)$ beginning in $1$.  See example~\ref{eg.1} for an illustration.

{\em Positive} and {\em negative} steps are defined with respect to the periodic orientation of the hyperplanes as follows:
\begin{equation}\label{eqn.posnegsteps}    
	\beginpicture
	\setcoordinatesystem units <0.75cm,0.75cm>         
	\setplotarea x from -1 to 1, y from -1 to 1.2 
	\plot 0 -0.7 0  1  / 
	\setplotsymbol(.)
	\arrow <5pt> [.2,.67] from -0.7 0.3 to 0.7 0.3
	\put{$\scriptstyle{+}$} at 0.4 0.9
	\put{$\scriptstyle{-}$} at -0.4 0.9
	\put{\scriptsize $v$}[t] at -0.7 -0.5 
	\put{\scriptsize $vs_i$}[t] at 0.7 -0.5
	\put{positive crossing,} at 0 -1.5
	\endpicture
\qquad
\beginpicture
	\setcoordinatesystem units <0.75cm,0.75cm>         
	\setplotarea x from -1 to 1, y from -1 to 1.2
	\plot 0 -0.7 0  1  / 
	\setplotsymbol(.)
	\arrow <5pt> [.2,.67] from 0.7 0.3 to -0.7 0.3
	\put{$\scriptstyle{+}$} at 0.4 0.9
	\put{$\scriptstyle{-}$} at -0.4 0.9
	\put{\scriptsize $vs_i$}[t] at -0.7 -0.5 
	\put{\scriptsize $v$}[t] at 0.7 -0.5
	\put{negative crossing,} at 0 -1.5
	\endpicture
\qquad
    \beginpicture
	\setcoordinatesystem units <0.75cm,0.75cm>         
	\setplotarea x from -1 to 1, y from -1 to 1.2 
	\plot 0 -0.7 0  1  / 
	\setplotsymbol(.)
	\plot 0.75 0.15 0.05 0.15 / \plot 0.02 0.15 0.02 0.35 /
	\arrow <5pt> [.2,.67] from 0.05 0.35 to 0.75 0.35
	\put{$\scriptstyle{+}$} at 0.4 0.9
	\put{$\scriptstyle{-}$} at -0.4 0.9	
	\put{\scriptsize $vs_i$}[t] at -0.7 -0.5 
	\put{\scriptsize $v$}[t] at 0.7 -0.5
	\put{positive folding,} at 0 -1.5
    \endpicture
\qquad
    \beginpicture
	\setcoordinatesystem units <0.75cm,0.75cm>         
	\setplotarea x from -1 to 1, y from -1 to 1.2 
	\plot 0 -0.7 0  1  /
	\setplotsymbol(.) 
	\plot -0.75 0.15 -0.05 0.15 / \plot -0.02 0.15 -0.02 0.35 /
	\arrow <5pt> [.2,.67] from -0.05 0.35 to -0.75 0.35	
	\put{$\scriptstyle{+}$} at 0.4 0.9
	\put{$\scriptstyle{-}$} at -0.4 0.9
	\put{\scriptsize $v$}[t] at -0.7 -0.5 
	\put{\scriptsize $vs_i$}[t] at 0.7 -0.5
	\put{negative folding.} at 0 -1.5
    \endpicture
\end{equation}
Moreover,
\begin{equation}\label{eqn.longshortsteps}
\hbox{a step is an {\em ascent} if $vs_i > v$, and it is a {\em descent} if $vs_i < v$ in the Bruhat order.}
\end{equation}


\begin{example}\label{eg.1}{\em  
See section~\ref{sec.eg} for more details on the alcove picture of type $\fsl_2$.

Let $w= (s_1 s_0^\vee)^4 = x^{-8\omega}$.  The following is an alcove walk $h \in \Gamma(\vec{w}, 1)$.
$$
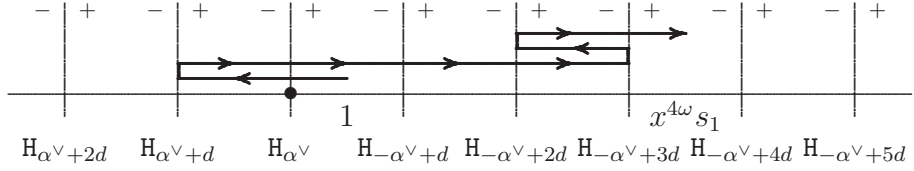
\beginpicture
\setcoordinatesystem units <1.5cm,1cm>         
\setplotarea x from -3 to 6, y from -1.5 to 1.5  
\put{$\bullet$} at 0 0 
{\small
\put{$\ttH_{-\alpha^\vee+2d}$}[t] at 2 -0.6
\put{$\ttH_{\alpha^\vee}$}[t] at 0 -0.6
\put{$\ttH_{\alpha^\vee+2d}$}[t] at -2 -0.6 
\put{$\ttH_{-\alpha^\vee+d}$}[t] at 1 -0.6
\put{$\ttH_{-\alpha^\vee+3d}$}[t] at 3 -0.6
\put{$\ttH_{\alpha^\vee+d}$}[t] at -1 -0.6
\put{$\ttH_{-\alpha^\vee+4d}$}[t] at 4 -0.6
\put{$\ttH_{-\alpha^\vee+5d}$}[t] at 5 -0.6
}\plot -2.5 0  5.5 0 /  
\plot  -2 1.2  -2 -0.3 /
\plot  -1 1.2  -1 -0.3 / \plot  0 1.2  0 -0.3 / \plot  1 1.2  1 -0.3
/ \plot  2 1.2  2 -0.3 / \plot  3 1.2  3 -0.3 / \plot  4 1.2  4 -0.3
/ \plot 5 1.2 5 -0.3 /
\put{$\scriptstyle{+}$} at -1.8 1.1 \put{$\scriptstyle{-}$} at -2.2 1.1
\put{$\scriptstyle{+}$} at -0.8 1.1 \put{$\scriptstyle{-}$} at -1.2 1.1
\put{$\scriptstyle{+}$} at 0.2 1.1 \put{$\scriptstyle{-}$} at -0.2 1.1
\put{$\scriptstyle{+}$} at 1.2 1.1 \put{$\scriptstyle{-}$} at 0.8 1.1
\put{$\scriptstyle{+}$} at 2.2 1.1 \put{$\scriptstyle{-}$} at 1.8 1.1
\put{$\scriptstyle{+}$} at 3.2 1.1 \put{$\scriptstyle{-}$} at 2.8 1.1
\put{$\scriptstyle{+}$} at 4.2 1.1 \put{$\scriptstyle{-}$} at 3.8 1.1
\put{$\scriptstyle{+}$} at 5.2 1.1 \put{$\scriptstyle{-}$} at 4.8 1.1
\put{$1$} at 0.5 -0.3
\put{$x^{4\omega}s_1$} at 3.5 -0.3
\setplotsymbol(.)
\arrow <6pt> [.2,.67] from 0.5 0.2 to -0.5 0.2
\plot -0.5 0.2 -1 0.2 / \plot -0.99 0.2 -0.99 0.4 /
\arrow <6pt> [.2,.67] from -1 0.4 to -0.5 0.4
\arrow <6pt> [.2,.67] from -0.5 0.4 to 0.5 0.4
\arrow <6pt> [.2,.67] from 0.5 0.4 to 1.5 0.4
\arrow <6pt> [.2,.67] from 1.5 0.4 to 2.5 0.4
\plot 2.5 0.4 3 0.4 / \plot 2.99 0.4 2.99 0.6 /
\arrow <6pt> [.2,.67] from 3 0.6 to 2.5 0.6
\plot 2.5 0.6 2 0.6 / \plot 2.01 0.6 2.01 0.8 /
\arrow <6pt> [.2,.67] from 2 0.8 to 2.5 0.8
\arrow <6pt> [.2,.67] from 2.5 0.8 to 3.5 0.8
\endpicture
$$
This walk of length eight has type $\vec{w}=(1,0,1,0,1,0,1,0)$.  The coroots $h_k^\vee$ are
\begin{align*}
h_1^\vee, \ldots, h_8^\vee 
&= 
\alpha^\vee,\ \alpha^\vee+d,\ \alpha^\vee,\ -\alpha^\vee+d,\ -\alpha^\vee+2d,\ -\alpha^\vee+3d,\ -\alpha^\vee+2d,\ -\alpha^\vee +3d.
\end{align*}
For $k=1,\ldots, 8$, the coroots $b_k^\vee= s_{i_8}^\vee \cdots s_{i_{k+1}}^\vee \alpha_{i_k}^\vee = -\alpha^\vee+(9-k)d$ are the labels of hyperplanes crossed by the walk of type $\vec{w}^{-1}$ starting in $1$.
$$
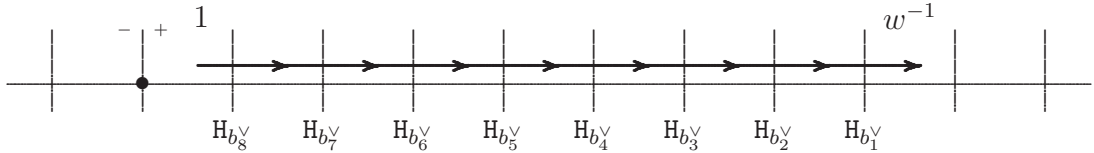
\beginpicture
\setcoordinatesystem units <1.2cm,1.2cm>         
\setplotarea x from -2 to 11, y from -1 to 0.8  
\put{$\bullet$} at 0 0 
\plot -1.5 0  10.5 0 / 
\plot  -1 0.6  -1 -0.3 / \plot  0 0.6  0 -0.3 / \plot  1 0.6  1 -0.3
/ \plot  2 0.6  2 -0.3 / \plot  3 0.6  3 -0.3 / \plot  4 0.6  4 -0.3
/ \plot  5 0.6  5 -0.3 / \plot  6 0.6  6 -0.3 / \plot  7 0.6  7 -0.3
/ \plot  8 0.6  8 -0.3 / \plot  9 0.6  9 -0.3 / \plot  10 0.6  10 -0.3
/
\put{\tiny$\scriptstyle{+}$} at 0.2 0.6 
\put{\tiny$\scriptstyle{-}$} at -0.2 0.6
\put{$w^{-1}$} at 8.5 0.75
\put{$1$} at 0.65 0.75
{\small
\put{$\ttH_{b_8^\vee}$}[t] at 1 -0.4 
\put{$\ttH_{b_7^\vee}$}[t] at 2 -0.4 
\put{$\ttH_{b_6^\vee}$}[t] at 3 -0.4
\put{$\ttH_{b_5^\vee}$}[t] at 4 -0.4 
\put{$\ttH_{b_4^\vee}$}[t] at 5 -0.4
\put{$\ttH_{b_3^\vee}$}[t] at 6 -0.4
\put{$\ttH_{b_2^\vee}$}[t] at 7 -0.4
\put{$\ttH_{b_1^\vee}$}[t] at 8 -0.4
}
\setplotsymbol(.)
\arrow <6pt> [.2,.67] from 0.5 0.2 to 1.5 0.2
\arrow <6pt> [.2,.67] from 1.5 0.2 to 2.5 0.2
\arrow <6pt> [.2,.67] from 2.5 0.2 to 3.5 0.2
\arrow <6pt> [.2,.67] from 3.5 0.2 to 4.5 0.2
\arrow <6pt> [.2,.67] from 4.5 0.2 to 5.5 0.2
\arrow <6pt> [.2,.67] from 5.5 0.2 to 6.5 0.2
\arrow <6pt> [.2,.67] from 6.5 0.2 to 7.5 0.2
\arrow <6pt> [.2,.67] from 7.5 0.2 to 8.5 0.2
\endpicture
$$
\hfill$\diamond$
}\end{example}

\section{Macdonald polynomials}
\label{sec.poly}

\subsection{Nonsymmetric Macdonald polynomials}\label{sec.nonsymmMac}
Given a weight $\mu \in \fh_\bbZ^*$, let $m_\mu \in W^\vee$ be the shortest element in the coset $x^\mu W_0$.  See~\eqref{eqn.mlcr}.

For $w\in W^\vee$, define $t_w = t_{i_1}\cdots t_{i_r}$, if $w=s_{i_1}^\vee \cdots s_{i_r}^\vee \in W^\vee$ is a reduced expression.

\begin{defn}{\em Let $m_\mu = \pi_j^\vee s_{i_1}^\vee\cdots s_{i_r}^\vee$ be a reduced expression.
The {\em nonsymmetric Macdonald polynomial} $E_\mu\in \bbK[X]$ indexed by $\mu \in \fh_\bbZ^*$ is defined by
\begin{equation}\label{eqn.E}
E_\mu \mathbf{1} = \tau_{m_\mu}^\vee \mathbf{1} 
= \pi_j^\vee \tau_{i_1}^\vee \cdots \tau_{i_r}^\vee\mathbf{1}.
\end{equation}
}\end{defn}

\begin{remark}{\em
In equation~\eqref{eqn.E}, the coefficient of $X^\mu$ in $E_\mu$ is $t_{v_\mu^{-1}}^{1/2}$, where $x^\mu = m_\mu v_\mu$~\eqref{eqn.mlcr}. 
In the literature, $E_\mu$ is often normalized so that the coefficient of $X^\mu$ in $E_\mu$ is $1$.  

By~\eqref{eqn.intertwiningrelns}, the nonsymmetric Macdonald polynomials are eigenfunctions for the operators $\{Y^{\lambda^\vee} \mid\lambda^\vee \in \fh_\bbZ\}$, and the set $\{E_\mu\mathbf{1} \mid \mu\in \fh_\bbZ^*\}$ is a basis for the polynomial representation $\bbK[X]\mathbf{1}$.  
}\end{remark}

The following result gives an expansion of a product of monomials $X^\nu$ and intertwining operators in terms of monomials in $\widetilde H$.  This leads to an expression for $E_\mu$ in the monomial basis.

For a walk $p\in \Gamma(\vec{w},z)$ of type $\vec{w}$ beginning in $z$, let
\begin{equation}
\begin{array}{rl}
\phi(p)\!\!\!&=  \{k\ |\ \hbox{the $k$th step of $p$ is a fold}\},\\
\phi_-(p)\!\!\!&=  \{k\ |\ \hbox{the $k$th step of $p$ is a negative fold}\},
\end{array}
\end{equation}
\begin{equation}\label{endpoint}
\begin{array}{l}
\ttb(p) \in W^\vee \hbox{ be the alcove where $p$ begins,}\\
\tte(p) \in W^\vee \hbox{ be the alcove where $p$ ends.}
\end{array}
\end{equation}
The {\em weight} $\ttwt(p) \in \fh_\bbZ^*$ and {\em final direction} $\ttd(p)\in W_0$ of a walk $p$ is defined by 
\begin{equation}
X^{\tte(p)} = X^{\ttwt(p)}T_{\ttd(p)}.
\end{equation}

\begin{theorem}\label{thm.Xtau-X}{\ }
\begin{enumerate}
\item[(a)] \cite[Theorem 2.2]{RY10}
Let $z, w\in W^\vee$, and fix a reduced expression $w = \pi_j^\vee s_{i_1}^\vee\cdots s_{i_r}^\vee$.  Then
$$ X^z \tau_w^\vee = \sum_{p \in \Gamma(\vec{w},z)} X^{\tte(p)} 
\left(\prod_{k\in \phi(p)}\frac{t_{b_k^\vee}^{-1/2}-t_{b_k^\vee}^{1/2}}{1-Y^{-b_k^\vee}}\right)
\left( \prod_{k\in \phi_-(p)} Y^{-b_k^\vee} \right)
,$$
where 
$b_k^\vee = s_{i_r}^\vee\cdots s_{i_{k+1}}^\vee \alpha_{i_k}^\vee$ and $t_{b_k^\vee} = t_{i_k}$, see~\eqref{eqn.bkvee}.\\

\item[(b)] \cite[Theorem 3.1]{RY10}
Let $\mu \in \fh_\bbZ^*$ and fix a reduced expression $m_\mu = \pi_j^\vee s_{i_1}^\vee \cdots s_{i_r}^\vee\in W^\vee$ for the minimal length representative of the coset $x^\mu W_0$.  The nonsymmetric Macdonald polynomial
$$E_\mu
= \sum_{p\in \Gamma(\vec{m}_\mu)} 
\left(\prod_{k\in \phi(p)}\frac{t_{b_k^\vee}^{-1/2}-t_{b_k^\vee}^{1/2}}{1-q^{\shf(-b_k^\vee)}t^{\hgt(-b_k^\vee)}}\right)
\left( \prod_{k\in \phi_-(p)} q^{\shf(-b_k^\vee)}t^{\hgt(-b_k^\vee)}\right)
t_{\ttd(p)}^{1/2}  X^{\ttwt(p)},$$
where the sum is over the set $\Gamma(\vec{m}_\mu) = \Gamma(\vec{m}_\mu, 1)$ of alcove walks of type $\vec{m}_\mu=  (\pi_j^\vee, i_1,\ldots , i_r)$ beginning in the fundamental alcove.
\end{enumerate}

\end{theorem}
\proof {\ }
\begin{enumerate}
\item[(a)] This is proved by using the definition of $\tau_i^\vee$ and the commutation relations~\eqref{eqn.intertwiningrelns}.\\

\item[(b)] Since $E_\mu\mathbf{1} = \tau_{m_\mu}^\vee\mathbf{1}$, then
the result follows from (a) by setting $w=m_\mu$, $z=1$, using
$$Y^{-b_k^\vee}\mathbf{1} = q^{\shf(-b_k^\vee)}t^{-\hgt(-b_k^\vee)}\mathbf{1}, \qquad
X^{\tte(p)}\mathbf{1} = X^{\ttwt(p)}T_{\ttd(p)} \mathbf{1}
= t_{\ttd(p)}^{1/2} X^{\ttwt(p)}\mathbf{1}.$$
\end{enumerate}
\qed

The next result gives an expansion of a product of monomials and intertwining operators in terms of intertwining operators in $\tilde H$.

Recall the definition of $X^z$ from Equation~\eqref{eqn.epsilonk}: identify a reduced expression $z = s_{i_r}^\vee \cdots s_{i_1}^\vee (\pi_j^\vee)^{-1} \in W^\vee$ with the minimal walk $b$ from $z\ttA$ to $\ttA$ of type $\vec{z}^{-1}$ via the sequence of alcoves
$$z\ttA,\quad
z\pi_j^\vee s_{i_1}^\vee \ttA,\quad
z\pi_j^\vee s_{i_1}^\vee s_{i_2}^\vee \ttA,\quad
\ldots,\quad
(z\pi_j^\vee s_{i_1}^\vee \cdots s_{i_r}^\vee) \ttA= \ttA, \qquad
\hbox{so}$$
$$X^z = (T_{i_r}^\vee)^{\epsilon_r} \cdots (T_{i_1}^\vee)^{\epsilon_1} (\pi_j^\vee)^{-1}, \quad
\hbox{where }
\epsilon_k = \begin{cases}
+1, &\hbox{if the $k$th step of $b$ is 
	\beginpicture
	\setcoordinatesystem units <1cm,1cm>                    
	\setplotarea x from -0.7 to 0.7, y from -0.5 to 0.5     
	\put{$\scriptstyle{-}$}[b] at -0.3 0.25 
	\put{$\scriptstyle{+}$}[b] at 0.3 0.25
	\plot  0 -0.4  0 0.5 /
	\arrow <5pt> [.2,.67] from -0.5 0 to 0.5 0   %
	\endpicture
,}\\
-1, &\hbox{if the $k$th step of $b$ is
		\beginpicture
	\setcoordinatesystem units <1cm,1cm>                    
	\setplotarea x from -0.7 to 0.7, y from -0.5 to 0.5     
	\put{$\scriptstyle{-}$}[b] at -0.3 0.25 
	\put{$\scriptstyle{+}$}[b] at 0.3 0.25
	\plot  0 -0.4  0 0.5 /
	\arrow <5pt> [.2,.67] from 0.5 0 to -0.5 0   %
	\endpicture
.}
\end{cases}
$$

Given a walk $h\in \Gamma(\vec{z}^{-1},w^{-1})$, let
\begin{equation}
\begin{array}{rl}
\xi_{\rmdes}(h)\!\!\!&=  \{k\mid \hbox{the $k$th step of $h$ is a descending crossing}\},\\
\phi_{\rmasc}(h)\!\!\!&=  \{k\mid \hbox{the $k$th step of $h$ is an ascending fold}\},
\end{array}
\end{equation}
\begin{equation}\label{eqn.psip}
\psi(h) = \left\{k\ \bigg|
\begin{array}{c}
\hbox{the $k$th step of $h$ is an ascending fold and } \epsilon_k = -1 \\
\hbox{or the $k$th step of $h$ is a descending fold and } \epsilon_k = +1
\end{array} \right\}.
\end{equation}

\begin{theorem}\label{thm.X-tau}
Let $z, w\in W^\vee$, and fix a reduced expression $z = s_{i_r}^\vee \cdots s_{i_1}^\vee (\pi_j^\vee)^{-1}$.  
Then
$$ X^z \tau_w^\vee 
	=\sum_{h\in \Gamma(\vec{z}^{-1},w^{-1})} 
	(-1)^{|\phi_{\rmasc}(h)|}\tau_{\tte(h)^{-1}}^\vee g_h(Y) n_h(Y),$$
where the sum is over all alcove walks of
type $\vec{z}^{-1} = (i_1,\ldots, i_r)$ beginning in $w^{-1}(\pi_j^\vee)^{-1}$, 
$$g_h(Y) 
	= 
	\left(\prod_{k\in\phi(h)} 
	\frac{t_{i_k}^{-1/2}-t_{i_k}^{1/2}}{1-Y^{-h_k^\vee}}\right)\!\!\!
	\left(\prod_{k\in\psi(h)} Y^{-h_k^\vee}\right),\quad
n_h(Y) 
	= \prod_{k\in \xi_\rmdes(h)} 
	\frac{1-t_{i_k}^{-1} Y^{-h_k^\vee}}{1-Y^{-h_k^\vee}}
	\frac{1-t_{i_k} Y^{-h_k^\vee}}{1-Y^{-h_k^\vee}},
$$
and $h_k^\vee$ indexes the $k$th active hyperplane of $h$, see~\eqref{eqn.pkvee}.
\end{theorem}

\proof The proof is by induction on the length of $z$.  The base case is the formulas
$$(T_{i}^\vee)^{\epsilon_i} 
= \tau_i^\vee - \frac{(t_i^{-1/2}-t_i^{1/2})(Y^{-\alpha_i^\vee})^{\frac12(1-\epsilon_i)}}{1-Y^{-\alpha_i^\vee}} 
= \tau_i^\vee + \frac{(t_i^{-1/2}-t_i^{1/2})(Y^{\alpha_i^\vee})^{\frac12(1+\epsilon_i)}}{1-Y^{\alpha_i^\vee}},$$
where $\epsilon_i \in \{\pm1\}$.
For the induction step, 
let $h \in \Gamma(\vec{z}^{-1}, w^{-1})$,
$$g_h(Y) 
	= 
	\prod_{k\in\phi(h)} 
	\frac{t_{i_k}^{-1/2}-t_{i_k}^{1/2}}{1-Y^{-h_k^\vee}} 
	\prod_{k\in\psi(h)} Y^{-h_k^\vee},\qquad
n_h(Y) 
	= \prod_{k\in \xi_\rmdes(h)} 
	\frac{1-t_{i_k}^{-1} Y^{-h_k^\vee}}{1-Y^{-h_k^\vee}}
	\frac{1-t_{i_k} Y^{-h_k^\vee}}{1-Y^{-h_k^\vee}},
$$
and let $h_1, h_2 \in \Gamma(\vec{v}^{-1}s_{i}^\vee, w^{-1})$ be the two extensions of $h$ by a crossing and a folding of type $i$, respectively.


By induction, a term in $(T_i^\vee)^{\epsilon_i}X^z\tau_w^\vee$ is
\begin{align*}
(&-1)^{|\phi_a(h)|}
	\left(T_i^\vee\right)^{\epsilon_i} \tau_{\tte(h)^{-1}}^\vee 
	g_h(Y)n_h(Y)\\
&= \begin{cases}
(-1)^{|\phi_a(h)|} 
	\displaystyle\left(\tau_{i}^\vee - \frac{(t_{i}^{-1/2}-t_{i}^{1/2}) 
	\left(Y^{-\alpha_{i}^\vee} \right)^{\frac12(1-\epsilon_i)}}
	{1-Y^{-\alpha_{i}^\vee}} \right) 
	\tau_{\tte(h)^{-1}}^\vee g_h(Y) n_h(Y),
		&\hbox{if } \tte(h)s_i^\vee > \tte(h),\\ 
(-1)^{|\phi_a(h)|}
	\displaystyle\left(\tau_{i}^\vee + \frac{(t_{i}^{-1/2}-t_{i}^{1/2}) 
	\left(Y^{\alpha_{i}^\vee} \right)^{\frac12(1+\epsilon_i)}}
	{1-Y^{\alpha_{i}^\vee}} \right) 
	\tau_{\tte(h)^{-1}}^\vee g_h(Y) n_h(Y), 
		&\hbox{if } \tte(h)s_i^\vee < \tte(h),
\end{cases}\\
&=
(-1)^{|\phi_a(h_1)|} \tau_{\tte(h_1)^{-1}}^\vee g_{h_1}(Y) n_{h_1}(Y)
+ (-1)^{|\phi_a(h_2)|} \tau_{\tte(h_2)^{-1}}^\vee g_{h_2}(Y) n_{h_2}(Y),
\end{align*}
since if $\tte(h)s_i^\vee < \tte(h)$ (the last step is descending), then 
$$\tau_i^\vee \tau_{\tte(h)^{-1}}^\vee
= (\tau_i^\vee)^2 \tau_{(s_i^\vee\tte(h_1))^{-1}}^\vee
= \tau_{\tte(h_1)^{-1}}^\vee
\left(\frac{1-t_{i}^{-1}Y^{-\tte(h_1)\alpha_{i}^\vee}}
	{1-Y^{-\tte(h_1)\alpha_{i}^\vee}} \right)
	\left(\frac{1-t_{i}Y^{-\tte(h_1)\alpha_{i}^\vee}}
	{1-Y^{-\tte(h_1)\alpha_{i}^\vee}} \right)
,$$
and the hyperplane crossed by the last step of $h_1$ is indexed by
$$h_{r+1}^\vee = -\tte(h)\alpha_i^\vee = -\tte(h_1)s_i^\vee \alpha_i^\vee = \tte(h_1)\alpha_i^\vee.$$
\qed

The bijection between left and right $W_0$-cosets of $W^\vee$ gives a bijection between 
minimal length left coset representatives and alcoves in the dominant chamber (minimal length right coset representatives) via taking inverses:
\begin{eqnarray*}
W^\vee/W_0 & \longleftrightarrow & W_0\backslash W^\vee\\
m_\mu & \leftrightarrow & m_\mu^{-1}.
\end{eqnarray*}
If $h$ is a walk whose endpoint $\tte(h)$ is in the dominant chamber, define $\varpi(h)\in\fh_\bbZ^*$ by
\begin{equation}\label{eqn.varweight}
\tte(h)^{-1}=m_{\varpi(h)}.
\end{equation}
Let
\begin{equation}
\begin{array}{rl}
\phi_{\mathrm{aff}}(h)\!\!\! &= \{ k \mid \hbox{the $k$th step of $h$ is a fold touching an affine hyperplane}\},\\
\phi_{\mathrm{o}}(h)\!\!\! &= \{ k \mid \hbox{the $k$th step of $h$ is a fold touching a hyperplane containing $0$}\},
\end{array}
\end{equation}
\begin{equation}\label{eqn.psiCaff}
\psi_{\mathrm{aff}}(h) = \left\{k\in \phi_\mathrm{aff}(h) \bigg|
\begin{array}{c}
\hbox{the $k$th step of $h$ is a negative fold and } \epsilon_k = -1 \\
\hbox{or the $k$th step of $h$ is a positive fold and } \epsilon_k = +1
\end{array} \right\}.
\end{equation}
The closure of the dominant chamber is $\overline{\ttC}= \{x\in \fh_\bbR^* \mid 0 \leq \langle x, \alpha^\vee \rangle \hbox{ for } \alpha \in R_+\}$.

\begin{corollary}\label{cor.X-tau}
Let $\mu\in \fh_\bbZ^*$, and fix a reduced expression 
$m_\mu = s_{i_r}^\vee \cdots s_{i_1}^\vee (\pi_j^\vee)^{-1}$ for the minimal length element $m_\mu$ in the coset $x^\mu W_0$.  
The monomial $X^\mu$ as a linear combination of nonsymmetric Macdonald polynomials is
$$ X^\mu 
= t_{v_\mu^{-1}}^{-1/2}\sum_{h\in \Gamma^{\overline\ttC}(\vec{m}_\mu^{-1})}
(-1)^{\phi_-(h)}g_h n_h E_{\varpi(h)},$$
where the sum is over all alcove walks of type $\vec{m_\mu}^{-1}=(\pi_j^\vee, i_1,\ldots, i_r)$ beginning in the fundamental alcove and contained in $\overline{\ttC}$, 
\begin{align*}
g_h &=\left(\prod_{k\in\phi_{\mathrm{o}}(h)} t_{i_k}^{\epsilon_k/2}\right)
	\left(\prod_{k\in\phi_\mathrm{aff}(h)} 
	\frac{t_{i_k}^{-1/2}-t_{i_k}^{1/2}}
		{1-q^{\shf(-h_k^\vee)}t^{\hgt(-h_k^\vee)}} \right)	
	\left(\prod_{k\in \psi_\mathrm{aff}(h)} q^{\shf(-h_k^\vee)}t^{\hgt(-h_k^\vee)}\right),\\
n_h &= \prod_{k\in \xi_-(h)} 
	\frac{1-q^{\shf(-h_k^\vee)}t^{\hgt(-h_k^\vee)}t_{i_k}^{-1}}
	{1-q^{\shf(-h_k^\vee)}t^{\hgt(-h_k^\vee)}}
	\frac{1-q^{\shf(-h_k^\vee)}t^{\hgt(-h_k^\vee)}t_{i_k}}
	{1-q^{\shf(-h_k^\vee)}t^{\hgt(-h_k^\vee)}},
\end{align*}
$h_k^\vee$ indexes the $k$th active hyperplane of $h$, and 
$\epsilon_1,\ldots, \epsilon_r$ are as defined in~\eqref{eqn.epsilonk}.
\end{corollary}
\proof
Since $X^\mu\mathbf{1} 
= X^{m_\mu} T_{v_\mu^{-1}}^{-1}\mathbf{1} 
= t_{v_\mu^{-1}}^{-1/2} X^{m_\mu}\mathbf{1}$, then set $z = m_\mu$, $w=1$ in Theorem~\ref{thm.X-tau} and use 
$$Y^{-h_k^\vee}\mathbf{1}= q^{\shf(-h_k^\vee)}t^{\hgt(-h_k^\vee)}\mathbf{1}, \quad \hbox{and} \quad \tau_{\tte(h)^{-1}}^\vee\mathbf{1} = \tau_{m_{\varpi(h)}}^\vee\mathbf{1} = E_{\varpi(h)}\mathbf{1}. 
$$ 

The main difference in applying the formula for $X^\mu$ in Theorem~\ref{thm.X-tau} to $\mathbf{1}$ is that if a walk $h\in \Gamma(\vec{m}_\mu^{-1})$ is not contained in the dominant chamber, then it has zero contribution to the sum.
Identify the expression $\tau_{\tte(h)^{-1}}^\vee n_h(Y)$ with the product of interwining operators which correspond to the crossing steps of $h$.
That is, 
$$\tau_{\tte(h)^{-1}}^\vee n_h(Y) = \tau_{i_{c_l}}^\vee \cdots \tau_{i_{c_1}}^\vee(\pi_j^\vee)^{-1},$$  where the $c_d$th step of $h$ is a crossing for $d=1,\ldots, l$.
If $h$ leaves the dominant chamber at the $k$th crossing, then $\left(s_{i_{c_k}}^\vee \cdots s_{i_{c_1}}^\vee (\pi_j^\vee)^{-1}\right)^{-1}$ is not a minimal length coset representative, and
\begin{align*}
\tau_{\tte(h)^{-1}}^\vee n_h(Y)g_h(Y)\mathbf{1} 
&= 
	g_h \ 
	\tau_{i_{c_l}}^\vee \cdots \tau_{i_{c_{k+1}}}^\vee
	\left(\tau_{i_{c_k}}^\vee \cdots \tau_{i_{c_1}}^\vee 
	(\pi_j^\vee)^{-1}\mathbf{1}\right) = 0.
\end{align*}

Since all walks which are now under consideration are contained in $\overline\ttC$, a few simplifications can be made.  The steps of a walk $h$ contained in $\overline\ttC$ may be analyzed according to whether the active hyperplane is a wall of the dominant chamber (ie.\! is one of $\ttH_{\alpha_1^\vee},\ldots, \ttH_{\alpha_n^\vee}$), or is an affine hyperplane.  There are three possibilities:
\begin{enumerate}
\item[(a)] $k\in \xi_\mathrm{aff}(h)$: the $k$th step of $h$ crosses an affine hyperplane.  

In this case, ascending crossings are equivalent to positive crossings, and descending crossings are equivalent to negative crossings.  So $\xi_\rmdes(h) = \xi_-(h)$.

\item[(b)] $k\in \phi_\mathrm{o}(h)$: the $k$th step of $h$ is a fold touching a wall of $\ttC$.

In this case, $Y^{-h_k^\vee}\mathbf{1} = q^0t^{-\langle \alpha_{i_k}^\vee,\rho \rangle} = t_{i_k}^{-1}$, and
the fold is necessarily ascending and positive.  In particular, 
$$
-\frac{(t_{i_k}^{-1/2}-t_{i_k}^{1/2})(Y^{-h_k^\vee})^{\frac12(1-\epsilon_k)}}
{1-Y^{-h_k^\vee}}\mathbf{1} = t_{i_k}^{\epsilon_k/2}\mathbf{1}.
$$

\item[(c)] $k\in \phi_\mathrm{aff}(h)$: the $k$th step of $h$ is a fold touching an affine hyperplane.

In this case, ascending folds are equivalent to negative folds, and descending folds are equivalent to positive folds.  
\end{enumerate}

Therefore, (b) and (c) give 
\begin{align*} 
(-1)&^{|\phi_\rmasc(h)|}
	\left(\prod_{k\in\phi(h)} 
	\frac{t_{i_k}^{-1/2}-t_{i_k}^{1/2}}{1-Y^{-h_k^\vee}} \right)
	\left(\prod_{k\in\psi(h)} Y^{-h_k^\vee}\right) \mathbf{1} \\
&=  \left(\prod_{k\in\phi_{\mathrm{o}}(h)} t_{i_k}^{\epsilon_k/2}\right)
\left((-1)^{\phi_-(h)}
\prod_{k\in\phi_{\mathrm{aff}}(h)}
	\frac{t_{i_k}^{-1/2}-t_{i_k}^{1/2}}
		{1-q^{\shf(-h_k^\vee)}t^{\hgt(-h_k^\vee)}} 	
	\prod_{k\in \psi_{\mathrm{aff}}(h)} q^{\shf(-h_k^\vee)}t^{\hgt(-h_k^\vee)}\right)\mathbf{1}.
\end{align*}
\qed

\begin{remark}{\label{rem.b}\ }
{\em 
\begin{enumerate}
\item[(a)]
If the walk $h$ is contained in the dominant chamber and the $k$th step of $h$ is a fold against an affine hyperplane indexed by $h_k^\vee = -\beta^\vee + jd$, then
$q^{\shf(-h_k^\vee)}t^{\hgt(-h_k^\vee)}\mathbf{1} 
= q^j t^{\langle \beta^\vee, \rho \rangle}\mathbf{1},$ where $j \in \bbN$ and $\langle \beta^\vee, \rho\rangle \in \bbN$, since $\beta^\vee \in R_+^\vee$.
\item[(b)]
Moreover, if $\mu$ is a dominant weight, then $\epsilon_k = -1$ for all $k=1,\ldots, r$ in equation~\eqref{eqn.epsilonk}, and so
$$\psi_\mathrm{aff}(h) = \phi_-(h).$$
\end{enumerate}
}\end{remark}

\subsection{Symmetric Macdonald polynomials}\label{sec.symmMac}
The {\em dominant weights} are 
$$(\fh_\bbZ^*)^+ = \{\mu \in \fh_\bbZ^* \mid \langle \mu, \alpha_i^\vee \rangle \geq 0 \hbox{ for } i=1,\ldots, n \}.$$
Given a dominant weight $\mu \in (\fh_\bbZ^*)^+$, let $m_\mu \in W^\vee$ be the shortest element in the coset $x^\mu W_0$ (see equation~\eqref{eqn.mlcr}).
Let $W_\mu \subseteq W_0$ be the stabilizer of $\mu$, $W^\mu$ be the minimal length representatives of $W_0/W_\mu$, and $w_\mu\in W_\mu$, $v_\mu \in W^\mu$ be the longest element in their respective set (see~\cite{BB05}).
Also let $W_\mu(t)= \sum_{u\in W_\mu} t_u$ be the {\em Poincar\'{e} polynomial} of $W_\mu$.

\begin{lemma}\label{lem.Poincare}  
Let $\mu \in (\fh_\bbZ^*)^+$, and let $W_\mu(t)$ be the Poincar\'e polynomial of the stabilizer of $\mu$.  Then
$$\sum_{u\in W_\mu} 
\left( \prod_{\alpha^\vee \in \calL(1,u)} t_{\alpha^\vee}^{\frac12}
	\frac{1-t^{\hgt(-\alpha^\vee)}t_{\alpha^\vee}^{-1}}	
	{1-t^{\hgt(-\alpha^\vee)}} \right)
\left( \prod_{\alpha^\vee \in \calL(u,w_\mu)}
t_{\alpha^\vee}^{-\frac12}
\frac{1-t^{\hgt(-\alpha^\vee)}t_{\alpha^\vee}}
	{1-t^{\hgt(-\alpha^\vee)}} \right)
= t_{w_\mu}^{-\frac12} W_\mu(t).$$
\end{lemma}
\proof  
For $u\in W_\mu$, since 
$$\calL(1,u) = \left\{u\beta^\vee \in R_+^\vee: \beta^\vee \in R_-^\vee\right\}, 
\quad\hbox{ and }\quad 
\calL(u,w_\mu) = \left\{u\beta^\vee \in R_-^\vee: \beta^\vee \in R_-^\vee\right\},$$
then
\begin{align*}
\sum_{u\in W_\mu} & 
	\left( \prod_{\alpha^\vee \in \calL(1,u)} t_{\alpha^\vee}^{\frac12}
	\frac{1-t^{\hgt(-\alpha^\vee)}t_{\alpha^\vee}^{-1}}	
	{1-t^{\hgt(-\alpha^\vee)}} \right)
\left( \prod_{\alpha^\vee \in \calL(u,w_\mu)}
	t_{\alpha^\vee}^{-\frac12}
	\frac{1-t^{\hgt(-\alpha^\vee)}t_{\alpha^\vee}}
	{1-t^{\hgt(-\alpha^\vee)}} \right)\mathbf{1}\\
&= \sum_{u\in W_\mu}
\left( \prod_{\beta^\vee \in R_-^\vee: \atop
		u\beta^\vee \in R_+^\vee} t_{-u\beta^\vee}^{\frac12}
	\frac{1-Y^{-u\beta^\vee}t_{-u\beta^\vee}^{-1}}	
	{1-Y^{-u\beta^\vee}} \right)
\left( \prod_{\beta^\vee \in R_-^\vee: \atop
		u\beta^\vee \in R_-^\vee} t_{u\beta^\vee}^{-\frac12}
	\frac{1-Y^{u\beta^\vee}t_{u\beta^\vee}}
	{1-Y^{u\beta^\vee}} \right)\mathbf{1}\\
&= 	t_{w_\mu}^{-\frac12} \sum_{u\in W_\mu} u
	\left( \prod_{\beta^\vee \in R_-^\vee} 
	\frac{1-Y^{\beta^\vee}t_{\beta^\vee}}
	{1-Y^{\beta^\vee}} \right)\mathbf{1}\\
&= t_{w_\mu}^{-\frac12} W_\mu(t)\mathbf{1},
\end{align*}
where the last equality is~\cite[Corollary 2.6(a)]{NR03}.
\qed

\begin{defn}{\em 
The {\em symmetric Macdonald polynomial} $P_\mu\in \bbK[X]^{W_0}$ indexed by $\mu \in (\fh_\bbZ^*)^+$ is defined by
\begin{equation}\label{eqn.defnP}
P_\mu \mathbf{1} 
	= \frac{1}{t_{w_\mu}^{-1/2}W_\mu(t)}\ \mathbf{1}_0\tau_{m_\mu}^\vee \mathbf{1},\qquad
\hbox{where }
\mathbf{1}_0 = \sum_{w\in W_0} t_{w_0w}^{-1/2}T_w.
\end{equation}
See~\cite[(5.5.7) and (5.7.10)]{M03}.
}\end{defn}
The symmetric Macdonald polynomials $\{P_\mu\mathbf{1} \mid \mu\in (\fh_\bbZ^*)^+\}$ is a basis for the $W_0$-invariant polynomials 
$$\bbK[X]^{W_0}\mathbf{1}= \{f\mathbf{1} \mid wf\mathbf{1} = f\mathbf{1}\hbox{ for all } w\in W_0\}.$$
\begin{remark}{\em
The polynomial $P_\mu$ is normalized so that the coefficient of $X^\mu$ in the monomial expansion of $P_\mu$ is $1$.  The normalization chosen here is different from that in~\cite{RY10}.
}\end{remark}


The following result provides an expression for $\mathbf{1}_0$ in terms of intertwining operators.  From this, $P_\mu$ can be expressed as a linear combination of $E_\nu$ for $\nu \in W_0\mu$.  
\begin{prop}\label{prop.10-tau}\label{cor.P-E} \cite[p.\! 203]{M94}. Also see~\cite[(5.7.8)]{M03} and~\cite[(4.13)]{C95b}.

Recall $\calL(v,w) = \{a^\vee \in S_+^\vee \mid \ttH_{a^\vee} \hbox{ separates } v\ttA \hbox{ and } w\ttA\}$.
\begin{enumerate}
\item[(a)]
The symmetrizing operator is
$$\mathbf{1}_0 = \sum_{w\in W_0} \tau_w^\vee 
	\left(\prod_{a^\vee \in \calL(w^{-1},w_0)} 
	t_{a^\vee}^{1/2}\left( 
	\frac{1-t_{a^\vee}^{-1}Y^{-a^\vee}}{1-Y^{-a^\vee}} \right)\right),$$

\item[(b)]
For $\mu\in (\fh_\bbZ^*)^+$,
$$P_\mu = 
\sum_{v\in W^\mu} 
\prod_{a^\vee \in m_\mu^{-1}\calL(v^{-1},v_\mu^{-1})} 
t_{a^\vee}^{1/2} \frac{1-q^{\shf(-a^\vee)}t^{\hgt(-a^\vee)} t_{a^\vee}^{-1}}
{1-q^{\shf(-a^\vee)}t^{\hgt(-a^\vee)}} \ E_{v\mu}.
$$
\end{enumerate}
\end{prop}
\proof \
(a) Since $\mathbf{1}_0= \sum_{w\in W_0} t_{w_0w}^{-1/2} T_w$, then $\mathbf{1}_0$ can be written in the form $\sum_{w\in W_0} \tau_w^\vee b_w$, where each $b_w$ is a rational function in $Y$.
The coefficient of $T_{w_0}$ in $\mathbf{1}_0$ is $1$, hence $b_{w_0}=1$.
The other $b_w$ can be computed by induction on the length of $w$, comparing coefficients in
$$\sum_{w\in W_0}\tau_w^\vee b_w t_i^{1/2}
= \mathbf{1}_0 t_i^{1/2}
= \mathbf{1}_0 T_i
= \sum_{w\in W_0}\tau_w^\vee b_w T_i
= \sum_{w\in W_0}\tau_w^\vee  
	b_w\left(\tau_i^\vee -\frac{t_i^{-1/2}-t_i^{1/2}}{1-Y^{-\alpha_i^\vee}}\right),
$$
using $\mathbf{1}_0T_i = \mathbf{1}_0 t_i^{1/2}$ for $i=1,\ldots, n$
from~\cite[(5.5.9)]{M03}.

(b) 
Each element $w$ in $W_0$ has a unique factorization of the form $w=vu$ where $v\in W^\mu$ and $u\in W_\mu$, so following (a),
 \begin{align*}
\mathbf{1}_0 
&= \sum_{w\in W_0} \tau_w^\vee b_{(w^{-1}, w_0^{-1})}
= \left(\sum_{v\in W^\mu} \tau_v^\vee  
	b_{(v^{-1},v_\mu^{-1})}\right)
	\left(\sum_{u\in W_\mu}\tau_u^\vee 
	b_{(u^{-1},w_\mu^{-1})}\right),
\end{align*}
where 
\begin{equation}\label{eqn.bvw}
b_{(v,w)} = \prod_{a^\vee \in \calL(v,w)} t_{a^\vee}^{1/2} \frac{1-t_{a^\vee}^{-1}Y^{-a^\vee}}{1-Y^{-a^\vee}}.
\end{equation}
Applying $\mathbf{1}_0$ to $\tau_{m_\mu}^\vee\mathbf{1}$, and using the fact that 
$T_u f\mathbf{1} = t_u^{1/2} f\mathbf{1}$ if $u (f\mathbf{1}) = f\mathbf{1}$, then 
\begin{align*}
\mathbf{1}_0 \tau_{m_\mu}^\vee\mathbf{1}
&= \left(\sum_{v\in W^\mu} \tau_v^\vee  
	b_{(v^{-1},v_\mu^{-1})}\right)
	\left(\sum_{u\in W_\mu} t_{w_\mu u}^{-1/2}T_u \right) 
	\tau_{m_\mu}^\vee \mathbf{1}
= t_{w_\mu}^{-1/2} W_\mu(t)
	\left(\sum_{v\in W^\mu} \tau_v^\vee  b_{(v^{-1},v_\mu^{-1})}
	\tau_{m_\mu}^\vee  \right)\mathbf{1}\\
&= t_{w_\mu}^{-1/2} W_\mu(t) \sum_{v\in W^\mu} 
	\left(\prod_{a^\vee \in m_\mu^{-1}\calL(v^{-1},v_\mu^{-1})} 
	t_{a^\vee}^{1/2}
	\frac{1-q^{\shf(-a^\vee)}t^{\hgt(-a^\vee)}t_{a^\vee}^{-1}}
	{1-q^{\shf(-a^\vee)}t^{\hgt(-a^\vee)}} \right)
	E_{v\mu}\mathbf{1}.
\end{align*}
Observe that $m_\mu^{-1} v_\mu^{-1} = (v_\mu m_\mu)^{-1} = x^{-w_0\mu}$.
\qed

The following calculation shows that if the weights $\mu$ and $\nu$ lie in the same $W_0$-orbit, then $\mathbf{1}_0 \tau_{m_\mu}^\vee\mathbf{1}$ and $\mathbf{1}_0 \tau_{m_\nu}^\vee\mathbf{1}$ differ by a scalar multiple.
\begin{prop}\label{prop.domP}
Let $\mu\in (\fh_\bbZ^*)^+$ and $v\in W^\mu$ with $vm_\mu > m_\mu$.  
Then
$$\mathbf{1}_0\tau_v^\vee \tau_{m_\mu}^\vee\mathbf{1}
= 
\left(\prod_{a^\vee \in \calL(m_\mu^{-1}, m_\mu^{-1}v^{-1})} 
t_{a^\vee}^{-\frac12} \frac{1-q^{\shf(-a^\vee)}t^{\hgt(-a^\vee)}t_{a^\vee}}{1-q^{\shf(-a^\vee)}t^{\hgt(-a^\vee)}}\right) \mathbf{1}_0 \tau_{m_\mu}^\vee\mathbf{1},
$$
where $\calL(v,w) = \{a^\vee\in S_+^\vee \mid \ttH_{a^\vee} \hbox{ separates $v$ and $w$} \}$.
\end{prop}
\proof 
By~\cite[(5.5.9)]{M03}, 
$\mathbf{1}_0T_i = \mathbf{1}_0 t_i^{\frac12}$ for $i=1,\ldots, n$, so
$$\mathbf{1}_0 \tau_i^\vee \tau_w^\vee \mathbf{1} 
= \mathbf{1}_0 \left(t_i^{\frac12} + \frac{t_i^{-\frac12}-t_i^{\frac12}}{1-Y^{-\alpha_i^\vee}}\right) \tau_w^\vee \mathbf{1} 
= \mathbf{1}_0 \tau_w^\vee \left(t_i^{\frac12} + \frac{t_i^{-\frac12}-t_i^{\frac12}}{1-q^{\shf(-w^{-1}\alpha_i^\vee)} t^{\hgt(-w^{-1}\alpha_i^\vee)}} \right)\mathbf{1}.$$
Let $v=s_{i_1} \cdots s_{i_r}\in W^\mu$ be a reduced expression, so that $\mathbf{1}_0 \tau_v^\vee \tau_{m_\mu}^\vee\mathbf{1}
= \mathbf{1}_0 \tau_{i_1}^\vee \cdots \tau_{i_r}^\vee \tau_{m_\mu}^\vee \mathbf{1}.$
Starting with $\tau_{i_1}^\vee$, commuting each $\tau_{i_j}^\vee$ from left to right gives
\begin{align*}
\mathbf{1}_0\tau_v^\vee \tau_{m_\mu}^\vee \mathbf{1}
&= \mathbf{1}_0 \tau_{m_\mu}^\vee \prod_{j=1}^r 
	\left(t_{i_j}^{\frac12}+\frac{t_{i_j}^{-\frac12}-t_{i_j}^{\frac12}}
	{1-Y^{-m_\mu^{-1}s_{i_r} \cdots s_{i_{j+1}} 
	\alpha_{i_j}^\vee}} \right)\mathbf{1} \\
&= \mathbf{1}_0 \tau_{m_\mu}^\vee 
	\prod_{a^\vee\in \calL(m_\mu^{-1},m_\mu^{-1}v^{-1})} 
	\left(t_{a^\vee}^{\frac12}
	+\frac{t_{a^\vee}^{-\frac12}-t_{a^\vee}^{\frac12}}
	{1-q^{\shf(-a^\vee)}t^{\hgt(-a^\vee)}} \right)\mathbf{1}.
\end{align*}
\qed

See Remark~\ref{rem.goodeh} for an alternate formulation of this statement.

\section{Multiplication formulas}\label{sec.lrrule}
In this section, all alcove walks under consideration are contained in the closure $\overline{\ttC}$ of the dominant chamber.

\subsection{Littlewood-Richardson formulas}
Various multiplication formulas for Macdonald polynomials can be obtained by applying Theorem~\ref{thm.X-tau} to $\mathbf{1}$.

\begin{corollary}\label{cor.XE-E}
Let $\mu,\lambda \in \fh_\bbZ^*$, and fix a reduced expression 
$x^\mu = s_{i_r}^\vee \cdots s_{i_1}^\vee (\pi_j^\vee)^{-1}$.
Then
$$ X^\mu E_\lambda
= \sum_{h\in \Gamma^{\overline\ttC}(\vec{x^{-\mu}}, m_\lambda^{-1})}
(-1)^{\phi_-(h)}g_h n_h E_{\varpi(h)},$$
where the sum is over the set of alcove walks of type $\vec{x^{-\mu}}=(\pi_j^\vee, i_1,\ldots, i_r)$ beginning in $m_\lambda^{-1}$ and contained in
$\overline\ttC$, 
\begin{align*}
g_h &=\left(\prod_{k\in\phi_{\mathrm{o}}(h)} t_{i_k}^{\epsilon_k/2}\right)
	\left(\prod_{k\in\phi_\mathrm{aff}(h)} 
	\frac{t_{i_k}^{-1/2}-t_{i_k}^{1/2}}
		{1-q^{\shf(-h_k^\vee)}t^{\hgt(-h_k^\vee)}} \right)	
	\left(\prod_{k\in \psi_\mathrm{aff}(h)} q^{\shf(-h_k^\vee)}t^{\hgt(-h_k^\vee)}\right),\\
n_h &= \prod_{k\in \xi_-(h)} 
	\frac{1-q^{\shf(-h_k^\vee)}t^{\hgt(-h_k^\vee)}t_{i_k}^{-1}}
	{1-q^{\shf(-h_k^\vee)}t^{\hgt(-h_k^\vee)}}
	\frac{1-q^{\shf(-h_k^\vee)}t^{\hgt(-h_k^\vee)}t_{i_k}}
	{1-q^{\shf(-h_k^\vee)}t^{\hgt(-h_k^\vee)}},
\end{align*}
$h_k^\vee$ indexes the $k$th active hyperplane of $h$, and 
$\epsilon_1,\ldots, \epsilon_r$ are as defined in~\eqref{eqn.epsilonk}.
\end{corollary}
\proof
Apply Theorem~\ref{thm.X-tau} to $\mathbf{1}$ with $z=x^\mu$, $w=m_\lambda$ , and the proof is the same as for Corollary~\ref{cor.X-tau}. \qed

Next, we examine multiplication involving symmetric Macdonald polynomials.  The following notation will be used throughout the rest of this section.

For $\mu \in \fh_\bbZ^*$, $\lambda \in (\fh_\bbZ^*)^+$, fix a reduced expression for the minimal length representative $m_\mu = s_{i_r}^\vee \cdots s_{i_1}^\vee (\pi_j^\vee)^{-1}$ of the coset $x^\mu W_0$.  For $v\in W^\lambda$, let $\Gamma_2^{\overline\ttC}\left(\vec{m}_\mu^{-1}, (vm_\lambda)^{-1}\right)$ denote the set of alcove walks satisfying the following properties:
\begin{enumerate}
\item has type $\vec{m}_\mu^{-1} = (\pi_j^\vee, i_1,\ldots, i_r)$,
\item begins in $(vm_\lambda)^{-1}$,
\item is contained in $\overline\ttC$,
\item each fold is coloured either black or grey.
\end{enumerate}

Given $h\in \Gamma_2^{\overline\ttC}$, let
\begin{eqnarray}
&&\label{eqn.ph}
\begin{array}{l}
	\hbox{$p(h)$ be the walk obtained by straightening all grey folds of $h$,}\\
	\qquad\hbox{and translated so that it ends in $1$,} 
	\end{array}\\
&&\label{eqn.epsilonkh}
\epsilon_k(h) = \begin{cases} 
	+1 & \hbox{if the $k$th step of $p(h)$ is positive,}\\
	-1 & \hbox{if the $k$th  step of $p(h)$ is negative.}
	\end{cases}
\end{eqnarray}
See Example~\ref{eg.2} for an illustration.

\begin{theorem}\label{thm.EP-E}
Let $\mu \in \fh_\bbZ^*$, $\lambda \in (\fh_\bbZ^*)^+$, and fix a reduced expression for the minimal length representative $m_\mu = s_{i_r}^\vee \cdots s_{i_1}^\vee (\pi_j^\vee)^{-1}$ of the coset $x^\mu W_0$.  Then
$$E_\mu P_\lambda
= \sum_{v\in W^\lambda} 
\sum_{h \in \Gamma_2^{\overline\ttC}\left(\vec{m}_\mu^{-1}, (vm_\lambda)^{-1}\right)
} 
\left(
(-1)^{|\phi_{\mathrm{grey}}^-(h)|} b_h\prod_{k=1}^r c_{k(h)}\right) 
E_{\varpi(h)},$$
where $|\phi_{\mathrm{grey}}^-(h)|$ is the number negative grey folds of $h$,
\begin{equation}
b_h =  \prod_{a^\vee \in \calL(\ttb(h),x^{-w_0\lambda})} 
	t_{a^\vee}^{1/2}
	\frac{1-q^{\shf(-a^\vee)}t^{\hgt(-a^\vee)}t_{a^\vee}^{-1}}
	{1-q^{\shf(-a^\vee)}t^{\hgt(-a^\vee)}},
\qquad
\ttb(h) \hbox{ is where $h$ begins,}
\end{equation}
and the coefficient $c_{k(h)}$ arising from the $k$th step of $h$ is
\begin{equation}\label{eqn.verylonglist}
\begin{array}{ll}
1,	
	& \begin{array}{l}
		\hbox{for a positive crossing,}
		\end{array}\\
\displaystyle 
	\frac{1-q^{\shf(-h_k^\vee)}t^{\hgt(-h_k^\vee)}t_{i_k}^{-1}}
	{1-q^{\shf(-h_k^\vee)}t^{\hgt(-h_k^\vee)}}
	\frac{1-q^{\shf(-h_k^\vee)}t^{\hgt(-h_k^\vee)}t_{i_k}}
	{1-q^{\shf(-h_k^\vee)}t^{\hgt(-h_k^\vee)}},
	& \begin{array}{l}
		\hbox{for a negative crossing,}
		\end{array}\\
\\
\displaystyle t_{i_k}^{\epsilon_k/2},
	& \begin{array}{l}
		\hbox{for a grey fold}\\
		\quad\hbox{touching a wall of }\ttC,\\
		\end{array}\\
\\
\displaystyle
	\frac{(t_{i_k}^{-1/2}-t_{i_k}^{1/2})
		(q^{\shf(-h_k^\vee)}t^{\hgt(-h_k^\vee)})^{\frac12(1-\epsilon_k(h))}}
		{1-q^{\shf(-h_k^\vee)}t^{\hgt(-h_k^\vee)}},
	& \begin{array}{l}
		\hbox{for a negative grey fold}\\
		\quad\hbox{touching an affine hyperplane,}\\
		\end{array}\\
\displaystyle
	\frac{(t_{i_k}^{-1/2}-t_{i_k}^{1/2})
		(q^{\shf(-h_k^\vee)}t^{\hgt(-h_k^\vee)})^{\frac12(1+\epsilon_k(h))}}
		{1-q^{\shf(-h_k^\vee)}t^{\hgt(-h_k^\vee)}},
	& \begin{array}{l}
		\hbox{for a positive grey fold}\\
		\quad\hbox{touching an affine hyperplane,}\\
		\end{array}\\
\\
\displaystyle 
	\frac{t_{b_k^\vee}^{-1/2}-t_{b_k^\vee}^{1/2}}
		{1-q^{\shf(-b_k^\vee)}t^{\hgt(-b_k^\vee)}},
	& \begin{array}{l}
		\hbox{for a black fold such that the}\\
		\quad\hbox{$k$th step of $p(h)$ is positive,}\end{array}\\
\displaystyle 
	\frac{(t_{b_k^\vee}^{-1/2}-t_{b_k^\vee}^{1/2})
		(q^{\shf(-b_k^\vee)}t^{\hgt(-b_k^\vee)})}
		{1-q^{\shf(-b_k^\vee)}t^{\hgt(-b_k^\vee)}},	
	& \begin{array}{l}
		\hbox{for a black fold such that the}\\
		\quad\hbox{$k$th step of $p(h)$ is negative,}\end{array}
\end{array}
\end{equation}
$\varpi(h)$, $h_k^\vee$, $\epsilon_k(h)$ are respectively defined in~\eqref{eqn.varweight}, \eqref{eqn.pkvee}, \eqref{eqn.epsilonkh}, and $b_k^\vee = s_{i_1}^\vee s_{i_2}^\vee\cdots s_{i_{k-1}}^\vee \alpha_{i_k}^\vee$.
\end{theorem}
\proof
The idea is to expand $E_\mu$ in terms of monomials $X^\nu$ (Theorem~\ref{thm.Xtau-X}), and then expand $X^\nu P_\lambda$ in terms of nonsymmetric Macdonald polynomials (Proposition~\ref{prop.10-tau} and Corollary~\ref{cor.XE-E}).

Since
$$E_\mu = \sum_{p\in \Gamma(\vec{m}_\mu)} f_p t_{\ttd(p)}^{1/2} X^{\ttwt(p)}, \quad\hbox{where} \quad
f_p= \prod_{k\in\phi(p)} \frac{t_{b_k^\vee}^{-1/2}-t_{b_k^\vee}^{1/2}}{1-q^{\shf(-b_k^\vee)}t^{\hgt(-b_k^\vee)}} \prod_{k\in \phi_-(p)} q^{\shf(-b_k^\vee)}t^{\hgt(-b_k^\vee)},$$
with notation from Theorem~\ref{thm.Xtau-X}.  Then
\begin{align*}
E_\mu \mathbf{1}_0 \tau_{m_\lambda}^\vee\mathbf{1}
&= \left(\sum_{p\in \Gamma(\vec{m}_\mu)} 
	f_p X^{\ttwt(p)} T_{\ttd(p)}^{1/2}\right)
	\mathbf{1}_0 \tau_{m_\lambda}^\vee\mathbf{1}
= \sum_{p\in \Gamma(\vec{m}_\mu)} 
	f_p X^{\tte(p)} 
	\mathbf{1}_0 \tau_{m_\lambda}^\vee\mathbf{1}\\
&= t_{w_\lambda}^{-1/2}W_{\lambda}(t) \sum_{v\in W^\lambda} b_{((vm_\lambda)^{-1},x^{-w_0\lambda})}
	\sum_{p\in \Gamma(\vec{m}_\mu)} 
	f_p X^{\tte(p)} 
	\tau_{vm_\lambda}^\vee\mathbf{1},
\end{align*}
where 
$\displaystyle 
b_{((vm_\lambda)^{-1}, x^{-w_0\lambda})} 
= \prod_{a^\vee \in m_\lambda^{-1}\calL(v^{-1},v_\lambda^{-1})} 
	t_{a^\vee}^{1/2}
	\frac{1-q^{\shf(-a^\vee)}t^{\hgt(-a^\vee)}t_{a^\vee}^{-1}}
	{1-q^{\shf(-a^\vee)}t^{\hgt(-a^\vee)}},
$
as in Proposition~\ref{prop.10-tau}.

The next step is to express $X^{\tte(p)}\tau_{vm_\lambda}^\vee\mathbf{1}$ in terms of nonsymmetric Macdonald polynomials.  The way this is achieved in Corollary~\ref{cor.XE-E} is to interpret $X^{\tte(p)}$ as a minimal length walk (without any folds) of type $\tte(p)^{-1}$ and beginning in the alcove $(vm_\lambda)^{-1}$, then the terms in the expansion of $X^{\tte(p)}\tau_{vm_\lambda}^\vee\mathbf{1}$ are generated by folding this walk in all possible ways (and only keeping those which are contained in the dominant chamber).

Instead, we interpret $f_p X^{\tte(p)}$ as a walk (possibly having folds) of type $m_\mu^{-1}$ and beginning in $(vm_\lambda)^{-1}$.  Then the terms in the expansion of $f_pX^{\tte(p)}\tau_{vm_\lambda}^\vee\mathbf{1}$ in the nonsymmetric Macdonald polynomial basis are generated by folding this walk in all possible ways (and keeping those contained in the dominant chamber).  These newly introduced folds contribute coefficients which are different from those coming from the folds of the seminal walk, so we keep track of the new folds by colouring them grey.  What remains to be done is to unify the previous concepts of various kinds of crossings and foldings from Theorem~\ref{thm.Xtau-X} and Corollary~\ref{cor.XE-E}.

That is,
\begin{align*}
f_pX^{\tte(p)}\tau_{vm_\lambda}^\vee\mathbf{1}
&= f_{p(h)} \sum_{h} (-1)^{|\phi_{\mathrm{grey}}^-(h)|} g_h n_h E_{\varpi(h)}\mathbf{1}
\end{align*}
where the sum is over all alcove walks of type $\vec{m}_\mu^{-1}$ beginning in $(vm_\lambda)^{-1}$, which are contained in $\overline\ttC$, such that the $k$th step is a black fold if and only if the $k$th step of $p(h)$ is a fold (see~\eqref{eqn.ph}), and all other folds are grey.  Here, the coefficient $g_h$ involves grey folds only (see Corollary~\ref{cor.XE-E}), and
$$f_{p(h)} = 
\prod_{k\in\phi(p(h))} \frac{t_{b_k^\vee}^{-1/2}-t_{b_k^\vee}^{1/2}}{1-q^{\shf(-b_k^\vee)}t^{\hgt(-b_k^\vee)}} \prod_{k\in \phi_-(p(h))} q^{\shf(-b_k^\vee)}t^{\hgt(-b_k^\vee)}, 
\qquad\hbox{for }
b_k^\vee = s_{i_1}^\vee \cdots s_{i_{k-1}}^\vee \alpha_{i_k}^\vee.
$$
Thus,
\begin{align*}
E_\mu P_\lambda\mathbf{1}
&= \sum_{v\in W^\lambda} 
	\sum_{p\in \Gamma(\vec{m}_\mu)} 
	b_{((vm_\lambda)^{-1},x^{-w_0\lambda})} f_p X^{\tte(p)} 
	\tau_{vm_\lambda}^\vee\mathbf{1},\\
&= \sum_{v\in W^\lambda} 
	\sum_{h\in \Gamma_2^{\overline\ttC}(\vec{m}_\mu^{-1}, (vm_\lambda)^{-1})} 
	\left((-1)^{|\phi_\mathrm{grey}^-(h)|} 
	b_h  
	\sum_{k=1}^r c_{k(h)}\right) E_{\varpi(h)}\mathbf{1}.
\end{align*}

The steps of $h$ may be a crossing, a grey fold, or a black fold, and depending on the step, the coefficient $c_k(h)$ contributed by the $k$th step of $h$ is summarized in table~\eqref{eqn.verylonglist}.
\qed

The next result expresses the product of two symmetric Macdonald polynomials in terms of symmetric Macdonald polynomials.  For the purposes of simplifying the notation, we assume for the remainder of this section that the walls of $\ttC$ are labeled by the {\em negative} coroots $\{-\alpha_1^\vee, \ldots, -\alpha_n^\vee\}$ instead of the usual positive coroots.
The labeling of the affine hyperplanes remain unchanged from before; they are labeled by positive affine coroots $\{ -\alpha^\vee + jd \mid \alpha^\vee \in R_+^\vee, \ j\in \bbZ_{\geq1}\}$.

\begin{remark}\label{rem.goodeh}{\em  
It is now useful to phrase Proposition~\ref{prop.domP} another way.  
Given $\mu \in (\fh_\bbZ^*)^+$, 
\begin{align*}
\mathbf{1}_0 \tau_{m_\mu}^\vee\mathbf{1} 
&= \left(t_{w_\mu}^{-\frac12}W_\mu(t)\right) P_\mu\mathbf{1}
= \left(\prod_{\alpha^\vee\in \calL(m_\mu^{-1}, m_\mu^{-1}w_\mu^{-1})} 
	t^{-\frac12}_{\alpha^\vee} 
	\frac{1-q^{\shf(-\alpha^\vee)}t^{\hgt(-\alpha^\vee)}t_{\alpha^\vee}}
	{1-q^{\shf(-\alpha^\vee)}t^{\hgt(-\alpha^\vee)}}\right)
	P_\mu \mathbf{1},
\end{align*}
(this is where we use the convention that the walls of $\ttC$ are labeled by negative coroots rather than positive coroots), then for $z\in W^\mu$,
\begin{equation}\label{eqn.goodeh}
\mathbf{1}_0 \tau_z^\vee \tau_{m_\mu}^\vee\mathbf{1}
=	\left(\prod_{a^\vee \in\calL(m_\mu^{-1}w_\mu^{-1}, m_\mu^{-1}z^{-1})} 
	t^{-\frac12}_{a^\vee}
	\frac{1-q^{\shf(-a^\vee)}t^{\hgt(-a^\vee)}t_{a^\vee}}
	{1-q^{\shf(-a^\vee)}t^{\hgt(-a^\vee)}} \right)
	P_\mu\mathbf{1}.
\end{equation}
Observe that $m_\mu^{-1} w_\mu^{-1} = \left(v_\mu x^{-\mu}\right) w_\mu^{-1} = x^{-w_0\mu} w_0.$
}\end{remark}

\begin{theorem}\label{thm.PP-P}
Let $\mu, \lambda \in (\fh_\bbZ^*)^+$, and fix a reduced expression for the minimal length representative $m_\mu = s_{i_r}^\vee \cdots s_{i_1}^\vee (\pi_j^\vee)^{-1}$ of the coset $x^\mu W_0$.  With the same notation from Theorem~\ref{thm.EP-E}, then
$$P_\mu P_\lambda
= \frac{1}{t_{w_\mu}^{-1/2}W_\mu(t)}
	\sum_{v\in W^\lambda} 
	\sum_{h \in \Gamma_2^{\overline\ttC}\left( \vec{m}_\mu^{-1}, (vm_\lambda)^{-1}\right)}
	\left((-1)^{\phi_{grey}^-(h)} 
	b_h e_h \prod_{k=1}^r c_{k(h)} \right)
	P_{-w_0\ttwt(h)},$$
where
\begin{equation}
e_h = \prod_{a^\vee \in \calL(x^{\ttwt(h)}w_0,\tte(h))} 
	t_{a^\vee}^{-\frac12}
	\frac{1-q^{\shf(-a^\vee)}t^{\hgt(-a^\vee)}t_{a^\vee}}
	{1-q^{\shf(-a^\vee)}t^{\hgt(-a^\vee)}},
\qquad
\tte(h) \hbox{ is where $h$ ends.}
\end{equation}
\end{theorem}

\proof
Since $t_{w_\mu}^{-\frac12}W_\mu(t) P_\mu\mathbf{1}_0 = \mathbf{1}_0 E_\mu\mathbf{1}_0$, then the result can obtained by applying the operator $\mathbf{1}_0$ to both sides of the equation in Theorem~\ref{thm.EP-E}.

The weight $\varpi(h)$ defined by $\tte(h)^{-1} = m_{\varpi(h)}$ is not necessarily dominant, so $\mathbf{1}_0 E_{\varpi(h)}\mathbf{1}$ is a certain scalar multiple of a symmetric Macdonald polynomial.  To find this scalar, let $z\in W_0$ be of minimal length such that $z\varpi(h)_+ = \varpi(h)$.  Note that $m_{\varpi(h)_+}^{-1} = m_{-w_0\varpi(h)_+}$, as $\varpi(h)_+$ is a dominant weight.  Moreover, $m_{-w_0\varpi(h)_+}$ is the minimal length representative of the coset $x^{\ttwt(h)}W_0$, because the endpoint $\tte(h)$ is in that coset by definition, therefore 
$$\ttwt(h) = -w_0\varpi(h)_+. $$
Hence remark~\ref{rem.goodeh} gives
\begin{align*}
\mathbf{1}_0 E_{\varpi(h)} \mathbf{1}
&=\mathbf{1}_0\tau_z^\vee \tau_{m_{\varpi(h)_+}}^\vee \mathbf{1}
= \left(\prod_{a^\vee \in \calL(x^{\ttwt(h)}w_0, \tte(h))} 
	t_{a^\vee}^{-\frac12}
	\frac{1-q^{\shf(-a^\vee)}t^{\hgt(-a^\vee)}t_{a^\vee}}
	{1-q^{\shf(-a^\vee)}t^{\hgt(-a^\vee)}}\right)
	P_{-w_0\ttwt(h)} \mathbf{1},
\end{align*}
since $\tte(h) = m_{\varpi(h)}^{-1} = m_{\varpi(h)_+}^{-1}z^{-1}$.
and
\begin{align*}
P_\mu P_\lambda \mathbf{1}
&= \frac{1}{t_{w_\mu}^{-\frac12}W_\mu(t)} \mathbf{1}_0 E_\mu P_\lambda \mathbf{1}
= \frac{1}{t_{w_\mu}^{-\frac12}W_\mu(t)}
	\sum_{v\in W^\lambda} 
	\sum_{h \in \Gamma_2^{\overline\ttC}\left( \vec{m}_\mu^{-1}, (vm_\lambda)^{-1}\right)} a_h(q,t)e_h
	P_{-w_0\ttwt(h)}\mathbf{1},
\end{align*}
where $a_h(q,t)=	\left((-1)^{\phi_{grey}^-(h)} 
	b_h
	\prod_{k=1}^r c_{k(h)} \right)$.
\qed

\begin{remark}{\em
Theorem~\ref{thm.Xtau-X} gives the expansion of $E_\mu$ in terms of monomials, so in theory, the product of nonsymmetric Macdonald polynomials may be computed as a sum over pairs of paths.
That is, since 
$$E_\mu = \sum_{p\in \Gamma(\vec{m}_\mu)} f_p t_{\ttd(p)}^{\frac12} X^{\ttwt(p)}, \quad\hbox{where} \quad
f_p= \prod_{k\in\phi(p)} \frac{t_{b_k^\vee}^{-\frac12}-t_{b_k^\vee}^{\frac12}}{1-q^{\shf(-b_k^\vee)}t^{\hgt(-b_k^\vee)}} \prod_{k\in \phi_-(p)} q^{\shf(-b_k^\vee)}t^{\hgt(-b_k^\vee)},$$
then
\begin{align*}
E_\mu E_\lambda\mathbf{1}
= E_\mu \tau_{m_\lambda}^\vee \mathbf{1}
&= \sum_{p\in \Gamma(\vec{m}_\mu)} f_p t_{\ttd(p)}^{\frac12} X^{\ttwt(p)} \tau_{m_\lambda}^\vee \mathbf{1}\\
&=\sum_{p\in \Gamma(\vec{m}_\mu)} f_p t_{\ttd(p)}^{\frac12} 
 \sum_{h\in \Gamma(\vec{x}^{-\ttwt(p)}, m_\lambda^{-1})}
(-1)^{\phi_-(h)}g_h n_h E_{\varpi(h)} 
\mathbf{1},
\end{align*} 
with the same notations from Theorem~\ref{thm.Xtau-X} and Corollary~\ref{cor.XE-E}.
}\end{remark}

\begin{remark}{\em \
\begin{enumerate}
\item Using interpolation Macdonald polynomials, Baratta obtained type A Pieri-type formulas~\cite[Proposition 8, 10]{B09} for the expansion of $E_\mu(q,t) P_{\omega_1}(0,0)$ and $E_\mu(q,t) P_{\omega_n}(0,0)$ in terms of nonsymmetric Macdonald polynomials.  Also see~\cite{La08}.
\item Haglund, Luoto, Mason, and van Willigenburg considered the type A case at $q=t=0$, when $P_\lambda(0,0)$ is a Schur polynomial and $E_\lambda(0,0)$ is a Demazure character, and obtained a formula for the expansion of $E_\mu(0,0) P_\lambda(0,0)$ in terms of $E_\mu(0,0)$ with positive coefficients~\cite[Theorem 6.1]{HLMvW09}.  The coefficients $\sum_{h: \varpi(h) = \gamma} a_h(0,0)$ count certain fillings of skew tableau-like diagrams called skyline diagrams.
\item Recently, a generalized Macdonald operator was introduced by van Diejen and Emsiz~\cite{vDE10}, who used it to obtain Pieri formulas for Macondald polynomials of arbitrary type, and also a Littlewood-Richardson formula for `small weights'.
\end{enumerate}
}\end{remark}

\begin{example}\label{eg.2}{\em
See section \ref{sec.eg} for more details on the alcove picture of type $\fsl_2$.

The following walks $h_1$ and $h_2 \in \Gamma_2^{\overline\ttC}(\vec{m}_{-8\omega}^{-1} , (s_1m_{2\omega})^{-1})$ are used in the expansion of $E_{-8\omega}P_{2\omega}$ in terms of nonsymmetric Macdonald polynomials.  Each is of type $\vec{m}_{-8\omega}^{-1} = (0,1,0,1,0,1,0,1)$ and begins in $(s_1 m_{2\omega})^{-1} = x^{2\omega}$.
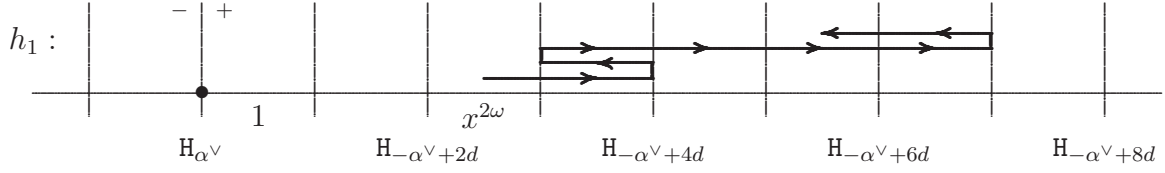
$$\beginpicture
\setcoordinatesystem units <1.5cm,1cm>         
\setplotarea x from -2 to 9, y from -1.2 to 1.5  
\put{$\bullet$} at 0 0 
{\small
\put{$\ttH_{\alpha^\vee}$}[t] at 0 -0.6
\put{$\ttH_{-\alpha^\vee+2d}$}[t] at 2 -0.6
\put{$\ttH_{-\alpha^\vee+4d}$}[t] at 4 -0.6
\put{$\ttH_{-\alpha^\vee+6d}$}[t] at 6 -0.6 
\put{$\ttH_{-\alpha^\vee+8d}$}[t] at 8 -0.6 
}\plot -1.5 0  8.5 0 /  
\plot  -1 1.2  -1 -0.3 / \plot  0 1.2  0 -0.3 / \plot  1 1.2  1 -0.3
/ \plot  2 1.2  2 -0.3 / \plot  3 1.2  3 -0.3 / \plot  4 1.2  4 -0.3
/ \plot 5 1.2 5 -0.3 / \plot  6 1.2  6 -0.3 /
\plot  7 1.2  7 -0.3 / \plot  8 1.2  8 -0.3 /
\put{$\scriptstyle{+}$} at 0.2 1.1 \put{$\scriptstyle{-}$} at -0.2 1.1
\put{$1$} at 0.5 -0.3
\put{$x^{2\omega}$} at 2.5 -0.3
\put{$h_1:$} at -1.5 0.7
\setplotsymbol(.)
\arrow <6pt> [.2,.67] from 2.5 0.2 to 3.5 0.2
\plot 3.5 0.2 4 0.2 / \plot 3.99 0.2 3.99 0.4 /
\arrow <6pt> [.2,.67] from 4 0.4 to 3.5 0.4
\plot 3.5 0.4 3 0.4 / \plot 3.01 0.4 3.01 0.6 /
\arrow <6pt> [.2,.67] from 3 0.6 to 3.5 0.6
\arrow <6pt> [.2,.67] from 3.5 0.6 to 4.5 0.6
\arrow <6pt> [.2,.67] from 4.5 0.6 to 5.5 0.6
\arrow <6pt> [.2,.67] from 5.5 0.6 to 6.5 0.6
\plot 6.5 0.6 7 0.6 / \plot 6.99 0.6 6.99 0.8 /
\arrow <6pt> [.2,.67] from 7 0.8 to 6.5 0.8
\arrow <6pt> [.2,.67] from 6.5 0.8 to 5.5 0.8
\endpicture
$$
The endpoint is in $\tte(h_1) = s_0^\vee s_1 s_0^\vee s_1 s_0^\vee = m_{6\omega}^{-1}$, so $\varpi(h_1)= 6\omega$.

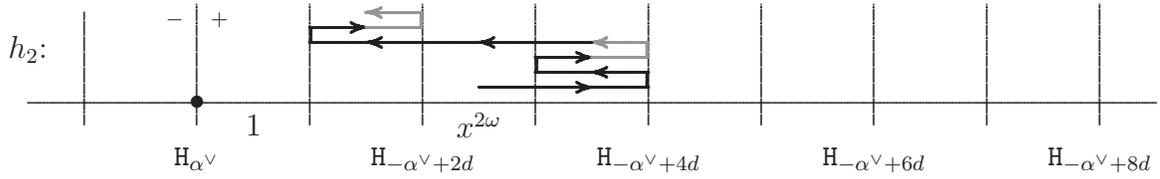
$$\beginpicture
\setcoordinatesystem units <1.5cm,1cm>         
\setplotarea x from -2 to 9, y from -1 to 1.5  
\put{$\bullet$} at 0 0 
{\small
\put{$\ttH_{\alpha^\vee}$}[t] at 0 -0.6
\put{$\ttH_{-\alpha^\vee+2d}$}[t] at 2 -0.6
\put{$\ttH_{-\alpha^\vee+4d}$}[t] at 4 -0.6
\put{$\ttH_{-\alpha^\vee+6d}$}[t] at 6 -0.6 
\put{$\ttH_{-\alpha^\vee+8d}$}[t] at 8 -0.6 
}\plot -1.5 0  8.5 0 /  
\plot -1 1.2  -1 -0.3 / \plot 0 1.3  0 -0.3 / \plot 1 1.3  1 -0.3
/ \plot 2 1.3  2 -0.3 / \plot 3 1.3  3 -0.3 / \plot 4 1.3  4 -0.3
/ \plot 5 1.3 5 -0.3 / \plot 6 1.3  6 -0.3 /
\plot 7 1.3  7 -0.3 / \plot 8 1.3  8 -0.3 /
\put{$\scriptstyle{+}$} at 0.2 1.1 \put{$\scriptstyle{-}$} at -0.2 1.1
\put{$1$} at 0.5 -0.3
\put{$x^{2\omega}$} at 2.5 -0.3
\put{$h_2$:} at -1.5 0.7
\setplotsymbol(.)
\arrow <6pt> [.2,.67] from 2.5 0.2 to 3.5 0.2
\plot 3.5 0.2 4 0.2 / \plot 3.99 0.2 3.99 0.4 /
\arrow <6pt> [.2,.67] from 4 0.4 to 3.5 0.4
\plot 3.5 0.4 3 0.4 / \plot 3.01 0.4 3.01 0.6 /
\arrow <6pt> [.2,.67] from 3 0.6 to 3.5 0.6
\textcolor{Gray}{\plot 3.5 0.6 4 0.6 / \plot 3.99 0.6 3.99 0.8 /
\arrow <6pt> [.2,.67] from 4 0.8 to 3.5 0.8
}\arrow <6pt> [.2,.67] from 3.5 0.8 to 2.5 0.8
\arrow <6pt> [.2,.67] from 2.5 0.8 to 1.5 0.8
\plot 1.5 0.8 1 0.8 / \plot 1.01 0.8 1.01 1 /
\arrow <6pt> [.2,.67] from 1 1 to 1.5 1
\textcolor{Gray}{\plot 1.5 1 2 1 / \plot 1.99 1 1.99 1.2 /
\arrow <6pt> [.2,.67] from 2 1.2 to 1.5 1.2
}\endpicture
$$
The endpoint is in $\tte(h_2) = s_0^\vee = m_{2\omega}^{-1}$, so $\varpi(h_2) = 2\omega$.  Notice that $h_2$ is obtained from $h_1$ by folding the fourth and eighth steps; these are indicated in grey.

Since $h_2$ is generated from $h_1$, then $p=p(h_1)=p(h_2)$ defined in~\eqref{eqn.ph}, is 
$$\beginpicture
\setcoordinatesystem units <1.5cm,1cm>         
\setplotarea x from -2 to 9, y from -1 to 1.5  
\put{$\bullet$} at 0 0 
{\small
\put{$\ttH_{\alpha^\vee}$}[t] at 0 -0.6
\put{$\ttH_{-\alpha^\vee+2d}$}[t] at 2 -0.6
\put{$\ttH_{-\alpha^\vee+4d}$}[t] at 4 -0.6
\put{$\ttH_{-\alpha^\vee+6d}$}[t] at 6 -0.6 
\put{$\ttH_{-\alpha^\vee+8d}$}[t] at 8 -0.6 
}\plot -1.5 0  8.5 0 /  
\plot  -1 1.2  -1 -0.3 / \plot  0 1.2  0 -0.3 / \plot  1 1.2  1 -0.3
/ \plot  2 1.2  2 -0.3 / \plot  3 1.2  3 -0.3 / \plot  4 1.2  4 -0.3
/ \plot 5 1.2 5 -0.3 / \plot 6 1.2 6 -0.3 / \plot 7 1.2  7 -0.3 /
\plot 8 1.2 8 -0.3 /
\put{$\scriptstyle{+}$} at 0.2 1.1 \put{$\scriptstyle{-}$} at -0.2 1.1
\put{$1$} at 0.5 -0.3
\put{$x^{4\omega}s_1$} at 3.5 -0.3
\put{$p:$} at -1.5 0.7
\setplotsymbol(.)
\arrow <6pt> [.2,.67] from 3.5 0.8 to 2.5 0.8
\plot 2.5 0.8 2 0.8 / \plot 2.01 0.6 2.01 0.8 /
\arrow <6pt> [.2,.67] from 2 0.6 to 2.5 0.6
\plot 2.5 0.6 3 0.6 / \plot 2.99 0.4 2.99 0.6 /
\arrow <6pt> [.2,.67] from 3 0.4 to 2.5 0.4
\arrow <6pt> [.2,.67] from 2.5 0.4 to 1.5 0.4
\arrow <6pt> [.2,.67] from 1.5 0.4 to 0.5 0.4
\arrow <6pt> [.2,.67] from 0.5 0.4 to -0.5 0.4
\plot -0.5 0.4 -1 0.4 / \plot -0.99 0.2 -0.99 0.4 /
\arrow <6pt> [.2,.67] from -1 0.2 to -0.5 0.2
\arrow <6pt> [.2,.67] from -0.5 0.2 to 0.5 0.2
\endpicture
$$
and the crossing steps give $(\epsilon_1, \epsilon_4, \epsilon_5, \epsilon_6, \epsilon_8) = (-1 ,-1, -1, -1, +1)$.

For $k=1,\ldots, 8$, the coroots $b_k^\vee = s_{i_1}^\vee \cdots s_{i_{k-1}}^\vee \alpha_{i_k}^\vee = -\alpha^\vee + kd$ is the sequence of labels of hyperplanes crossed by the walk of type $\vec{m}_{-8\omega}^{-1}$ beginning in $1$:
$$
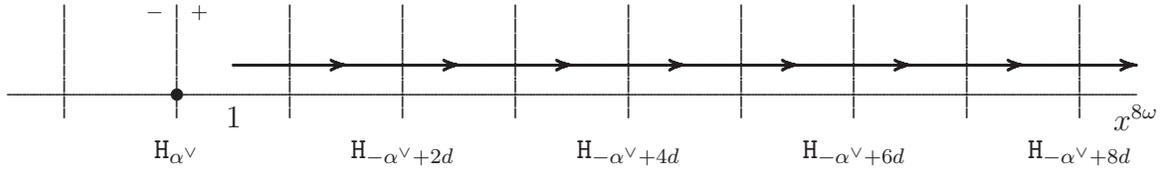
\beginpicture
\setcoordinatesystem units <1.5cm,1cm>         
\setplotarea x from -2 to 9, y from -1.2 to 1.5  
\put{$\bullet$} at 0 0 
{\small
\put{$\ttH_{\alpha^\vee}$}[t] at 0 -0.6
\put{$\ttH_{-\alpha^\vee+2d}$}[t] at 2 -0.6
\put{$\ttH_{-\alpha^\vee+4d}$}[t] at 4 -0.6
\put{$\ttH_{-\alpha^\vee+6d}$}[t] at 6 -0.6 
\put{$\ttH_{-\alpha^\vee+8d}$}[t] at 8 -0.6 
}\plot -1.5 0  8.5 0 /  
\plot  -1 1.2  -1 -0.3 / \plot  0 1.2  0 -0.3 / \plot  1 1.2  1 -0.3
/ \plot  2 1.2  2 -0.3 / \plot  3 1.2  3 -0.3 / \plot  4 1.2  4 -0.3
/ \plot 5 1.2 5 -0.3 / \plot 6 1.2 6 -0.3 / \plot 7 1.2  7 -0.3 /
\plot 8 1.2 8 -0.3 /
\put{$\scriptstyle{+}$} at 0.2 1.1 \put{$\scriptstyle{-}$} at -0.2 1.1
\put{$1$} at 0.5 -0.3
\put{$x^{8\omega}$} at 8.5 -0.3
\setplotsymbol(.)
\arrow <6pt> [.2,.67] from 0.5 0.4 to 1.5 0.4
\arrow <6pt> [.2,.67] from 1.5 0.4 to 2.5 0.4
\arrow <6pt> [.2,.67] from 2.5 0.4 to 3.5 0.4
\arrow <6pt> [.2,.67] from 3.5 0.4 to 4.5 0.4
\arrow <6pt> [.2,.67] from 4.5 0.4 to 5.5 0.4
\arrow <6pt> [.2,.67] from 5.5 0.4 to 6.5 0.4
\arrow <6pt> [.2,.67] from 6.5 0.4 to 7.5 0.4
\arrow <6pt> [.2,.67] from 7.5 0.4 to 8.5 0.4
\endpicture
$$

Using Theorem~\ref{thm.EP-E}, we compute the terms in the expansion of $E_{-8\omega}P_{2\omega}$ arising from $h_1$ and $h_2$.  Both walks begin in $x^{2\omega}$, so $b_{(\ttb(h_j), x^{2\omega})} =1 $ for $j=1,2$. 
 
From $h_1$, 
$$\begin{array}{cccc}
c_1 = 1, &
c_2 = \displaystyle \frac{t^{-\frac12}(1-t)}{1-q^2t}, &
c_3 = \displaystyle \frac{t^{-\frac12}(1-t)q^3t}{1-q^3t}, &
c_4 = 1,\\
c_5 = 1, &
c_6 = 1, &
c_7 = \displaystyle \frac{t^{-\frac12}(1-t)}{1-q^7t}, &
c_8 = \displaystyle \frac{1-q^6}{1-q^6t} \frac{1-q^6t^2}{1-q^6t},
\end{array}$$
so $h_1$ gives rise to the term
$$\frac{1-t}{1-q^2t}
\frac{1-t}{1-q^3t}
\frac{1-t}{1-q^7t}
\frac{1-q^6}{1-q^6t} \frac{1-q^6t^2}{1-q^6t}
q^3t^{-\frac12}E_{6\omega}.$$

From $h_2$,
$$\begin{array}{cccc}
c_1 = 1, &
c_2 = \displaystyle \frac{t^{-\frac12}(1-t)}{1-q^2t}, &
c_3 = \displaystyle \frac{t^{-\frac12}(1-t)q^3t}{1-q^3t}, &
c_4 = \displaystyle \frac{t^{-\frac12}(1-t)q^4t}{1-q^4t},\\
c_5 = \displaystyle \frac{1-q^3}{1-q^3t} \frac{1-q^3t^2}{1-q^3t}, &
c_6 = \displaystyle \frac{1-q^2}{1-q^2t} \frac{1-q^2t^2}{1-q^2t}, &
c_7 = \displaystyle \frac{t^{-\frac12}(1-t)}{1-q^7t}, &
c_8 = \displaystyle \frac{t^{\frac12}(1-t)}{1-q^2t},
\end{array}$$
and $h_2$ has two negative grey folds, so it gives rise to the term
$$(-1)^2 \frac{1-t}{1-q^2t}
\frac{1-t}{1-q^3t}
\frac{1-t}{1-q^4t}
\frac{1-q^3}{1-q^3t} \frac{1-q^3t^2}{1-q^3t}
\frac{1-q^2}{1-q^2t} \frac{1-q^2t^2}{1-q^2t}
\frac{1-t}{1-q^7t}
\frac{1-t}{1-q^2t}
q^7t^{-\frac12}E_{2\omega}.$$
\hfill$\diamond$
}\end{example}


\subsection{Pieri formulas} \label{sec.pieri}
This section concerns the special cases of Theorems~\ref{thm.EP-E} and~\ref{thm.PP-P} when $\mu = \omega_j\in (\fh_\bbZ^*)^+$ is a minuscule weight.   Recall that $x^{\omega_j} = m_{\omega_j} v_{\omega_j}$.  In this section, we will use the notation $v_j = v_{\omega_j}$.
Also, in this case, the minimal length coset representative $m_{\omega_j} = \pi_j^\vee \in \Pi^\vee$.  
So
$$E_{\omega_j}\mathbf{1} 
= \pi_j^\vee \mathbf{1} 
= X^{\omega_j}T_{v_j^{-1}}\mathbf{1} 
= t_{v_j^{-1}}^{\frac12}X^{\omega_j}\mathbf{1}.
$$
Moreover, the walks appearing in Theorems~\ref{thm.EP-E} and~\ref{thm.PP-P} have type $(\pi_j^\vee)^{-1}$, which is a ``change in sheets''.  Such walks do not have crossings or foldings, so the product formulas simplify significantly. 
For $\lambda \in (\fh_\bbZ^*)^+$,
\begin{align}
\label{eqn.PieriEP} 
E_{\omega_j} P_\lambda
&= \sum_{v\in W^\lambda}
\sum_{h \in \Gamma\left((\pi_j^\vee)^{-1}, (vm_\lambda)^{-1}\right)}
	\!\!\!\!
	b_{(\ttb(h),x^{-w_0\lambda})} E_{\varpi(h)},\\
\label{eqn.PieriPP} 
P_{\omega_j} P_\lambda
&= \frac{1}{t_{w_{\omega_j}}^{-\frac12} W_{\omega_j}(t)}
	\sum_{v\in W^\lambda}
	\sum_{h \in \Gamma\left((\pi_j^\vee)^{-1}, (vm_\lambda)^{-1}\right)}
	\!\!\!\!
	b_{(\ttb(h),x^{-w_0\lambda})} 
	e_{(x^{\ttwt(h)}w_0, \tte(h))}
	P_{-w_0\ttwt(h)},
\end{align}
where each sum is over the set of alcove walks of type $(\pi_j^\vee)^{-1}$ beginning in $(vm_\lambda)^{-1}$ for $v\in W^\lambda$. 
Recall that the weight $\varpi(h)$ is defined by $\tte(h)^{-1} = m_{\varpi(h)}$, and
\begin{align}\label{eqn.bhqt}
b_{\left(\ttb(h),x^{-w_0\lambda}\right)} 
&= \prod_{a^\vee\in\calL\left(\ttb(h), x^{-w_0\lambda}\right)} 
	t_{a^\vee}^{\frac12} 
	\frac{1-q^{\shf(-a^\vee)}t^{\hgt(-a^\vee)} t_{a^\vee}^{-1}} 
	{1-q^{\shf(-a^\vee)}t^{\hgt(-a^\vee)}},\\
\label{eqn.ehqt}
e_{\left(x^{\ttwt(h)}w_0, \tte(h)\right)}
&= \prod_{a^\vee\in\calL\left(x^{\ttwt(h)}w_0, \tte(h)\right)} 
	t_{a^\vee}^{-\frac12} 
	\frac{1-q^{\shf(-a^\vee)}t^{\hgt(-a^\vee)} t_{a^\vee}} 
	{1-q^{\shf(-a^\vee)}t^{\hgt(-a^\vee)}}.
\end{align}
Given a walk $h$ beginning in $\ttb(h) = m_\lambda^{-1}v^{-1}$,
its endpoint is in $\tte(h) = m_\lambda^{-1} v^{-1} (\pi_j^\vee)^{-1}$, so $m_{\varpi(h)} = \tte(h)^{-1} = \pi_j^\vee v m_\lambda$ means that
$\varpi(h)=\pi_j^\vee v\lambda = x^{\omega_j}v_j^{-1}v\lambda =  v_j^{-1}v\lambda + \omega_j$. 
Thus we may also write
\begin{align*}
E_{\omega_j}P_\lambda
&= \sum_{v\in W^\lambda}  
	\left( \prod_{a^\vee\in m_\lambda^{-1}
	\calL\left(v^{-1}, v_\lambda^{-1}\right)} 
	t_{a^\vee}^{\frac12} 
	\frac{1-q^{\shf(-a^\vee)}t^{\hgt(-a^\vee)} t_{a^\vee}^{-1}} 
	{1-q^{\shf(-a^\vee)}t^{\hgt(-a^\vee)}}\right)
	E_{v_j^{-1}v\lambda+\omega_j}.
\end{align*}

\subsection{Compression of the Pieri formula}
Formula~\eqref{eqn.PieriPP} is a sum over $|W^\lambda|$ walks, and many of the walks have the same weight.  By imposing a condition on the final direction of the walks and modifying the coefficients appropriately, the formula can be compressed to contain the minimal number of terms.

\begin{corollary}\label{cor.PPshort} Let $\omega_j$ be a minuscule weight, and $\lambda$ be a dominant weight.  Then with the same notation as~\eqref{eqn.bhqt} and~\eqref{eqn.ehqt},
$$P_{\omega_j}P_\lambda 
= 
\sum_{{h\in \Gamma((\pi_j^\vee)^{-1},(zm_\lambda)^{-1})}  
	\atop {z\in W^\lambda,\ \ttd(h) \in W^{\omega_j}}}
	b_{\left(\ttb(h),x^{-w_0\lambda}\right)}  
	e_{\left(x^{\ttwt(h)}w_0,\tte(h)w_{\omega_j}\right)}  
	P_{-w_0\ttwt(h)}, $$
where the sum is over the set of alcove walks of type $(\pi_j^\vee)^{-1}$, beginning in $m_\lambda^{-1} z^{-1}$ for $z\in W^\lambda$, and has final direction $\ttd(h) \in W^{\omega_j}$.
\end{corollary}

\begin{remark}{\em $\ $
\begin{enumerate}
\item The {\em initial direction} of $h$ is defined by $X^{-w_0\lambda}T_{\tti(h)}$.  It follows from equation~\eqref{eqn.tteh} that
the final direction condition $\ttd(h) \in W^{\omega_j}$ is equivalent to the initial direction condition $\tti(h) \in W^{-w_0\omega_j}$.  
\item Since the walks which are under consideration do not have folds, then the condition that the walks begin in the alcoves $m_\lambda^{-1}z^{-1}$ for $z \in W^\lambda$ is enough to guarantee that the walks are contained in the closure of the dominant chamber.
\end{enumerate}
}\end{remark}

\proof
If a walk $h$ begins in the alcove 
\begin{equation}\label{eqn.ttbh}
\ttb(h) = m_\lambda^{-1}z^{-1} = m_{-w_0\lambda} z^{-1} 
= x^{-w_0\lambda} \left(v_\lambda z^{-1}\right),
\end{equation}
for some $z\in W^\lambda$, then $h$ has initial direction $\tti(h) = v_\lambda z^{-1}$.  And since $h$ is a walk of type $(\pi_j^\vee)^{-1}$, then $h$ ends in the alcove
\begin{equation}\label{eqn.tteh}
\tte(h) = \ttb(h)(\pi_j^\vee)^{-1} = x^{-w_0\lambda}\tti(h) (x^{\omega_j}v_j^{-1})^{-1} = x^{-w_0\lambda} x^{-\tti(h)v_j\omega_j} \left(\tti(h) v_j\right),
\end{equation}
so $h$ has final direction $\ttd(h) = \tti(h) v_j$, and weight 
\begin{equation}\label{eqn.ttwth}
\ttwt(h) 
= -w_0\left(\lambda + w_0^{-1}\tti(h) v_j\omega_j\right)
= -w_0\left(\lambda + w_0^{-1}\ttd(h)\omega_j\right).
\end{equation}

The weights of the walks $h\in \Gamma((\pi_j^\vee)^{-1}, (vm_\lambda)^{-1})$ are not distinct, since the stabilizer of $\omega_j$ is not trivial.  The idea is to group together the walks with the same weight by factoring $\tti(h) = vw$ for $v\in W^{-w_0\omega_j}$ and $w\in W_{-w_0\omega_j}$.

We make one more observation before proceeding with the calculation.  Suppose $h$ is a walk beginning in an alcove $x^{-w_0\lambda} u$, for some $u\in W_0$, which is not contained in the dominant chamber, and that the hyperplane $\ttH_{-\alpha_k^\vee}$ separates the alcove from the dominant chamber.  Notice that one of the factors of the coefficient 
\begin{equation}\label{eqn.startoutsideC}
b_{\left(\ttb(h),x^{-w_0\lambda}\right)} 
= \prod_{a^\vee\in\calL\left(\ttb(h), x^{-w_0\lambda}\right)} 
	t_{a^\vee}^{\frac12} 
	\frac{1-q^{\shf(-a^\vee)}t^{\hgt(-a^\vee)} t_{a^\vee}^{-1}} 
	{1-q^{\shf(-a^\vee)}t^{\hgt(-a^\vee)}}, 
\end{equation}
is $\displaystyle t_k^{\frac12} \frac{1-t_k t_k^{-1}}{1-t_k} =0$, so that $b_{\left(\ttb(h),x^{-w_0\lambda}\right)} =0$ in this case.  Therefore, although such walks are not counted in equation~\eqref{eqn.PieriPP}, it will still make sense to include them in the following calculation.

Consider the bijections
$$\bar{\quad} : W_{-w_0{\omega_j}} \rightarrow  W_{\omega_j}: w \mapsto \bar{w}= v_j^{-1} w v_j,
\quad\hbox{and}\quad
\bar{\quad} :W^{-w_0{\omega_j}} \rightarrow  W^{\omega_j} : v \mapsto  \bar{v}=vv_j.$$
With this notation, then $\tti(h) = vw$ if and only if $\ttd(h) = \bar{v}\bar{w}$.
By equations~\eqref{eqn.ttbh} and~\eqref{eqn.tteh}, the walks $h$ have the same weight only if
$$\ttb(h) \in \calB_v =
\left\{x^{-w_0\lambda}vw \mid w\in W_{-w_0\omega_j}\right\},
$$
for a fixed $v\in W^{-w_0\omega_j}$, or equivalently, if
$$\tte(h) \in \calE_{\bar{v}} = 
\left\{x^{\ttwt(h)}\bar{v}\bar{w} \mid \bar{w}\in W_{\omega_j}\right\},
$$
for a fixed $\bar{v}\in W^{\omega_j}$.  
Since $b_{(x^{-w_0\lambda}vw, x^{-w_0\lambda}v)} = b_{(x^{\ttwt(h)}\bar{v}\bar{w}, x^{\ttwt(h)}\bar{v})}$, then 
\begin{align*}
\sum_{h: \tte(h) \in \calE_{\bar{v}}} &
	b_{\left(\ttb(h),x^{-w_0\lambda}\right)} 
	e_{\left(x^{\ttwt(h)}w_0, \tte(h)\right)} \mathbf{1} \\
&= b_{\left(x^{-w_0\lambda}v,x^{-w_0\lambda}\right)} 
	e_{\left(x^{\ttwt(h)}w_0, x^{\ttwt(h)}\bar{v} w_{\omega_j}\right)}
	\sum_{\bar{w}\in W_{\omega_j}} 
	b_{\left(x^{\ttwt(h)}\bar{v}\bar{w}, x^{\ttwt(h)}\bar{v}\right)}
	e_{\left(x^{\ttwt(h)}\bar{v} w_{\omega_j}, x^{\ttwt(h)}\bar{v}\bar{w}
	\right)} \mathbf{1}\\
&= b_{\left(x^{-w_0\lambda}v,x^{-w_0\lambda}\right)} 
	e_{\left(x^{\ttwt(h)}w_0, x^{\ttwt(h)}\bar{v} w_{\omega_j}\right)}
	\left( x^{\ttwt(h)}\bar{v}  \cdot
    \sum_{\bar{w}\in W_{\omega_j}} 
	b_{(\bar{w}, 1)}	e_{(w_{\omega_j}, \bar{w})} \right)\mathbf{1}\\
&= b_{\left(x^{-w_0\lambda}v,x^{-w_0\lambda}\right)} 
	e_{\left(x^{\ttwt(h)}w_0, x^{\ttwt(h)}\bar{v} w_{\omega_j}\right)}
	t_{w_{\omega_j}}^{-\frac12} W_{\omega_j}(t)\mathbf{1},
\end{align*}
by Lemma~\ref{lem.Poincare}.
Therefore, we can restrict to walks which have final direction $\ttd(h) = \bar{v} \in W^{\omega_j}$,
$$
P_{\omega_j}P_\lambda
= 
\sum_{h: \ttd(h) \in W^{\omega_j}} 
 b_{\left(\ttb(h),x^{-w_0\lambda}\right)} 
	e_{\left(x^{\ttwt(h)}w_0, \tte(h) w_{\omega_j}\right)}
	P_{-w_0\ttwt(h)}.$$
\qed

\subsection{Type $A$ Pieri formulas and partitions}

In the case of the type $A_n$ root systems, Macdonald gave Pieri formulas for symmetric Macdonald polynomials in terms of partitions.
This section is a brief sketch of the relation between the alcove walk combinatorics in Corollary~\ref{cor.PPshort} and the partition combinatorics in Macdonald's formula~\cite[(6.24)(iv)]{M88}, reproduced below.

In the type $A_n$ setting, the parameters $t=t_0=\cdots = t_n$ are necessarily all equal.  The correspondence between partitions with at most $n+1$ parts and type $A_n$ dominant weights is as follows:
\begin{eqnarray*}
\hbox{column partitions $(1^j)$}
	&\leftrightarrow& 
	\hbox{minuscule weights $\omega_j$ for $j=1,\ldots,n$,}\\
\hbox{column partition $(1^{n+1})$} 
	&\leftrightarrow& 
	\hbox{weight $0$,}
\end{eqnarray*} 
so
$$e_j = s_{(1^j)} = P_{\omega_j}(q,t) \hbox{ is the $j$th elementary symmetric polynomial.}$$  
Viewing a partition $\kappa$ as a collection of boxes justified to the top and the left,
the {\em arm-length} and {\em leg-length} of a box $s$ in the $i$th row and $j$th column of $\kappa$ are
\begin{align*}
a_\kappa(s) &= \kappa_i-j = \textrm{ number of boxes to the east of } s,\\
l_\kappa(s) &= \kappa_j'-i = \textrm{ number of boxes to the south of } s.
\end{align*}
If $\kappa \supseteq \lambda$, let $C_{\kappa-\lambda}$ (respectively $R_{\kappa-\lambda}$) be the set of columns (respectively rows) of $\kappa$ that intersect the skew partition $\kappa-\lambda$.

\begin{theorem} \cite[(6.24)(iv)]{M88} \label{eqn.MacPP}
$$P_{\omega_j}P_\lambda = \sum_\kappa 
\left( 
\prod_{s\in C_{\kappa-\lambda} - R_{\kappa-\lambda}} 
	\frac{1-q^{a_\lambda(s)+1}t^{l_\lambda(s)}}
	{1-q^{a_\kappa(s)+1}t^{l_\kappa(s)}}
	\frac{1-q^{a_\kappa(s)}t^{l_\kappa(s)+1}}
	{1-q^{a_\lambda(s)}t^{l_\lambda(s)+1}}
	\right)
	P_\kappa, 
$$
where the sum is over the set of partitions $\kappa$ with at most $n+1$ parts, obtained from $\lambda$ by adding $j$ boxes, with no two in the same row. \qed
\end{theorem}

To draw the connection between walks and partitions, we begin by considering Corollary~\ref{cor.PPshort} with parameters $t_i=t$:

{\em \noindent{\bf Corollary 4.8.}
\begin{align*}
&P_{\omega_j}P_\lambda\\
&=\sum_h \left(
	\prod_{\alpha^\vee \in \calL(\tti(h))}
	\left(x^\lambda\!\cdot\!
	\frac{1-q^{\shf(-\alpha^\vee)}t^{\hgt(-\alpha^\vee)-1} } 
	{1-q^{\shf(-\alpha^\vee)}t^{\hgt(-\alpha^\vee)}}\right)
\!\!
\left(x^{\ttwt(h)}\!\cdot\!	
	\frac{1-q^{\shf(-\alpha^\vee)}t^{\hgt(-\alpha^\vee)+1} } 
	{1-q^{\shf(-\alpha^\vee)}t^{\hgt(-\alpha^\vee)}}\right)
 \right) P_{\ttwt(h)},
\end{align*}
where the sum is over the set of alcove walks of type $\pi_j^\vee$, beginning in the alcoves $x^{\lambda}W^{\omega_j}$.
}
\proof
First, notice that
$\calL\left(\ttb(h),x^{-w_0\lambda}\right)
= x^{-w_0\lambda}\calL\left(1,\tti(h)\right)$, while
\begin{align*}
\calL\left(x^{\ttwt(h)}w_0, \tte(h)w_{\omega_j}\right)
&=\calL\left(x^{\ttwt(h)}, \tte(h) w_{\omega_j}w_0^{-1}\right)
=x^{\ttwt(h)} \calL\left(1, \ttd(h)v_j^{-1}\right)
=x^{\ttwt(h)} \calL\left(1, \tti(h)\right),
\end{align*}
so the coefficients in Corollary~\ref{cor.PPshort} can be rewritten as
\begin{align*}
b_{\left(\ttb(h),x^{-w_0\lambda}\right)} 
	&e_{\left(x^{\ttwt(h)}w_0, \tte(h) w_{\omega_j}\right)}\mathbf{1}
= \prod_{\alpha^\vee\in \calL(1,\tti(h))} 
\left(x^{-w_0\lambda}\cdot
	\frac{1-Y^{-\alpha^\vee}t^{-1} } {1-Y^{-\alpha^\vee}} \right)
\left(x^{\ttwt(h)}\cdot	
	\frac{1-Y^{-\alpha^\vee}t } 
	{1-Y^{-\alpha^\vee}}\right)\mathbf{1},
\end{align*}
where $x^{-w_0\lambda} = x^{\ttwt(h) + \ttd(h)\omega_j}$ by~\eqref{eqn.ttwth}.

Next, consider the action of the automorphism $-w_0:\fh_\bbR^* \rightarrow \fh_\bbR^*$, which acts on the (type $A_n$) lattice $\fh_\bbZ^*$ by sending $\omega_j \mapsto \omega_{n+1-j}$.  Hence, this automorphism also acts on the set of alcoves, with action given by
$$-w_0(x^\beta w\ttA) = x^{-w_0\beta} w_0^{-1}ww_0\ttA, \qquad
\hbox{for } \beta\in \fh_\bbZ^*, w\in W_0.$$
Since $\lambda$ is a dominant weight, then
$$(zm_\lambda)^{-1} = m_\lambda^{-1}z^{-1} = m_{-w_0\lambda}z^{-1}
= x^{-w_0\lambda}v_{-w_0\lambda}^{-1} z^{-1} = x^{-w_0\lambda} (zv_{-w_0\lambda})^{-1},$$
where $zv_{-w_0\lambda}\in W^{-w_0\omega_j}$.  The Pieri formula stated above is the result of applying $-w_0$ to the set of alcoves.

The last observation is, while it seems that the condition that walks must be contained in the dominant chamber is dropped in this version of Corollary~\ref{cor.PPshort}, in fact if the walk $h$ begins outside the the dominant chamber, then the coefficient $b_h =0$ for the same reason as explained in~\eqref{eqn.startoutsideC}.  So the Pieri formula remains unchanged even if the set of walks on which the sum runs over is expanded to include those which begin outside the dominant chamber.
\qed

In what follows, we fix $\lambda\in (\fh_\bbZ^*)^+$ to be a regular dominant weight, so that when viewed as a partition, $\lambda$ does not have equal parts.  Also fix a minuscule weight $\omega_j$.

Let $\calK$ be the lattice of partitions which are obtained from $\lambda$ by adding $j$ boxes with no two in the same row.  Given $\kappa\in \calK$, any box not contained in $\lambda$ will be called an {\em added box}.  Since we assume that  $\lambda$ has no equal parts, $\calK$ may be viewed as the set of $j$-subsets of $[n+1]=\{1,\ldots, n+1\}$, identifying $\kappa$ with the subset $\{r_1,\ldots, r_j\}$ if the $j$ boxes were added in rows $r_1,\ldots, r_j$ of $\lambda$ to obtain $\kappa$.  The order in $\calK$ is defined by the covering relation:
$\kappa$ covers $\theta$ if and only if
$\theta$ can be obtained from $\kappa$ by moving one of the added boxes down one row.  Thus, the unique maximal element $\hat{1}$ in $\calK$ is $\lambda$ with boxes added in rows $1$ through $j$, and the unique minimal element $\hat{0}$ in $\calK$ is $\lambda$ with boxes added in rows $(n+1)-(j-1)$ through $n+1$.

Further, in the Hasse diagram representing the lattice $\calK$, if $\kappa$ covers $\theta$ and a box was moved from row $i$ to row $i+1$, then the edge between $\kappa$ and $\theta$ is labeled $s_i$.

Recall $\Gamma(\pi_j^\vee, x^\lambda W^{\omega_j})$ is the set of alcove walks of type $\pi_j^\vee$ which begin in the alcoves $x^\lambda W^{\omega_j}$.
\begin{lemma}  Assume $\lambda$ is a regular dominant weight, or equivalently, a partition with no equal parts. 
There is a bijection 
\begin{eqnarray*}
\calK &\longleftrightarrow & \Gamma\left(\pi_j^\vee, x^\lambda W^{\omega_j}\right) \\
\kappa &\leftrightarrow& h_\kappa
\end{eqnarray*}
such that
\begin{enumerate}
\item $h_\kappa$ is the walk with $\ttwt(h_\kappa)
= \kappa$,
\item if a shortest path from $\kappa$ to $\hat{1}$ in $\calK$ traverses the sequence of edges $s_{i_1},\ldots, s_{i_\ell}$, then $\tti(h) = s_{i_1}\cdots s_{i_\ell}$.
\item Moreover, for each pair $\kappa \leftrightarrow h_\kappa$, there is a bijection
\begin{eqnarray*}
C_{\kappa-\lambda}-R_{\kappa-\lambda} & \longleftrightarrow & \calL(\tti(h_\kappa))\\
\Box_{\alpha^\vee} & \leftrightarrow & \alpha^\vee
\end{eqnarray*} 
satisfying
$$
l_\kappa(\Box_{\alpha^\vee}) = \hgt(-x^{\kappa}\alpha^\vee),
\qquad\hbox{and}\qquad
a_\kappa(\Box_{\alpha^\vee}) = \shf(-x^{\kappa}\alpha^\vee).
$$
\end{enumerate}
\end{lemma}
\proof

\noindent (1)
Since $\lambda$ is a regular weight, then no matter which $j$ rows we choose to add boxes, the resulting shape corresponds to a partition, hence the lattice $\calK$ has $\binom{n+1}{j}$ elements.  In the type $A_n$ root system, $W_0 = \fS_{n+1}$ and $W_{\omega_j} \cong \fS_j \times \fS_{n+1-j}$, thus $|W^{\omega_j}| = (n+1)!/j!(n+1-j)! = \binom{n+1}{j}.$  So the sets in question have the same cardinality.  To each partition $\kappa\in \calK$, we assign the walk $h_\kappa$ such that $\ttwt(h_\kappa) = \kappa$.

\noindent (2)
For $i=1,\ldots, n+1$, adding a box in row $i$ of the partition $\lambda$ yields the composition corresponding to the weight $\lambda + \omega_i-\omega_{i-1}$, where we use the notation $\omega_0 = \omega_{n+1}=0$.  Hence, adding a box in rows $1$ through $i$ of $\lambda$ yields the partition $\lambda + \omega_i$.

If $\kappa$ is the partition obtained from $\lambda$ by adding $j$ boxes in rows $r_1,\ldots, r_j$ (assuming $r_1< r_2< \cdots < r_j$), then $\kappa = \lambda + \sum_{i=1}^j (\omega_{r_i}-\omega_{r_i-1})$.  Observe that $\hat{1} = \lambda+\omega_j$, while $\hat{0}=\lambda+v_j\omega_j$. 
A minimal path in $\calK$ from $\hat{1}$ to $\kappa$ can be constructed as follows: begin by traveling along the edges $s_j, s_{j+1}, \ldots, s_{r_j-2}, s_{r_j-1}$ to reach the partition in which the last box in row $j$ of $\hat{1}$ has been moved to row $r_j$, then travel along the edges $s_{j-1}, s_{j}, \ldots, s_{r_{j-1}-1}$ to reach the partition in which the last box in row $j-1$ of the previous partition has been moved to row $r_{j-1}$, and so on, and end by traveling along the edges $s_1, s_2, \ldots, s_{r_1-1}$ to reach the partition $\kappa$.

Inverting these $j$ sequences of edges gives a path from $\kappa$ to $\hat{1}$.  Let
$$\sigma_k = s_{r_k-1} s_{r_k-2} \cdots s_{k+1}s_k, 
\hbox{ for } k=1,\ldots, j,
\qquad
\hbox{ and let } 
\sigma = \sigma_1 \sigma_2 \cdots \sigma_j.$$
Then
$\lambda+ \tti(h)\omega_j 
= \ttwt(h_\kappa) 
= \kappa 
= \lambda + \sigma\omega_j.$
Since this is a path of minimal length in $\calK$, then the above factorization is a reduced expression for $\tti(h) = \sigma \in W^{\omega_j}$.

\noindent (3)
By equation~\eqref{calLw}, 
$$
\begin{array}{l}
\calL(\tti(h_\kappa)) 
= \Big\{
\alpha_{r_1-1}^\vee,\
	(s_{r_1-1}\alpha_{r_1-2}^\vee),\ \ldots,\ 
	(s_{r_1-1} s_{r_1-2} \cdots s_2\alpha_1^\vee),\\ \\
\qquad\qquad\qquad\qquad
(\sigma_1\alpha_{r_2-1}^\vee),\ 
	(\sigma_1 s_{r_2-1}\alpha_{r_2-2}^\vee),\ \ldots,\ 
	(\sigma_1 s_{r_2-1} s_{r_2-2} \cdots s_3\alpha_2^\vee),\\ \\
\qquad\qquad\qquad\qquad
\qquad \qquad\qquad\qquad\qquad \vdots \\ \\
\qquad\qquad\qquad\qquad
(\sigma_1\cdots \sigma_{j-1}\alpha_{r_j-1}^\vee),\
	\ldots,\ 
	(\sigma_1\cdots \sigma_{j-1}s_{r_j-1} s_{r_j-2} \cdots s_{j+1}\alpha_j^\vee)\Big\}.\end{array}
$$ 
In contrast, $C_{\kappa-\lambda}- R_{\kappa-\lambda}$ is the set of boxes in $\kappa$ which are in columns containing an added box, and in rows not containing an added box.  For example, when $\kappa=\hat{1}$, $C_{\kappa-\lambda}- R_{\kappa-\lambda} = \emptyset$.  

Consider the box that was added to $\lambda$ in row $r_k$ (and column $c_k$).  Column $c_k$ contains $(r_k-1) - (k-1)$ boxes which are in $C_{\kappa-\lambda}- R_{\kappa-\lambda}$.   Altogether, 
$$|C_{\kappa-\lambda}- R_{\kappa-\lambda}| = \sum_{k=1}^j (r_k-k) = \sum_{k=1}^j\ell(\sigma_k) = \ell(\sigma) = \ell(\tti(h_\kappa))
= |\calL(\tti(h_\kappa))|,$$
noting that column $c_k$ of $\kappa$ has $\ell(\sigma_k)$ boxes.  Thus, define the correspondence $C_{\kappa-\lambda}- R_{\kappa-\lambda} \leftrightarrow \calL(\tti(h_\kappa))$ as follows: to the boxes in column $c_k$ which belong to $C_{\kappa-\lambda}- R_{\kappa-\lambda}$, assign in order from bottom to top, the coroots
$$(\sigma_1 \cdots \sigma_{k-1} \alpha_{r_k-1}^\vee),\ 
	(\sigma_1 \cdots \sigma_{k-1}s_{r_k-1}\alpha_{r_k-2}^\vee),\ \ldots,\ 
	(\sigma_1 \cdots \sigma_{k-1}
	s_{r_k-1} s_{r_k-2} \cdots s_{k+1}\alpha_k^\vee).$$
Using this particular choice of factorization for $\tti(h_\kappa)$, and noting that in the type $A_n$ root system
$$s_i\alpha_k^\vee 
=\begin{cases} 
\alpha_k^\vee + \alpha_i^\vee, &\hbox{if } i=k\pm1,\\
-\alpha_k^\vee, &\hbox{if } i=k,\\
\alpha_k^\vee, &\hbox{otherwise,} 
\end{cases}$$
we see that if the coroot 
$$\gamma^\vee 
= \left(\sigma_1\cdots \sigma_{k-1}\right)\left( s_{r_k-1}s_{r_k-2}\cdots s_{r_k-m+1}\right)\alpha_{r_k-m}^\vee 
$$ 
is the one which corresponds to the box $\Box_{\gamma^\vee}$ in column $c_k$ and row $t$ of $\kappa$ under the above bijection, then
$$\gamma^\vee 
= \alpha_{r_k-1}^\vee + \alpha_{r_k-2}^\vee + \cdots + \alpha_{r_k-m}^\vee + \cdots + \alpha_{t+1}^\vee +\alpha_{t}^\vee.$$
Hence it follows that
\begin{align*}
\hgt(-x^{\kappa}\gamma^\vee) 
	&= \hgt(\gamma^\vee)
	= r_k-t = l_\kappa(\Box_{\gamma^\vee}),\\
\shf(-x^{\kappa}\gamma^\vee)
	&= \left\langle \kappa, \gamma^\vee \right\rangle
	= \sum_{i=1}^{n} d_i\left\langle \omega_i, \gamma^\vee  \right\rangle
	= a_\kappa(\Box_{\gamma^\vee}),
\end{align*}
where $d_i$ in $\kappa = \sum_{i=1}^n d_i\omega_i$ is the number of columns in $\kappa$ with height $i$.
\qed

Under the assumption that $\lambda$ has no equal parts, then $a_\kappa(s) = a_\lambda(s)$ and $l_\kappa(s) = l_\lambda(s)+1$, for all $\kappa \in\calK$.  On the other hand,  $\hgt(-x^\lambda\alpha^\vee) = \hgt(\alpha^\vee) = \hgt(-x^{\ttwt(h)}\alpha^\vee)$.  Also, 
$\shf(-x^{\ttwt(h)}\alpha^\vee)
=\langle \ttwt(h), \alpha^\vee \rangle
= \langle \lambda + \tti(h)^{-1}\omega_j, \alpha^\vee \rangle
= \shf(-x^\lambda\alpha^\vee) + \langle \omega_j, \tti(h)\alpha^\vee \rangle$ for all $\alpha^\vee \in \calL(\tti(h))$, and since $\tti(h)^{-1}\alpha^\vee$ is of the form $s_j s_{i_\ell-1}\cdots s_{i_\ell-m+1}\alpha_{i_\ell-m}^\vee$, then $\langle \omega_j, \tti(h)\alpha^\vee\rangle = -1$, and so $\shf(-x^{\ttwt(h)}\alpha^\vee) = \shf(-x^\lambda\alpha^\vee)-1$.

It is now clear that Theorem~\ref{eqn.MacPP} and Corollary~\ref{cor.PPshort} are identical, in the case that $\lambda$ is a regular dominant weight (a partition with no equal parts).


\begin{remark}
{\em \ 
\begin{enumerate}
\item In the case when $\lambda$ is not a regular dominant weight (a partition having equal parts), then $|\calL(\tti(h))| \geq C_{\kappa-\lambda}-R_{\kappa-\lambda}$, and there are many canceling pairs of factors in the alcove walk formula which do not appear in the tableau formula.  Note that $a_\kappa(s) = a_\lambda(s)$, but $l_\kappa(s) \geq l_\lambda(s)+1$ in this case. 
\item The Pieri case involving multiplication with $P_{j\omega_1}$, where $P_{j\omega_1}(q,q) = s_{(j)} = h_j$ is the $j$th complete symmetric polynomial, is more involved, and is likely related to the compression phenomenon discussed in~\cite{Le}.
\end{enumerate}}\end{remark}

\subsection{Hall-Littlewood polynomials} \label{sec.specialization}
The Hall-Littlewood polynomials $P_\mu(t)=P_\mu(0,t)$ are the symmetric Macdonald polynomials under the specialization of the parameter $q=0$.  In this section, we assume all parameters $t_i = t$ for $i=0,\ldots, n$.

Combinatorial formulas for Hall-Littlewood polynomials were given in 
~\cite[Theorem 1.1]{Sc06}, and a Pieri formula was given in~\cite[p.\!\! 217]{M88}. 
We restate the alcove walk formula \cite[Theorem 4.9]{R06} for the Littlewood-Richardson rule for $P_\lambda(t)$ in the notation of this paper.  Recall that the final direction $\ttd(h)$ is defined by $X^{\tte(h)} = X^{\ttwt(h)}T_{\ttd(h)}$, and the initial direction $\tti(h)$ is defined by $X^{\ttb(h)} = X^{\lambda} T_{\tti(h)}$. 
Let $\ttf(h)$ be the number of folds in $h$, and $\ttf_0(h)$ be the number of folds in $h$ touching a wall of $\ttC$. 
A {\em positively folded walk} is a walk with no negative folds.

\begin{theorem}\cite[Theorem 4.9]{R06}. Also see~\cite[Theorem 1.3]{Sc06}.
\label{thm.HL1}

Let $\mu,\lambda \in (\fh_\bbZ^*)^+$, and fix a reduced expression for $x^\mu$.  Then
$$P_\mu(t) P_\lambda(t)
= \sum_{v\in W^\mu} \sum_{h\in \Gamma^{\ttC-\rho}_+
	\left(\vec{x}^{\mu}, x^\lambda v\right)} 
	t^{\frac12(\ell(\tti(h)) + \ell(\ttd(h)) - \ttf(h))} 
	(1-t)^{\ttf(h)-\ttf_0(h)} 
	P_{\ttwt(h)}(t),
$$
where the sum is over all positively folded alcove walks of type $\vec{x}^\mu$ beginning in $x^\lambda v$ for $v\in W^\mu$, contained in the shifted dominant chamber $\ttC-\rho$. \qed
\end{theorem}

As a corollary to Theorem~\ref{thm.PP-P}, we derive a different version of the above result.  The main differences are that the walks under consideration are of type $\vec{m}_\mu$ instead of $\vec{x}^\mu$, and the walks are contained in $\overline{\ttC}$ instead of $\ttC$ shifted by $-\rho$. 
\begin{corollary}\label{cor.HL2}

Let $\mu,\lambda \in (\fh_\bbZ^*)^+$, and fix a reduced expression for $m_\mu$.  Then
$$P_\mu(t) P_\lambda(t)
= 	t^{-\frac12\ell(v_\mu)} 
	\sum_{v\in W^\lambda}  \sum_{h \in \Gamma^{\overline{\ttC}}_+
	\left(\vec{m}_\mu, x^\lambda v^{-1}\right)} 
	t^{\frac12( \ell(\tti(h))+\ell(\ttd(h))-\ttf(h))}
	(1-t)^{\ttf(h)-\ttf_0(h)}
	\frac{W_{\ttwt(h)}(t)}{W_\mu(t)}
	P_{\ttwt(h)}(t),$$
where the sum is over all positively folded alcove walks of type $\vec{m}_\mu$ beginning in $x^\lambda v^{-1}$ for $v\in W^\lambda$, contained in $\overline\ttC$. 

\end{corollary}
\proof
Keeping the notation from Theorem~\ref{thm.PP-P}, let 
$$a_h(q,t) = (-1)^{\phi_{grey}^-(h)} 
	b_h	e_h
	\prod_{k=1}^r c_{k(h)}.
$$
The first step is to show that the walks $h\in \Gamma_2^\ttC(\vec{m}_\mu^{-1}, (vm_\lambda)^{-1})$ for which $a_h(0,t) \neq 0$ are precisely those whose folds must be positive and grey.  In this section, these are referred to as the `grey positively folded walks'.  

The monomial expansion of $E_\mu$ is a sum over the set of walks $\Gamma(\vec{m}_\mu)$.  There is a unique walk $l\in\Gamma(\vec{m}_\mu)$ without folds, and since $\mu$ is a dominant weight, then every step of $l$ is a positive crossing because $l$ is contained in the dominant chamber.  The walks in $\Gamma(\vec{m}_\mu)$ are generated by folding $l$ in all possible ways, so if a walk $p\in \Gamma(\vec{m}_\mu)$ has a fold, then it has at least one negative fold.  At $q=0$, Theorem~\ref{thm.Xtau-X} gives
$$E_\mu(0,t) = t_{v_\mu^{-1}}^{\frac12} X^\mu,  
\quad\hbox{if $\mu$ is dominant}.$$
In other words, if $h\in \Gamma_2^\ttC(\vec{m}_\mu^{-1}, (vm_\lambda)^{-1})$ has a black fold, then $p(h)$ has at least one negative black fold, so the last row of~\eqref{eqn.verylonglist} gives
$a_h(0,t)=0$ in this case.  Thus, any walk with a black fold does not contribute to the expansion of $P_\mu(t)P_\lambda(t)$ in the basis of Hall-Littlewood polynomials. 

Another way to put this is if $h$ survives setting $q=0$, then $p(h)$ is the unfolded walk of type $\vec{m}_\mu^{-1}$ ending in $1$.  Accordingly, $\epsilon_k(h) = -1$ for $k=1,\ldots, r$, since $\mu$ is a dominant weight.  Using this observation, if $h$ has a negative grey fold (so the fold must be against an affine wall), then the fourth row of~\eqref{eqn.verylonglist} gives $a_h(0,t)=0$ in this case.  

Since the walks now under consideration do not have negative grey folds, then the expansion of the product of Hall-Littlewood polynomials does not involve negative signs. Moreover, every fold is grey, so we may forget the colour of the folds for the remainder of this section.
The remaining coefficients in Theorem~\ref{thm.PP-P} simplify:
\begin{align*}
b_h &= t^{\frac12(\ell(x^{-w_0\lambda}) - \ell(\ttb(h)) )} 
	= t^{\frac12 \ell(\tti(h))},\\
e_h &= t^{-\frac12(\ell(\tte(h)) - \ell(m_{\ttwt(h)}) )}
	= t^{\frac12(\ell(\ttd(h)) - \ell(v_{\ttwt(h)}) )},\\
c_{k(h)} &= \begin{cases}
	t^{-1/2}, & \hbox{for (grey) folds touching a wall of }\ttC,\\
	t^{-1/2}(1-t), & 
		\hbox{for (grey) folds touching an affine hyperplane.}
	\end{cases}	\\
		&= t^{-\frac12 \ttf(h)} (1-t)^{\ttf(h)-\ttf_0(h)}.
\end{align*}

One last observation to make is that the action of the automorphism $-w_0: \fh_\bbR^* \rightarrow \fh_\bbR^*$ acts on the lattice $\fh_\bbZ^*$ by permuting $\omega_1,\ldots, \omega_n$, so $-w_0$ acts on the set of alcoves, given by 
$$-w_0(x^\beta w \ttA) =  x^{-w_0\beta} w_0^{-1} w w_0\ttA, \qquad\hbox{for }
\beta\in \fh_\bbZ^*, w\in W_0.$$
Since $\lambda$ is a dominant weight, then
$$(vm_\lambda)^{-1} = m_\lambda^{-1}v^{-1} = m_{-w_0\lambda} v^{-1}
= x^{-w_0\lambda} v_{-w_0\lambda}^{-1} v^{-1} = x^{-w_0\lambda} v_\lambda v^{-1}. $$
Therefore, if $h$ is a walk of type $\vec{m}_\mu^{-1} = \vec{m}_{-w_0\mu}$ beginning in $(vm_\lambda)^{-1}$, then $-w_0(h)$ is a walk of type $\vec{m}_\mu$ beginning in $x^\lambda w_0^{-1} v_\lambda v^{-1}w_0.$
Moreover, $\{v_\lambda v^{-1} \mid v\in W^\lambda\} = \{y^{-1} \mid y\in  W^{-w_0\lambda} \}$, so 
$\{ w_0 v_\lambda v^{-1} w_0^{-1} \mid v\in W^\lambda \} 
= \{ v^{-1} \mid v\in W^\lambda\}$.

Putting all the observations together, 
$$P_\mu(t) P_\lambda(t)
= \sum_{v\in W^\lambda}  \sum_{h \in \Gamma_+^{\overline{\ttC}}
	\left(\vec{m}_\mu, x^\lambda v^{-1}\right)} 
	t^{\frac12(\tti(h)+\ttd(h)-\ttf(h))}
	(1-t)^{\ttf(h)-\ttf_0(h)}
	\frac{ t^{-\frac12}_{v_{\ttwt(h)}}
	t_{w_{\ttwt(h)}}^{-\frac12}W_{\ttwt(h)}(t)}
		{t_{w_\mu}^{-\frac12}W_\mu(t)}
	P_{\ttwt(h)}(t),$$
where the sum is over the set of positively folded walks of type $\vec{m}_\mu$, beginning in $x^\lambda v^{-1}$ for $v\in W^\lambda$, which are contained in $\overline{\ttC}$.

To finish, note that 
$t^{-\frac12}_{v_{\ttwt(h)}}t_{w_{\ttwt(h)}}^{-\frac12}t_{w_\mu}^{\frac12} 
= t_{w_0}^{-\frac12} t_{w_\mu}^{\frac12} 
= t_{v_\mu}^{-\frac12}.$
\qed

%


%

\begin{example}{\em See section~\ref{sec.eg} for more details on the alcove picture of type $\fsl_3$.

Let $\varphi =\omega_1+\omega_2$.  Using either Theorem~\ref{thm.HL1} or Corollary~\ref{cor.HL2},
\begin{align*}
P_\varphi(t)^2 = P_{2\varphi}(t) &+ (1+t) P_{3\omega_1}(t) + (1+t)P_{3\omega_2}(t)\\ 
&+ (2 + t(1-t)) P_\varphi(t) + (1+2t+2t^2+t^3)P_0(t).
\end{align*}
We remark that 16 walks were used in the first formula, and seven walks were used in the second formula.
$\hfill\diamond$
}\end{example}

\section{Mainly Lots of Examples}\label{sec.eg}
\subsection{Type $\fsl_2$}
The root datum is
$$\fh_\bbZ^* =\bbZ\omega,\quad 
R = \{\pm\alpha\},\quad 
\fh_\bbZ = \bbZ\omega^\vee,\quad 
R^\vee = \{\pm\alpha^\vee\},$$
with $\alpha = 2\omega$ and $\alpha^\vee = 2\omega^\vee$.

Let $s_0^\vee = x^\alpha s_1^\vee$ and $\pi^\vee = x^\omega s_1^\vee$.
The extended affine Weyl group $W^\vee$ is generated by $s_1^\vee, s_0^\vee, \pi^\vee$, subject to the relations 
\begin{equation}\label{eqn.presn1}
(\pi^\vee)^2 =1,\qquad
\pi^\vee s_0^\vee = s_1^\vee \pi^\vee,\qquad
(s_i^\vee)^2 =1,\hbox{ for } i=0,1.
\end{equation} 
Alternatively, 
\begin{equation}\label{eqn.presn2}
W^\vee = \{x^{k\omega} w \mid k\in\bbZ, w\in \{1, s_1^\vee\} \}.
\end{equation}

The following is the alcove picture for the extended affine Weyl group $W^\vee$, showing the correspondence between the alcoves and the elements of $W^\vee$.  The periodic orientation is indicated by $\scriptstyle{+}$ and $\scriptstyle{-}$ on either side of the hyperplanes.  The two ways of indexing the alcoves correspond to the two presentations~\eqref{eqn.presn1} and~\eqref{eqn.presn2} for $W^\vee$.

$$\beginpicture
\setcoordinatesystem units <1.8cm,1cm>         
\setplotarea x from -4 to 4.5, y from -1.5 to 2.8  
\put{$\bullet$} at 0 0 
{\small
\put{$\ttH_{-\alpha^\vee+2d}$}[t] at 2 -0.9
\put{$\ttH_{\alpha^\vee}$}[t] at 0 -0.9 
\put{$\ttH_{\alpha^\vee+2d}$}[t] at -2 -0.9 
\put{$\ttH_{-\alpha^\vee+d}$}[t] at 1 -0.9
\put{$\ttH_{-\alpha^\vee+3d}$}[t] at 3 -0.9 
\put{$\ttH_{\alpha^\vee+d}$}[t] at -1 -0.9
}
\put{Sheet $1$} at -4.2 0
\plot -3.5 0  4.5 0 / \plot -3 0.7 -3 -0.4 / \plot  -2 0.7  -2 -0.4 /
\plot  -1 0.7  -1 -0.4 / \plot  0 0.7  0 -0.4 / \plot  1 0.7  1 -0.4
/ \plot  2 0.7  2 -0.4 / \plot  3 0.7  3 -0.4 / \plot  4 0.7  4 -0.4
/
\put{$\scriptstyle{+}$} at -1.9 0.7 \put{$\scriptstyle{-}$} at -2.1 0.7
\put{$\scriptstyle{+}$} at -0.9 0.7 \put{$\scriptstyle{-}$} at -1.1 0.7
\put{$\scriptstyle{+}$} at 0.1 0.7 \put{$\scriptstyle{-}$} at -0.1 0.7
\put{$\scriptstyle{+}$} at 1.1 0.7 \put{$\scriptstyle{-}$} at 0.9 0.7
\put{$\scriptstyle{+}$} at 2.1 0.7 \put{$\scriptstyle{-}$} at 1.9 0.7
\put{$\scriptstyle{+}$} at 3.1 0.7 \put{$\scriptstyle{-}$} at 2.9 0.7
\put{$1$}[b] at 0.5 0.15 \put{$s_0^\vee$}[b] at 1.5 0.15 
\put{$s_0^\vee s_1^\vee$}[b] at 2.5 0.15 
\put{$s_0^\vee s_1^\vee s_0^\vee$}[b] at 3.5 0.15 
\put{$s_1^\vee$}[b] at -0.5 0.15 
\put{$s_1^\vee s_0^\vee$}[b] at -1.5 0.15 
\put{$s_1^\vee s_0^\vee s_1^\vee$}[b] at -2.5 0.15
\put{$1$}[t] at 0.5 -0.15
\put{$x^{\alpha}s_1^\vee$}[t] at 1.5 -0.15
\put{$x^{\alpha}$}[t] at 2.5 -0.15
\put{$x^{2\alpha}s_1^\vee$}[t] at 3.5 -0.15 
\put{$s_1^\vee$}[t] at -0.5 -0.15
\put{$x^{-\alpha}$}[t] at -1.5 -0.15
\put{$x^{-\alpha}s_1^\vee$}[t] at -2.5 -0.15
\put{Sheet $\pi^\vee$} at -4.2 2
\plot -3.5 2 4.5 2 / 
\plot  -3 2.7  -3 1.5 / \plot  -2 2.7  -2 1.5 /
\plot  -1 2.7  -1 1.5 / \plot  0 2.7  0 1.5 / \plot  1 2.7  1 1.5 /
\plot  2 2.7  2 1.5 / \plot  3 2.7  3 1.5 / \plot  4 2.7  4 1.5 /
\put{$\scriptstyle{+}$} at -1.9 2.7 \put{$\scriptstyle{-}$} at -2.1 2.7
\put{$\scriptstyle{+}$} at -0.9 2.7 \put{$\scriptstyle{-}$} at -1.1 2.7
\put{$\scriptstyle{+}$} at 0.1 2.7 \put{$\scriptstyle{-}$} at -0.1 2.7
\put{$\scriptstyle{+}$} at 1.1 2.7 \put{$\scriptstyle{-}$} at 0.9 2.7
\put{$\scriptstyle{+}$} at 2.1 2.7 \put{$\scriptstyle{-}$} at 1.9 2.7
\put{$\scriptstyle{+}$} at 3.1 2.7 \put{$\scriptstyle{-}$} at 2.9 2.7
\put{$\pi^\vee$}[b] at 0.5 2.25 
\put{$\pi^\vee s_1^\vee$}[b] at 1.5 2.15
\put{$\pi^\vee s_1^\vee s_0^\vee$}[b] at 2.5 2.15
\put{$\pi^\vee s_1^\vee s_0^\vee s_1^\vee$}[b] at 3.5 2.15
\put{$\pi^\vee s_0^\vee$}[b] at -0.5 2.15
\put{$\pi^\vee s_0^\vee s_1^\vee$}[b] at -1.5 2.15
\put{$\pi^\vee s_0^\vee s_1^\vee s_0^\vee$}[b] at -2.5 2.15
\put{$x^{\omega}s_1^\vee$}[t] at 0.5 1.85 
\put{$x^{\omega}$}[t] at 1.5 1.85 
\put{$x^{3\omega}s_1^\vee$}[t] at 2.5 1.85
\put{$x^{3\omega}$}[t] at 3.5 1.85 
\put{$x^{-\omega}$}[t] at -0.5 1.85 
\put{$x^{-\omega}s_1^\vee$}[t] at -1.5 1.85
\put{$x^{-3\omega}$}[t] at -2.5 1.85
\endpicture
$$

The double affine Hecke algebra over the field $\bbK = \bbQ(q^{1/2},t^{1/2})$ is the algebra generated by $\pi, T_0, T_1$, and the group $X=\{q^kX^{j\omega} \mid k\in\hbox{$\frac12$}\bbZ,  j\in \bbZ\}$ subject to the relations
$$
\pi^2=1,\qquad
T_i^2 = (t^{1/2}-t^{-1/2}) T_i+1, \textrm{ for } i=0,1, \qquad
X^{j\omega}X^{k\omega} = X^{(j+k)\omega} \textrm{ for } j,k\in \bbZ,
$$
$$\pi T_0 \pi^{-1} = T_1, \qquad
T_1 X^\omega T_1 = X^{-\omega},\qquad
\pi X^\omega \pi^{-1} = q^{1/2} X^{-\omega}.
$$

Let $(T_0^\vee)^{-1} = X^\alpha T_1$, $Y^{-\alpha^\vee} = T_1^{-1}T_0^{-1}$, $Y^{-\alpha_0^\vee} = qY^{\alpha^\vee} = qT_0T_1$.
The intertwiners are $\pi^\vee = X^\omega T_1$, and
$$\tau_i^\vee 
= T_i^\vee + \frac{t^{-1/2}-t^{1/2}}{1-Y^{-\alpha_i^\vee}}
= (T_i^\vee)^{-1} + \frac{(t^{-1/2}-t^{1/2})Y^{-\alpha_i^\vee}}{1-Y^{-\alpha_i^\vee}}, \quad\textrm{ for } i=0,1.$$

For $\mu \in \fh_\bbZ^*$, the minimal coset representatives are given by
$$m_{k\omega} =\begin{cases}
x^{k\omega}s_1 = \pi^\vee(s_1\pi^\vee)^{k-1}, & k \geq 1,\\
x^{k\omega} = (s_1\pi^\vee)^k, & k \leq 0.
\end{cases}$$

\begin{example} 
{\em 
{\bf Littlewood-Richardson formulas.}
In this example, we calculate $E_{3\omega}P_{k\omega}$ and $P_{3\omega}P_{k\omega}$ for $k\geq 3$.

Let $\mu=3\omega$ and $\lambda = k\omega$.
The minimal length representative of the coset $x^{3\omega}W_0$ is $m_{3\omega} =  s_0^\vee s_1 \pi^\vee$, so Theorem~\ref{thm.EP-E} and Theorem~\ref{thm.PP-P} give
\begin{align*}
E_{3\omega}P_{k\omega}
&= \sum_{h\in \Gamma_2^{\overline{\ttC}}} 
\left( (-1)^{|\phi_{\mathrm{grey}}^-(h)|} b_h c_{1(h)} c_{2(h)} \right)
E_{\varpi(h)},\\
P_{3\omega}P_{k\omega}
&= \sum_{h\in \Gamma_2^{\overline{\ttC}}} 
\left( (-1)^{|\phi_{\mathrm{grey}}^-(h)|} b_h e_h c_{1(h)} c_{2(h)} \right)
P_{-w_0\ttwt(h)},\\
\end{align*}
where the walks are of type $\vec{m}_{3\omega}^{-1} = (\pi^\vee, 1, 0)$, and begin in $x^{k\omega}$ or $x^{k\omega}s_1$.  There are 18 walks in all.  

According to where each walk begins and ends (the alcoves $\ttb(h)$ and $\tte(h)$), we have
$$b_h = \begin{cases} 
1, & \hbox{if } \ttb(h) = x^{k\omega},\\
\displaystyle t^{\frac12}\frac{1-q^k}{1-q^kt}, & \hbox{if } \ttb(h) = x^{k\omega}s_1,
\end{cases}
\qquad\hbox{and}\qquad
e_h = \begin{cases}
1, & \hbox{if } \tte(h) = x^{j\omega}s_1,\\
\displaystyle t^{-\frac12}\frac{1-q^jt^2}{1-q^jt} & \hbox{if } \tte(h) = x^{j\omega}.
\end{cases}
$$
Also note that for these calculations, we have 
$$Y^{-b_1^\vee}\mathbf{1} = Y^{\alpha^\vee-d}\mathbf{1} = qt\mathbf{1}, \qquad\hbox{and}\qquad
Y^{-b_2^\vee}\mathbf{1} = Y^{\alpha^\vee-2d}\mathbf{1} = q^2 t\mathbf{1}.
$$

In the following figures, they are organized into four groups so that walks in a group are generated by the same walk $p(h)$ (see~\eqref{eqn.ph}).  
In this root system, $-w_0 \ttwt(h) = \varpi(h)_+$.

\newpage
\noindent{\bf Group I.} The walks in this group are generated by 
$$\beginpicture
\setcoordinatesystem units <0.75cm, 0.75cm>         
\setplotarea x from -4 to 10, y from -0.5 to 1.5   
    \plot -1.5 1 3 1 / 
    \plot -1 0.8 -1 1.4 / \plot  0 0.8  0 1.4 / 
    \plot  1 0.8  1 1.4 / \plot  2 0.8  2 1.4 / 
    \plot -1.5 0 3 0 / 
    \plot -1 -0.2 -1 0.6 / \plot  0 -0.2  0 0.6 / 
    \plot  1 -0.2  1 0.6 / \plot  2 -0.2  2 0.6 /
    \put{$\bullet$} at 0 0
    \put{$p(h)=$} at -3 0.5
    \setplotsymbol(.)
    \arrow<6pt> [.2,.67] from 1.5 0.2 to 0.5 0.2
    \arrow<6pt> [.2,.67] from 2.5 0.2 to 1.5 0.2
    	\setdots<3pt> \plot 2.5 0.3 2.5 1.3 / 
	\put{which gives $\epsilon_1=-1$, $\epsilon_2=-1$.}[l] at 4 0.5
\endpicture
$$ 

\hspace{-0.4cm}
\begin{tabular}{ccccc}
$h$ & $b_h$ & $e_h$ & $(-1)^{|\phi_{\mathrm{grey}}^-(h)|}c_{1(h)}c_{2(h)}$ & $\varpi(h)$ \\ \\
\beginpicture
\setcoordinatesystem units <0.5cm, 0.6cm>         
\setplotarea x from -2.5 to 3, y from -1 to 1   
    \plot -3 0 3 0 / \plot -3 -1 3 -1 / 
    \plot -2 -0.2 -2 0.6 /
    \plot -1 -0.2 -1 0.6 / \plot  0 -0.2  0 0.6 / 
    \plot  1 -0.2  1 0.6 / \plot  2 -0.2  2 0.6 / 
    \plot -2 -1.2 -2 -0.8 /
    \plot -1 -1.2 -1 -0.8 / \plot  0 -1.2  0 -0.8 / 
    \plot  1 -1.2  1 -0.8 / \plot  2 -1.2  2 -0.8 /
    \put{\tiny$k\omega$} at 0 -1.4
    \put{\tiny$(k+3)\omega$}[r] at 3 1
    \setplotsymbol(.)
    \arrow <4pt> [.2,.67] from 0.5 0.3 to 1.5 0.3
    \arrow <4pt> [.2,.67] from 1.5 0.3 to 2.5 0.3 
    \setdots<3pt> \plot 0.5 -0.8 0.5 0.3 / 
\endpicture
&$1$ &$1$ &$1$ &\scriptsize$(k+3)\omega$ \\ \\
\beginpicture
\setcoordinatesystem units <0.5cm, 0.6cm>         
\setplotarea x from -2.5 to 3, y from -1 to 1   
    \plot -3 0 3 0 / \plot -3 -1 3 -1 / 
    \plot -2 -0.2 -2 0.6 /
    \plot -1 -0.2 -1 0.6 / \plot  0 -0.2  0 0.6 / 
    \plot  1 -0.2  1 0.6 / \plot  2 -0.2  2 0.6 / 
    \plot -2 -1.2 -2 -0.8 /
    \plot -1 -1.2 -1 -0.8 / \plot  0 -1.2  0 -0.8 / 
    \plot  1 -1.2  1 -0.8 / \plot  2 -1.2  2 -0.8 /
    \put{\tiny$k\omega$} at 0 -1.4
    \put{\tiny$-(k+1)\omega$} at 1.5 1
    \setplotsymbol(.)
    \arrow <4pt> [.2,.67] from 0.5 0.3 to 1.5 0.3
    \textcolor{Gray}
    { \plot 1.5 0.3 2 0.3 /
    \plot 1.97 0.3 1.97 0.5 /
    \arrow <4pt> [.2,.67] from 2 0.5 to 1.5 0.5 
    }\setdots<3pt> \plot 0.5 -0.8 0.5 0.3 / 
\endpicture
&$1$ &\footnotesize$\displaystyle t^{-\frac12}\frac{1-q^{k+1}t^2}{1-q^{k+1}t}$ 
&\footnotesize$\displaystyle -\frac{(t^{-\frac12}-t^{\frac12})q^{k+2}t}{1-q^{k+2}t}$
&\scriptsize$-(k+1)\omega$ \\ \\
\beginpicture
\setcoordinatesystem units <0.5cm, 0.6cm>         
\setplotarea x from -2.5 to 3, y from -1 to 1   
    \plot -3 0 3 0 / \plot -3 -1 3 -1 / 
    \plot -2 -0.2 -2 0.6 /
    \plot -1 -0.2 -1 0.6 / \plot  0 -0.2  0 0.6 / 
    \plot  1 -0.2  1 0.6 / \plot  2 -0.2  2 0.6 / 
    \plot -2 -1.2 -2 -0.8 /
    \plot -1 -1.2 -1 -0.8 / \plot  0 -1.2  0 -0.8 / 
    \plot  1 -1.2  1 -0.8 / \plot  2 -1.2  2 -0.8 /
    \put{\tiny$k\omega$} at 0 -1.4
    \put{\tiny$-(k-1)\omega$} at -0.5 1
    \setplotsymbol(.)
    \textcolor{Gray}
    { \plot 0.5 0.3 1 0.3 /
    \plot 0.97 0.3 0.97 0.5 /
    \arrow <4pt> [.2,.67] from 1 0.5 to 0.5 0.5 
    }\arrow <4pt> [.2,.67] from 0.5 0.5 to -0.5 0.5
    \setdots<3pt> \plot 0.5 -0.8 0.5 0.3 / 
\endpicture
&$1$ &\footnotesize$\displaystyle t^{-\frac12}\frac{1-q^{k-1}t^2}{1-q^{k-1}t}$ 
&\footnotesize$\displaystyle
-\left(\frac{(t^{-\frac12}-t^{\frac12})q^{k+1}t}{1-q^{k+1}t}\right)
\left(\frac{1-q^k}{1-q^kt}\frac{1-q^kt^2}{1-q^kt}\right)$
&\scriptsize$-(k-3)\omega$ \\ \\
\beginpicture
\setcoordinatesystem units <0.5cm, 0.6cm>         
\setplotarea x from -2.5 to 3, y from -1 to 1   
    \plot -3 0 3 0 / \plot -3 -1 3 -1 / 
    \plot -2 -0.2 -2 0.6 /
    \plot -1 -0.2 -1 0.6 / \plot  0 -0.2  0 0.6 / 
    \plot  1 -0.2  1 0.6 / \plot  2 -0.2  2 0.6 / 
    \plot -2 -1.2 -2 -0.8 /
    \plot -1 -1.2 -1 -0.8 / \plot  0 -1.2  0 -0.8 / 
    \plot  1 -1.2  1 -0.8 / \plot  2 -1.2  2 -0.8 /
    \put{\tiny$k\omega$} at 0 -1.4
    \put{\tiny$(k+1)\omega$} at 0.5 1
    \setplotsymbol(.)
    \textcolor{Gray}
    {\plot 0.5 0.2 1 0.2 /
    \plot 0.97 0.2 0.97 0.35 /
    \arrow <4pt> [.2,.67] from 1 0.35 to 0.5 0.35 
    \plot 0.5 0.35 0 0.35 /
    \plot 0.03 0.35 0.03 0.5 /
    \arrow <4pt> [.2,.67] from 0 0.5 to 0.5 0.5 
    }\setdots<3pt> \plot 0.5 -0.8 0.5 0.2 /
    \endpicture
&$1$ &$1$ 
&\footnotesize$\displaystyle -
\left(\frac{(t^{-\frac12}-t^{\frac12})q^{k+1}t}{1-q^{k+1}t}\right)
\left(\frac{t^{-\frac12}-t^{\frac12}}{1-q^kt}\right)$
&\scriptsize$(k+1)\omega$ \\ \\
\beginpicture
\setcoordinatesystem units <0.5cm, 0.6cm>         
\setplotarea x from -2.5 to 3, y from -1 to 1   
    \plot -3 0 3 0 / \plot -3 -1 3 -1 / 
    \plot -2 -0.2 -2 0.6 /
    \plot -1 -0.2 -1 0.6 / \plot  0 -0.2  0 0.6 / 
    \plot  1 -0.2  1 0.6 / \plot  2 -0.2  2 0.6 / 
    \plot -2 -1.2 -2 -0.8 /
    \plot -1 -1.2 -1 -0.8 / \plot  0 -1.2  0 -0.8 / 
    \plot  1 -1.2  1 -0.8 / \plot  2 -1.2  2 -0.8 /
    \put{\tiny$k\omega$} at 0 -1.4
    \put{\tiny$-(k-3)\omega$}[l] at -2.5 1
    \setplotsymbol(.)
    \arrow <4pt> [.2,.67] from -0.5 0.3 to -1.5 0.3
    \arrow <4pt> [.2,.67] from -1.5 0.3 to -2.5 0.3 
    \setdots<3pt> \plot -0.5 -0.8 -0.5 0.3 /
\endpicture
&\footnotesize$\displaystyle t^{\frac12}\frac{1-q^k}{1-q^kt}$ 
&\footnotesize$\displaystyle t^{-\frac12}\frac{1-q^{k-3}t^2}{1-q^{k-3}t}$ 
&\footnotesize$\displaystyle 
\left(\!\frac{1-q^{k-1}}{1-q^{k-1}t}\frac{1-q^{k-1}t^2}{1-q^{k-1}t}\!\right)
\!\!\!
\left(\!\frac{1-q^{k-2}}{1-q^{k-2}t}\frac{1-q^{k-2}t^2}{1-q^{k-2}t}\!\right)$ 
&\scriptsize$-(k-3)\omega$ \\ \\
\beginpicture
\setcoordinatesystem units <0.5cm, 0.6cm>         
\setplotarea x from -2.5 to 3, y from -1 to 1   
    \plot -3 0 3 0 / \plot -3 -1 3 -1 / 
    \plot -2 -0.2 -2 0.6 /
    \plot -1 -0.2 -1 0.6 / \plot  0 -0.2  0 0.6 / 
    \plot  1 -0.2  1 0.6 / \plot  2 -0.2  2 0.6 / 
    \plot -2 -1.2 -2 -0.8 /
    \plot -1 -1.2 -1 -0.8 / \plot  0 -1.2  0 -0.8 / 
    \plot  1 -1.2  1 -0.8 / \plot  2 -1.2  2 -0.8 /
    \put{\tiny$k\omega$} at 0 -1.4
    \put{\tiny$(k-1)\omega$} at -1.5 1
    \setplotsymbol(.)
    \arrow <4pt> [.2,.67] from -0.5 0.3 to -1.5 0.3
    \textcolor{Gray}
    { \plot -1.5 0.3 -2 0.3 /
    \plot -1.97 0.3 -1.97 0.5 /
    \arrow <4pt> [.2,.67] from -2 0.5 to -1.5 0.5 
    }\setdots<3pt> \plot -0.5 -0.8 -0.5 0.3 /
    \endpicture
&\footnotesize$\displaystyle t^{\frac12}\frac{1-q^k}{1-q^kt}$ &$1$ 
&\footnotesize$\displaystyle 
\left(\frac{1-q^{k-1}}{1-q^{k-1}t}\frac{1-q^{k-1}t^2}{1-q^{k-1}t}\right) 
\left(\frac{t^{-\frac12}-t^{\frac12}}{1-q^{k-2}t}\right)$
&\scriptsize$(k-1)\omega$ \\ \\
\beginpicture
\setcoordinatesystem units <0.5cm, 0.6cm>         
\setplotarea x from -2.5 to 3, y from -1 to 1   
    \plot -3 0 3 0 / \plot -3 -1 3 -1 / 
    \plot -2 -0.2 -2 0.6 /
    \plot -1 -0.2 -1 0.6 / \plot  0 -0.2  0 0.6 / 
    \plot  1 -0.2  1 0.6 / \plot  2 -0.2  2 0.6 / 
    \plot -2 -1.2 -2 -0.8 /
    \plot -1 -1.2 -1 -0.8 / \plot  0 -1.2  0 -0.8 / 
    \plot  1 -1.2  1 -0.8 / \plot  2 -1.2  2 -0.8 /
    \put{\tiny$k\omega$} at 0 -1.4
    \put{\tiny$(k+1)\omega$} at 0.5 1
    \setplotsymbol(.)
    \textcolor{Gray}
    { \plot -0.5 0.3 -1 0.3 /
    \plot -0.97 0.3 -0.97 0.5 /
    \arrow <4pt> [.2,.67] from -1 0.5 to -0.5 0.5 
    }\arrow <4pt> [.2,.67] from -0.5 0.5 to 0.5 0.5
    \setdots<3pt> \plot -0.5 -0.8 -0.5 0.3 /
\endpicture
&\footnotesize$\displaystyle t^{\frac12}\frac{1-q^k}{1-q^kt}$ &$1$ 
&\footnotesize$\displaystyle \frac{t^{-\frac12}-t^{\frac12}}{1-q^{k-1}t}$
&\scriptsize$(k+1)\omega$ \\ \\
\beginpicture
\setcoordinatesystem units <0.5cm, 0.6cm>         
\setplotarea x from -2.5 to 3, y from -1 to 1   
    \plot -3 0 3 0 / \plot -3 -1 3 -1 / 
    \plot -2 -0.2 -2 0.6 /
    \plot -1 -0.2 -1 0.6 / \plot  0 -0.2  0 0.6 / 
    \plot  1 -0.2  1 0.6 / \plot  2 -0.2  2 0.6 / 
    \plot -2 -1.2 -2 -0.8 /
    \plot -1 -1.2 -1 -0.8 / \plot  0 -1.2  0 -0.8 / 
    \plot  1 -1.2  1 -0.8 / \plot  2 -1.2  2 -0.8 /
    \put{\tiny$k\omega$} at 0 -1.4
    \put{\tiny$-(k-3)\omega$} at -0.5 1
    \setplotsymbol(.)
    \textcolor{Gray}
    {\plot -0.5 0.2 -1 0.2 /
    \plot -0.97 0.2 -0.97 0.35 /
    \arrow <4pt> [.2,.67] from -1 0.35 to -0.5 0.35 
    \plot -0.5 0.35 0 0.35 /
    \plot -0.03 0.35 -0.03 0.5 /
    \arrow <4pt> [.2,.67] from 0 0.5 to -0.5 0.5 
    }\setdots<3pt> \plot -0.5 -0.8 -0.5 0.2 /
    \endpicture
&\footnotesize$\displaystyle t^{\frac12}\frac{1-q^k}{1-q^kt}$ 
&\footnotesize$\displaystyle t^{-\frac12}\frac{1-q^{k-1}t^2}{1-q^{k-1}t}$ 
&\footnotesize$\displaystyle 
-\left(\frac{t^{-\frac12}-t^{\frac12}}{1-q^{k-1}t}\right)
\left(\frac{(t^{-\frac12}-t^{\frac12})q^kt}{1-q^kt}\right)$
&\scriptsize$-(k-3)\omega$\\ \\
\end{tabular}

\noindent{\bf Group II.} The walks in this group are generated by 
$$\hspace{-0.4cm}
\beginpicture
\setcoordinatesystem units <0.75cm, 0.75cm>         
\setplotarea x from -4 to 10, y from -0.35 to 2   
    \plot -1.5 1 3 1 / 
    \plot -1 0.8 -1 1.4 / \plot  0 0.8  0 1.4 / 
    \plot  1 0.8  1 1.4 / \plot  2 0.8  2 1.4 / 
    \plot -1.5 0 3 0 / 
    \plot -1 -0.2 -1 0.6 / \plot  0 -0.2  0 0.6 / 
    \plot  1 -0.2  1 0.6 / \plot  2 -0.2  2 0.6 /
    \put{$\bullet$} at 0 0
    \put{$p(h)=$} at -3 0.5
    \setplotsymbol(.)
    \arrow <6pt> [.2,.67] from 1.5 0.2 to 0.5 0.2
    \plot 1.5 0.4 2 0.4 /
    \plot 1.97 0.2 1.97 0.4 /
    \arrow <6pt> [.2,.67] from 2 0.2 to 1.5 0.2 
    \setdots<3pt> \plot 1.5 0.4 1.5 1.3 / 
    \put{which gives $\epsilon_2=-1$.}[l] at 4 0.5
\endpicture
$$

\begin{center}
\begin{tabular}{ccccc} 
$h$ &  $b_h$ & $e_h$ & $(-1)^{\phi_{\mathrm{grey}}^-(h)}c_{1(h)}c_{2(h)}$
& $\varpi(h)$ \\  \\
\beginpicture
\setcoordinatesystem units <0.5cm, 0.6cm>         
\setplotarea x from -4 to 4, y from -1 to 1   
    \plot -3 0 3 0 / \plot -3 -1 3 -1 / 
    \plot -2 -0.2 -2 0.6 /
    \plot -1 -0.2 -1 0.6 / \plot  0 -0.2  0 0.6 / 
    \plot  1 -0.2  1 0.6 / \plot  2 -0.2  2 0.6 / 
    \plot -2 -1.2 -2 -0.8 /
    \plot -1 -1.2 -1 -0.8 / \plot  0 -1.2  0 -0.8 / 
    \plot  1 -1.2  1 -0.8 / \plot  2 -1.2  2 -0.8 /
    \put{\tiny$k\omega$} at 0 -1.4
    \put{\tiny$-(k-1)\omega$} at -0.5 1
    \setplotsymbol(.)
    \textcolor{Black}
    { \plot 0.5 0.3 1 0.3 /
    \plot 0.97 0.3 0.97 0.5 /
    \arrow <4pt> [.2,.67] from 1 0.5 to 0.5 0.5 
    }\arrow <4pt> [.2,.67] from 0.5 0.5 to -0.5 0.5
    \setdots<3pt> \plot 0.5 -0.8 0.5 0.3 /
\endpicture
&$1$ 
&\footnotesize$\displaystyle t^{-\frac12}\frac{1-q^{k-1}t^2}{1-q^{k-1}t}$ 
&\footnotesize$\displaystyle
\left(\frac{(t^{-\frac12}-t^{\frac12})qt}{1-qt}\right)
\left(\frac{1-q^k}{1-q^kt}\frac{1-q^kt^2}{1-q^kt}\right)$
&\footnotesize$-(k-1)\omega$ \\ \\
\beginpicture
\setcoordinatesystem units <0.5cm, 0.6cm>         
\setplotarea x from -4 to 4, y from -1 to 1   
    \plot -3 0 3 0 / \plot -3 -1 3 -1 / 
    \plot -2 -0.2 -2 0.6 /
    \plot -1 -0.2 -1 0.6 / \plot  0 -0.2  0 0.6 / 
    \plot  1 -0.2  1 0.6 / \plot  2 -0.2  2 0.6 / 
    \plot -2 -1.2 -2 -0.8 /
    \plot -1 -1.2 -1 -0.8 / \plot  0 -1.2  0 -0.8 / 
    \plot  1 -1.2  1 -0.8 / \plot  2 -1.2  2 -0.8 /
    \put{\tiny$k\omega$} at 0 -1.4
    \put{\tiny$(k+1)\omega$} at 0.5 1
    \setplotsymbol(.)
    \textcolor{Black}
    {\plot 0.5 0.2 1 0.2 /
    \plot 0.97 0.2 0.97 0.35 /
    \arrow <4pt> [.2,.67] from 1 0.35 to 0.5 0.35 
    }\textcolor{Gray} 
    {\plot 0.5 0.35 0 0.35 /
    \plot 0.03 0.35 0.03 0.5 /
    \arrow <4pt> [.2,.67] from 0 0.5 to 0.5 0.5 
    }\setdots<3pt> \plot 0.5 -0.8 0.5 0.2 /
    \endpicture
&$1$ &$1$ 
&\footnotesize$\displaystyle
\left(\frac{(t^{-\frac12}-t^{\frac12})qt}{1-qt}\right)
\left(\frac{t^{-\frac12}-t^{\frac12}}{1-q^kt}\right)$
&\footnotesize$(k+1)\omega$ \\ \\
\beginpicture
\setcoordinatesystem units <0.5cm, 0.6cm>         
\setplotarea x from -4 to 4, y from -1 to 1   
    \plot -3 0 3 0 / \plot -3 -1 3 -1 / 
    \plot -2 -0.2 -2 0.6 /
    \plot -1 -0.2 -1 0.6 / \plot  0 -0.2  0 0.6 / 
    \plot  1 -0.2  1 0.6 / \plot  2 -0.2  2 0.6 / 
    \plot -2 -1.2 -2 -0.8 /
    \plot -1 -1.2 -1 -0.8 / \plot  0 -1.2  0 -0.8 / 
    \plot  1 -1.2  1 -0.8 / \plot  2 -1.2  2 -0.8 /
    \put{\tiny$k\omega$} at 0 -1.4
    \put{\tiny$(k+1)\omega$} at 0.5 1
    \setplotsymbol(.)
    \textcolor{Black}
    { \plot -0.5 0.3 -1 0.3 /
    \plot -0.97 0.3 -0.97 0.5 /
    \arrow <4pt> [.2,.67] from -1 0.5 to -0.5 0.5 
    }\arrow <4pt> [.2,.67] from -0.5 0.5 to 0.5 0.5
    \setdots<3pt> \plot -0.5 -0.8 -0.5 0.3 /
\endpicture
&\footnotesize$\displaystyle t^{\frac12}\frac{1-q^k}{1-q^kt}$ &$1$ 
&\footnotesize$\displaystyle\frac{(t^{-\frac12}-t^{\frac12})qt}{1-qt}$
&\footnotesize$(k+1)\omega$ \\ \\
\beginpicture
\setcoordinatesystem units <0.5cm, 0.6cm>         
\setplotarea x from -4 to 4, y from -1 to 1   
    \plot -3 0 3 0 / \plot -3 -1 3 -1 / 
    \plot -2 -0.2 -2 0.6 /
    \plot -1 -0.2 -1 0.6 / \plot  0 -0.2  0 0.6 / 
    \plot  1 -0.2  1 0.6 / \plot  2 -0.2  2 0.6 / 
    \plot -2 -1.2 -2 -0.8 /
    \plot -1 -1.2 -1 -0.8 / \plot  0 -1.2  0 -0.8 / 
    \plot  1 -1.2  1 -0.8 / \plot  2 -1.2  2 -0.8 /
    \put{\tiny$k\omega$} at 0 -1.4
    \put{\tiny$-(k-1)\omega$} at -0.5 1
    \setplotsymbol(.)
    \textcolor{Black}
    {\plot -0.5 0.2 -1 0.2 /
    \plot -0.97 0.2 -0.97 0.35 /
    \arrow <4pt> [.2,.67] from -1 0.35 to -0.5 0.35
    }\textcolor{Gray} 
    {\plot -0.5 0.35 0 0.35 /
    \plot -0.03 0.35 -0.03 0.5 /
    \arrow <4pt> [.2,.67] from 0 0.5 to -0.5 0.5 
    }\setdots<3pt> \plot -0.5 -0.8 -0.5 0.2 /
    \endpicture
&\footnotesize$\displaystyle t^{\frac12}\frac{1-q^k}{1-q^kt}$ 
&\footnotesize$\displaystyle t^{-\frac12}\frac{1-q^{k-1}t^2}{1-q^{k-1}t}$ 
&\footnotesize$\displaystyle 
-\left(\frac{(t^{-\frac12}-t^{\frac12})qt}{1-qt}\right)
\left(\frac{(t^{-\frac12}-t^{\frac12})q^kt}{1-q^kt}\right)$
&\footnotesize$-(k-1)\omega$\\ \\
\end{tabular}
\end{center}

\newpage
\noindent{\bf Group III.} The walks in this group are generated by 
$$\hspace{-0.4cm}
\beginpicture
\setcoordinatesystem units <0.75cm, 0.75cm>         
\setplotarea x from -4 to 10, y from -0.35 to 2   
    \plot -1.5 1 3 1 / 
    \plot -1 0.8 -1 1.4 / \plot  0 0.8  0 1.4 / 
    \plot  1 0.8  1 1.4 / \plot  2 0.8  2 1.4 / 
    \plot -1.5 0 3 0 / 
    \plot -1 -0.2 -1 0.6 / \plot  0 -0.2  0 0.6 / 
    \plot  1 -0.2  1 0.6 / \plot  2 -0.2  2 0.6 /
    \put{$\bullet$} at 0 0
    \put{$p(h)=$} at -3 0.5
    \setplotsymbol(.)
    \plot 0.5 0.4 1 0.4 /
    \plot 0.97 0.2 0.97 0.4 /
    \arrow <6pt> [.2,.67] from 1 0.2 to 0.5 0.2 
    \arrow <6pt> [.2,.67] from -0.5 0.4 to 0.5 0.4
    	\setdots<3pt> \plot -0.5 0.4 -0.5 1.3 / 
	\put{which gives $\epsilon_1=+1$.}[l] at 4 0.5
\endpicture
$$

\begin{center}
\begin{tabular}{cccccc}
$h$ & $b_h$ & $e_h$ & $(-1)^{\phi_{\mathrm{grey}}^-(h)}c_{1(h)}c_{2(h)}$
&$\varpi(h)$ 
\\ \\
\beginpicture
\setcoordinatesystem units <0.5cm, 0.6cm>         
\setplotarea x from -3 to 4, y from -1 to 1   
    \plot -3 0 3 0 / \plot -3 -1 3 -1 / 
    \plot -2 -0.2 -2 0.6 /
    \plot -1 -0.2 -1 0.6 / \plot  0 -0.2  0 0.6 / 
    \plot  1 -0.2  1 0.6 / \plot  2 -0.2  2 0.6 / 
    \plot -2 -1.2 -2 -0.8 /
    \plot -1 -1.2 -1 -0.8 / \plot  0 -1.2  0 -0.8 / 
    \plot  1 -1.2  1 -0.8 / \plot  2 -1.2  2 -0.8 /
    \put{\tiny$k\omega$} at 0 -1.4
    \put{\tiny$-(k+1)\omega$} at 1.5 1
    \setplotsymbol(.)
    \arrow <4pt> [.2,.67] from 0.5 0.3 to 1.5 0.3
    \textcolor{Black}
    { \plot 1.5 0.3 2 0.3 /
    \plot 1.97 0.3 1.97 0.5 /
    \arrow <4pt> [.2,.67] from 2 0.5 to 1.5 0.5 
    }\setdots<3pt> \plot 0.5 -0.8 0.5 0.3 /
    \endpicture
&$1$ 
&\footnotesize$\displaystyle t^{-\frac12}\frac{1-q^{k+1}t^2}{1-q^{k+1}t}$ 
&\footnotesize$\displaystyle \frac{(t^{-\frac12}-t^{\frac12})q^2t}{1-q^2t}$
&\footnotesize$-(k+1)\omega$  \\ \\
\beginpicture
\setcoordinatesystem units <0.5cm, 0.6cm>         
\setplotarea x from -3 to 4, y from -1 to 1   
    \plot -3 0 3 0 / \plot -3 -1 3 -1 / 
    \plot -2 -0.2 -2 0.6 /
    \plot -1 -0.2 -1 0.6 / \plot  0 -0.2  0 0.6 / 
    \plot  1 -0.2  1 0.6 / \plot  2 -0.2  2 0.6 / 
    \plot -2 -1.2 -2 -0.8 /
    \plot -1 -1.2 -1 -0.8 / \plot  0 -1.2  0 -0.8 / 
    \plot  1 -1.2  1 -0.8 / \plot  2 -1.2  2 -0.8 /
    \put{\tiny$k\omega$} at 0 -1.4
    \put{\tiny$(k+1)\omega$} at 0.5 1
    \setplotsymbol(.)
    \textcolor{Gray}
    {\plot 0.5 0.2 1 0.2 /
    \plot 0.97 0.2 0.97 0.35 /
    \arrow <4pt> [.2,.67] from 1 0.35 to 0.5 0.35 
    }\textcolor{Black}
    {\plot 0.5 0.35 0 0.35 /
    \plot 0.03 0.35 0.03 0.5 /
    \arrow <4pt> [.2,.67] from 0 0.5 to 0.5 0.5 
    }\setdots<3pt> \plot 0.5 -0.8 0.5 0.2 /
    \endpicture
&$1$ &$1$ 
&\footnotesize$\displaystyle 
-\left(\frac{t^{-\frac12}-t^{\frac12}}{1-q^{k+1}t}\right)
\left(\frac{(t^{-\frac12}-t^{\frac12})q^2t}{1-q^2t}\right)$
&\footnotesize$(k+1)\omega$ \\ \\
\beginpicture
\setcoordinatesystem units <0.5cm, 0.6cm>         
\setplotarea x from -3 to 4, y from -1 to 1   
    \plot -3 0 3 0 / \plot -3 -1 3 -1 / 
    \plot -2 -0.2 -2 0.6 /
    \plot -1 -0.2 -1 0.6 / \plot  0 -0.2  0 0.6 / 
    \plot  1 -0.2  1 0.6 / \plot  2 -0.2  2 0.6 / 
    \plot -2 -1.2 -2 -0.8 /
    \plot -1 -1.2 -1 -0.8 / \plot  0 -1.2  0 -0.8 / 
    \plot  1 -1.2  1 -0.8 / \plot  2 -1.2  2 -0.8 /
    \put{\tiny$k\omega$} at 0 -1.4
    \put{\tiny$(k-1)\omega$} at -1.5 1
    \setplotsymbol(.)
    \arrow <4pt> [.2,.67] from -0.5 0.3 to -1.5 0.3
    \textcolor{Black}
    { \plot -1.5 0.3 -2 0.3 /
    \plot -1.97 0.3 -1.97 0.5 /
    \arrow <4pt> [.2,.67] from -2 0.5 to -1.5 0.5 
    }\setdots<3pt> \plot -0.5 -0.8 -0.5 0.3 /
    \endpicture
&\footnotesize$\displaystyle t^{\frac12}\frac{1-q^k}{1-q^kt}$ 
&$1$ 
&\footnotesize$\displaystyle 
\left(\frac{1-q^{k-1}}{1-q^{k-1}t}\frac{1-q^{k-1}t^2}{1-q^{k-1}t}\right)
\left(\frac{(t^{-\frac12}-t^{\frac12})q^2t}{1-q^2t}\right)$ 
&\footnotesize$(k-1)\omega$
\\ \\
\beginpicture
\setcoordinatesystem units <0.5cm, 0.6cm>         
\setplotarea x from -3 to 4, y from -1 to 1   
    \plot -3 0 3 0 / \plot -3 -1 3 -1 / 
    \plot -2 -0.2 -2 0.6 /
    \plot -1 -0.2 -1 0.6 / \plot  0 -0.2  0 0.6 / 
    \plot  1 -0.2  1 0.6 / \plot  2 -0.2  2 0.6 / 
    \plot -2 -1.2 -2 -0.8 /
    \plot -1 -1.2 -1 -0.8 / \plot  0 -1.2  0 -0.8 / 
    \plot  1 -1.2  1 -0.8 / \plot  2 -1.2  2 -0.8 /
    \put{\tiny$k\omega$} at 0 -1.4
    \put{\tiny$-(k-1)\omega$} at -0.5 1
    \setplotsymbol(.)
    \textcolor{Gray}
    {\plot -0.5 0.2 -1 0.2 /
    \plot -0.97 0.2 -0.97 0.35 /
    \arrow <4pt> [.2,.67] from -1 0.35 to -0.5 0.35 
    }\textcolor{Black}
    {\plot -0.5 0.35 0 0.35 /
    \plot -0.03 0.35 -0.03 0.5 /
    \arrow <4pt> [.2,.67] from 0 0.5 to -0.5 0.5 
    }\setdots<3pt> \plot -0.5 -0.8 -0.5 0.2 /
    \endpicture
&\footnotesize$\displaystyle t^{\frac12}\frac{1-q^k}{1-q^kt}$ 
&\footnotesize$\displaystyle t^{-\frac12}\frac{1-q^{k-1}t^2}{1-q^{k-1}t}$ 
&\footnotesize$\displaystyle
\left(\frac{(t^{-\frac12}-t^{\frac12})q^{k-1}t}{1-q^{k-1}t}\right)
\left(\frac{(t^{-\frac12}-t^{\frac12})q^2t}{1-q^2t}\right)$
&\footnotesize$-(k-1)\omega$\\ \\
\end{tabular}
\end{center}

\noindent{\bf Group IV.} The walks in this group are generated by 
$$\hspace{-0.4cm}
\beginpicture
\setcoordinatesystem units <0.75cm, 0.75cm>         
\setplotarea x from -4 to 4, y from -0.35 to 2   
    \plot -1.5 1 3 1 / 
    \plot -1 0.8 -1 1.4 / \plot  0 0.8  0 1.4 / 
    \plot  1 0.8  1 1.4 / \plot  2 0.8  2 1.4 / 
    \plot -1.5 0 3 0 / 
    \plot -1 -0.2 -1 0.6 / \plot  0 -0.2  0 0.6 / 
    \plot  1 -0.2  1 0.6 / \plot  2 -0.2  2 0.6 /
    \put{$\bullet$} at 0 0
    \put{$p(h)=$} at -3 0.5
    \setplotsymbol(.)
	\plot 0.5 0.35 1 0.35 /
    \plot 0.97 0.15 0.97 0.35 /
    \arrow <6pt> [.2,.67] from 1 0.15 to 0.5 0.15 
    \plot 0.5 0.55 0 0.55 /
    \plot 0.03 0.35 0.03 0.55 /
    \arrow <6pt> [.2,.67] from 0 0.35 to 0.5 0.35 
    \setdots<3pt> \plot 0.5 0.55 0.5 1.3 /
\endpicture
$$

\begin{center}
\begin{tabular}{cccccc}
$h$ & $b_h$ & $e_h$ & $(-1)^{\phi_{\mathrm{grey}}^-(h)}c_{1(h)}c_{2(h)}$ 
&$\varpi(h)$ \\
\beginpicture
\setcoordinatesystem units <0.5cm, 0.6cm>         
\setplotarea x from -4 to 4, y from -1 to 1   
    \plot -3 0 3 0 / \plot -3 -1 3 -1 / 
    \plot -2 -0.2 -2 0.6 /
    \plot -1 -0.2 -1 0.6 / \plot  0 -0.2  0 0.6 / 
    \plot  1 -0.2  1 0.6 / \plot  2 -0.2  2 0.6 / 
    \plot -2 -1.2 -2 -0.8 /
    \plot -1 -1.2 -1 -0.8 / \plot  0 -1.2  0 -0.8 / 
    \plot  1 -1.2  1 -0.8 / \plot  2 -1.2  2 -0.8 /
    \put{\tiny$k\omega$} at 0 -1.4
    \put{\tiny$(k+1)\omega$} at 0.5 1
    \setplotsymbol(.)
    \textcolor{Black}
    {\plot 0.5 0.2 1 0.2 /
    \plot 0.97 0.2 0.97 0.35 /
    \arrow <4pt> [.2,.67] from 1 0.35 to 0.5 0.35 
    \plot 0.5 0.35 0 0.35 /
    \plot 0.03 0.35 0.03 0.5 /
    \arrow <4pt> [.2,.67] from 0 0.5 to 0.5 0.5 
    }\setdots<3pt> \plot 0.5 -0.8 0.5 0.2 /
    \endpicture
&$1$ &$1$ 
&\footnotesize$\displaystyle 
\left(\frac{(t^{-\frac12}-t^{\frac12})q^2t}{1-q^2t}\right)
\left(\frac{(t^{-\frac12}-t^{\frac12})}{1-qt}\right)$
&\footnotesize$(k+1)\omega$ \\ \\
\beginpicture
\setcoordinatesystem units <0.5cm, 0.6cm>         
\setplotarea x from -4 to 4, y from -1 to 1   
    \plot -3 0 3 0 / \plot -3 -1 3 -1 / 
    \plot -2 -0.2 -2 0.6 /
    \plot -1 -0.2 -1 0.6 / \plot  0 -0.2  0 0.6 / 
    \plot  1 -0.2  1 0.6 / \plot  2 -0.2  2 0.6 / 
    \plot -2 -1.2 -2 -0.8 /
    \plot -1 -1.2 -1 -0.8 / \plot  0 -1.2  0 -0.8 / 
    \plot  1 -1.2  1 -0.8 / \plot  2 -1.2  2 -0.8 /
    \put{\tiny$k\omega$} at 0 -1.4
    \put{\tiny$-(k-1)\omega$} at -0.5 1
    \setplotsymbol(.)
    \textcolor{Black}
    {\plot -0.5 0.2 -1 0.2 /
    \plot -0.97 0.2 -0.97 0.35 /
    \arrow <4pt> [.2,.67] from -1 0.35 to -0.5 0.35 
    \plot -0.5 0.35 0 0.35 /
    \plot -0.03 0.35 -0.03 0.5 /
    \arrow <4pt> [.2,.67] from 0 0.5 to -0.5 0.5 
    }\setdots<3pt> \plot -0.5 -0.8 -0.5 0.2 /
    \endpicture
&\quad\footnotesize$\displaystyle t^{\frac12}\frac{1-q^k}{1-q^kt}$ 
&\quad\footnotesize$\displaystyle t^{-\frac12}\frac{1-q^{k-1}t^2}{1-q^{k-1}t}$ 
&\footnotesize$\displaystyle 
\left(\frac{(t^{-\frac12}-t^{\frac12})q^2t}{1-q^2t}\right)
\left(\frac{(t^{-\frac12}-t^{\frac12})}{1-qt}\right)$
&\footnotesize$(k+1)\omega$ \\ \\
\end{tabular}
\end{center}

Thus $E_{3\omega}P_{k\omega}$ is a linear combination of $E_{(k+3)\omega}$, $E_{-(k+1)\omega}$, $E_{(k+1)\omega}$, $E_{-(k-1)\omega}$, $E_{(k-1)\omega}$, $E_{-(k-3)\omega}$, and after simplification,
\begin{align*}
E_{3\omega}P_{k\omega}
= E_{(k+3)\omega}
&+ q^2\frac{1-t}{1-q^2t} \frac{1-q^k}{1-q^{k+2}t} t^{1/2} 
	E_{-(k+1)\omega}\\
&+ \frac{1-t}{1-q}\frac{1-q^2}{1-q^2t}\frac{1-q^k}{1-q^{k-1}t}
	\frac{1-q^{k+1}t^2}{1-q^{k+1}t} E_{(k+1)\omega}\\
&+ q\frac{1-t}{1-q} \frac{1-q^2}{1-q^2t} \frac{1-q^{k-1}}{1-q^{k-1}t}
	\frac{1-q^k}{1-q^kt} \frac{1-q^kt^2}{1-q^{k+1}t} 
	t^{1/2}E_{-(k-1)\omega}\\
&+ \frac{1-t}{1-q^2t}\frac{1-q^{k-1}}{1-q^{k-2}t}
	\frac{1-q^k}{1-q^{k-1}t}\frac{1-q^{k-1}t^2}{1-q^{k-1}t}
	\frac{1-q^kt^2}{1-q^kt} E_{(k-1)\omega}\\
&+ \frac{1-q^{k-2}}{1-q^{k-2}t}\frac{1-q^{k-2}t^2}{1-q^{k-2}t}
	\frac{1-q^{k-1}}{1-q^{k-1}t}\frac{1-q^{k-1}t^2}{1-q^{k-1}t}
	\frac{1-q^k}{1-q^kt} t^{1/2}E_{-(k-3)\omega}.
\end{align*}
Proposition~\ref{prop.domP} gives
$$P_{-j\omega} =  t^{-\frac12}\frac{1-q^jt^2}{1-q^jt} P_{j\omega} = e_{h} P_{j\omega}$$ 
if $\tte(h) = x^{j\omega}$ is a dominant weight.  Thus, after simplification,
\begin{align*}
P_{3\omega}P_{k\omega}
= P_{(k+3)\omega} 
&+ \frac{1-t}{1-q} \frac{1-q^3}{1-q^2t} 
	\frac{1-q^k}{1-q^{k-1}t} 
	\frac{1-q^{k+1}t^2}{1-q^{k+2}t} P_{(k+1)\omega} \\
&+ \frac{1-t}{1-q} \frac{1-q^3}{1-q^2t} 
	\frac{1-q^{k-1}}{1-q^{k-2}t} \frac{1-q^k}{1-q^{k-1}t}
	\frac{1-q^{k-1}t^2}{1-q^kt} \frac{1-q^kt^2}{1-q^{k+1}t}
	P_{(k-1)\omega}\\
&+ \frac{1-q^{k-2}}{1-q^{k-2}t}\frac{1-q^{k-2}t^2}{1-q^{k-2}t}
	\frac{1-q^{k-1}}{1-q^{k-1}t}\frac{1-q^{k-1}t^2}{1-q^{k-1}t}
	\frac{1-q^k}{1-q^kt} \frac{1-q^{k-3}t^2}{1-q^{k-3}t}
	P_{(k-3)\omega}.
\end{align*}

Consider the case $q=0$ where the symmetric Macdonald polynomials become Hall-Littlewood polynomials.  Following the discussion in Section~\ref{sec.specialization}, the walks giving a nonzero contribution to the sum are those whose only folds are positive and grey.  The expression
$$P_{3\omega}(t)P_{k\omega}(t) = P_{(k+3)\omega}(t)
+(1-t)P_{(k+1)\omega}(t)
+(1-t)P_{(k-1)\omega}(t)
+P_{(k-3)\omega}(t),$$
is given by four positively folded walks, and coincides with the Littlewood-Richardson formulas~\cite[Theorem 1.3]{Sc06} and~\cite[Theorem 4.9]{R06} for Hall-Littlewood polynomials.
We also mention that
$$E_{3\omega}(t)P_{k\omega}(t) = E_{(k+3)\omega}(t)  
+(1-t)E_{(k+1)\omega}(t)
+(1-t)E_{(k-1)\omega}(t)
+t^{1/2}E_{-(k-3)\omega}(t).$$

In the case $q=t=0$, the symmetric Macdonald polynomials are Schur polynomials $P_\mu(0,0)= s_\mu$, and the nonsymmetric Macdonald polynomials are Demazure characters $E_\mu(0,0) = \mathcal{A}_\mu$.  The dominant weights $k\omega$ correspond to the partitions $(k,0)$, and the weights $-k\omega$ correspond to the compositions $(0,k)$ for $k\geq 0$.  The four positively folded walks correspond to the four Littlewood-Richardson tableaux that give the classical Littlewood-Richardson formula for Schur functions:
$$s_{3\omega}s_{k\omega}= s_{(k+3)\omega} 
+s_{(k+1)\omega}
+s_{(k-1)\omega}
+s_{(k-3)\omega}.$$
By normalizing the nonsymmetric polynomials $E_\mu$ so that the coefficient of the monomial $X^\mu$ in $E_\mu$ is $1$, (see the paragraph following~\eqref{eqn.E}), then each term in
$$\calA_{3\omega}s_{k\omega} = \calA_{(k+3)\omega}
+\calA_{(k+1)\omega}
+\calA_{(k-1)\omega}
+\calA_{-(k-3)\omega},$$
corresponds to a skyline filling in the formula~\cite[Theorem 6.1]{HLMvW09}.
\hfill$\diamond$

}\end{example}

\begin{example}{\em
{\bf Pieri formulas.}  
The weight $\omega$ is minuscule in the $\fsl_2$ root system.
The two walks which give the expressions 
\begin{align*}
E_{\omega}P_{k\omega} &= E_{(k+1)\omega} + \frac{1-q^k}{1-q^kt}t^{1/2} E_{-(k-1)\omega},\\
P_{\omega}P_{k\omega} &= P_{(k+1)\omega} + \frac{1-q^k}{1-q^kt} \frac{1-q^{k-1}t^2}{1-q^{k-1}t} P_{(k-1)\omega},
\end{align*}
are walks of type $\vec{m}_{\omega}^{-1} = \pi^\vee$ corresponding to `changes in sheets', and begin in the alcoves $m_{k\omega}^{-1}$ or $m_{k\omega}^{-1}s_1$.
%
\hfill$\diamond$
}\end{example}


\subsection{Type $\fsl_3$}\label{sec.eg2}
Let $\{\varepsilon_1, \varepsilon_2, \varepsilon_3\}$ be an orthonormal basis for $\bbR^3$, and let $\{\varepsilon_1^\vee, \varepsilon_2^\vee, \varepsilon_3^\vee\}$ be its dual basis, where $\langle \varepsilon_i, \varepsilon_j^\vee \rangle = \delta_{ij}$.  
The simple roots and simple coroots are
\begin{eqnarray*}
&\alpha_1 = \varepsilon_1-\varepsilon_2,\quad
\alpha_2 = \varepsilon_2-\varepsilon_3,\quad
\varphi = \varepsilon_1 - \varepsilon_3,\\
&\alpha_1^\vee = \varepsilon_1^\vee-\varepsilon_2^\vee,\quad
\alpha_2^\vee = \varepsilon_2^\vee-\varepsilon_3^\vee,\quad
\varphi^\vee = \varepsilon_1^\vee - \varepsilon_3^\vee,
\end{eqnarray*}
so that the root and coroot lattices are $Q = \bbZ\alpha_1+\bbZ\alpha_2$ and $Q = \bbZ\alpha_1^\vee+\bbZ\alpha_2^\vee$.
The fundamental weights and fundamental coweights are
\begin{eqnarray*}
&\omega_1 = \varepsilon_1 - \hbox{$\frac13$}(\varepsilon_1+\varepsilon_2+\varepsilon_3), \quad
\omega_2 = \varepsilon_1+\varepsilon_2 - \hbox{$\frac23$}(\varepsilon_1+\varepsilon_2+\varepsilon_3),\\
&\omega_1^\vee = \varepsilon_1^\vee - \hbox{$\frac13$}(\varepsilon_1^\vee+\varepsilon_2^\vee+\varepsilon_3^\vee), \quad
\omega_2^\vee = \varepsilon_1^\vee+\varepsilon_2^\vee - \hbox{$\frac23$}(\varepsilon_1^\vee+\varepsilon_2^\vee+\varepsilon_3^\vee),
\end{eqnarray*}
and the weight and coweight lattices are
$\fh_\bbZ^* = \bbZ \omega_1\oplus \bbZ \omega_2$ and $\fh_\bbZ = \bbZ \omega_1^\vee \oplus \bbZ \omega_2^\vee$.  
The group $\Pi^\vee \cong \fh_\bbZ^*/Q$ is the cyclic group of order 3.  Note that $\langle \fh_\bbZ^*, \fh_\bbZ \rangle \subseteq \hbox{$\frac13$}\bbZ$.

The Weyl group of this root system is the symmetric group on three symbols 
$$W_0=\fS_3 = \left\langle  s_1, s_2, \mid s_1s_2s_1=s_2s_1s_2 = s_\varphi, s_i^2=1 \hbox{ for } i =1,2 \right\rangle,$$ 
and the extended affine Weyl group $W^\vee$ is generated by the group $W_0$ and the element $\pi^\vee= x^{\omega_1}s_1s_2$, subject to the relations
\begin{equation}
(\pi^\vee)^3=1, \quad 
\pi^\vee s_0^\vee = s_1\pi^\vee,\quad
\pi^\vee s_1 = s_2\pi^\vee,\quad
s_0s_1s_0 = s_1s_0s_1, \quad
s_0s_2s_0 = s_2s_0s_2,
\end{equation}
where $s_0^\vee= x^\varphi s_\varphi$ and $(\pi^\vee)^2= x^{\omega_2}s_2s_1$.
Alternatively, 
\begin{equation}
W^\vee = \{x^{k_1\omega_1+k_2\omega_2} w \mid k_1,k_2\in\bbZ, w\in \fS_3\}.
\end{equation}

What follows is the alcove picture for the extended affine Weyl group $W^\vee$, showing the correspondence between the alcoves and the elements of $W^\vee$.  The periodic orientation is indicated by $\scriptstyle{+}$ and $\scriptstyle{-}$ on either side of the hyperplanes.  
Since the $\fsl_n$ root system is self-dual, the dual alcove picture is identical.


\newpage
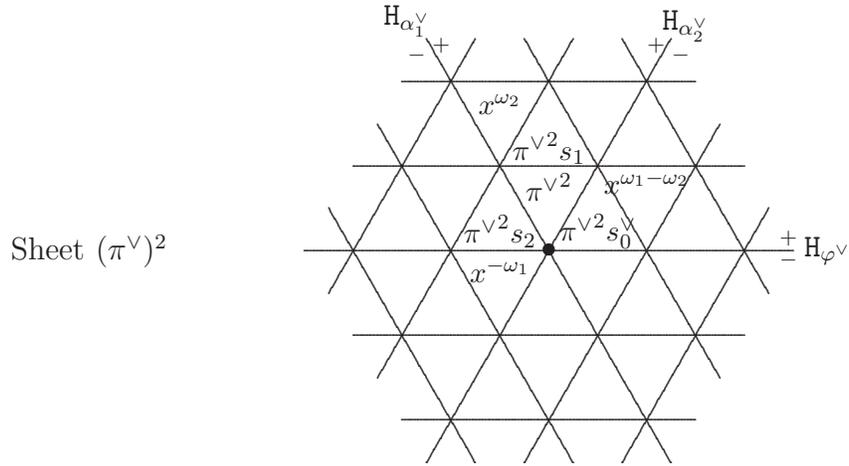
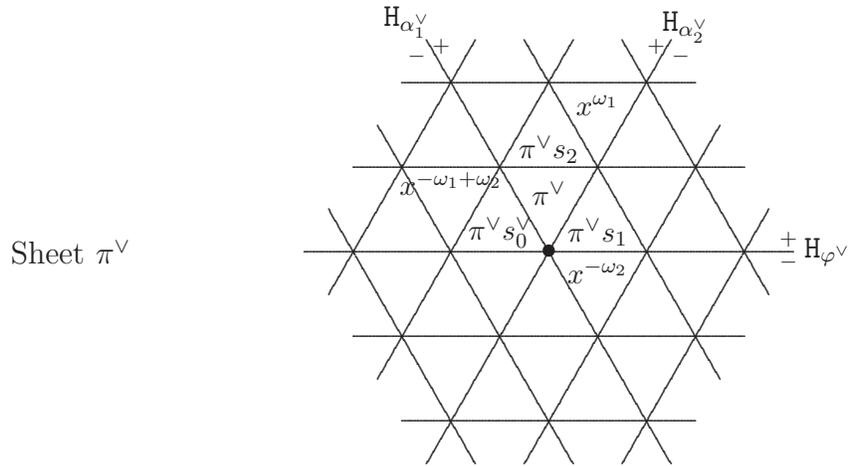
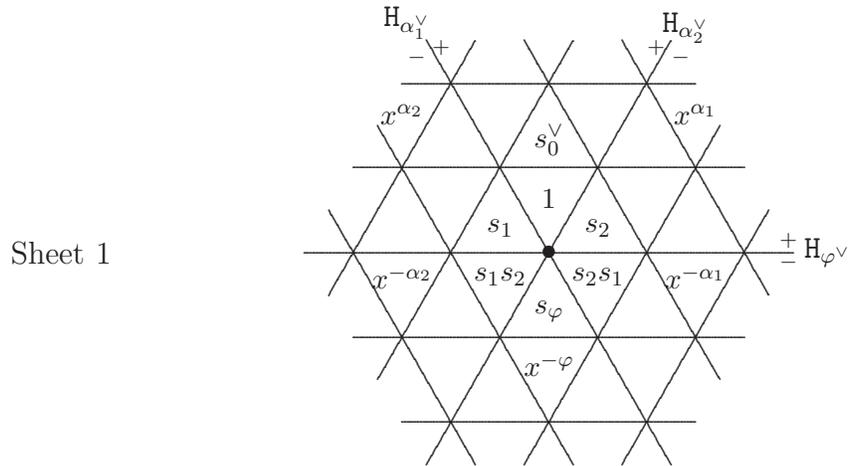
\begin{figure}\label{fig.sl3alcoves}
$$\beginpicture
\setcoordinatesystem units <0.65cm,0.65cm>         
\setplotarea x from -11 to 6, y from -5 to 5    
    \plot 1.5 -4.34 4.5 0.87 /
    \plot -0.5 -4.34 3.5 2.59 /
    \plot -2.5 -4.34 2.5 4.34 /	
    \put{$\ttH_{\alpha_2^\vee}$} at 2.8 4.6
    \plot -3.5 -2.6 0.5 4.34 /		
    \plot -4.5 -0.87 -1.5 4.34 /
    \plot -1.5 -4.34 -4.5 0.87 /
    \plot 0.5 -4.34 -3.5 2.59 /	
    \plot 2.5 -4.34 -2.5 4.34 /	
    \put{$\ttH_{\alpha_1^\vee}$} at -2.9 4.7
    \plot 3.5 -2.6 -0.5 4.34 /
    \plot 4.5  -0.87 1.5 4.34 /
    \plot -3 3.46  3 3.46 /
    \plot -4 1.73 4 1.73 /			
    \plot -5 0  5 0 /					
    \put{$\ttH_{\varphi^\vee}$} at 5.7 0
    \plot -4 -1.73  4 -1.73 /
    \plot -3 -3.46  3 -3.46 /
    \put{$\scriptstyle +$} at 2.2 4.2 \put{$\scriptstyle -$} at 2.7 4
    \put{$\scriptstyle +$} at -2.2 4.2 \put{$\scriptstyle -$} at -2.7 4
    \put{$\scriptstyle +$} at 4.9 0.25 \put{$\scriptstyle -$} at 4.9 -0.2
    \put{Sheet $(\pi^\vee)^2$}[l] at -11 0 
    \put{$\bullet$} at 0 0
    \put{\small ${\pi^\vee}^2$} at 0 1.3
    \put{\small${\pi^\vee}^2s_2$} at -1 0.4
    \put{\small${\pi^\vee}^2s_1$} at 0 2.1
    \put{\small${\pi^\vee}^2s_0^\vee$} at 1 0.4
    \put{\small$x^{\omega_2}$} at -1 3
    \put{\small $x^{\omega_1-\omega_2}$} at 2 1.4
    \put{\small $x^{-\omega_1}$} at -1 -0.4
\endpicture
$$
$$
\beginpicture
\setcoordinatesystem units <0.65cm,0.65cm>         
\setplotarea x from -11 to 6, y from -5 to 5    
    \plot 1.5 -4.34 4.5 0.87 /
    \plot -0.5 -4.34 3.5 2.59 /
    \plot -2.5 -4.34 2.5 4.34 /	
    \put{$\ttH_{\alpha_2^\vee}$} at 2.8 4.6
    \plot -3.5 -2.6 0.5 4.34 /		
    \plot -4.5 -0.87 -1.5 4.34 /
    \plot -1.5 -4.34 -4.5 0.87 /
    \plot 0.5 -4.34 -3.5 2.59 /	
    \plot 2.5 -4.34 -2.5 4.34 /	
    \put{$\ttH_{\alpha_1^\vee}$} at -2.9 4.7
    \plot 3.5 -2.6 -0.5 4.34 /
    \plot 4.5  -0.87 1.5 4.34 /
    \plot -3 3.46  3 3.46 /
    \plot -4 1.73 4 1.73 /			
    \plot -5 0  5 0 /					
    \put{$\ttH_{\varphi^\vee}$} at 5.7 0
    \plot -4 -1.73  4 -1.73 /
    \plot -3 -3.46  3 -3.46 /
    \put{$\scriptstyle +$} at 2.2 4.2 \put{$\scriptstyle -$} at 2.7 4
    \put{$\scriptstyle +$} at -2.2 4.2 \put{$\scriptstyle -$} at -2.7 4
    \put{$\scriptstyle +$} at 4.9 0.25 \put{$\scriptstyle -$} at 4.9 -0.2
    \put{Sheet $\pi^\vee$}[l] at -11 0
    \put{$\bullet$} at 0 0
    \put{\small $\pi^\vee$} at 0 1.2
    \put{\small $\pi^\vee s_1$} at 1 0.4 
    \put{\small $\pi^\vee s_2$} at 0 2.1
    \put{\small $\pi^\vee s_0^\vee$} at -1 0.4
    \put{\small $x^{\omega_1}$} at 1 3
    \put{\small $x^{-\omega_1+\omega_2}$} at -2 1.4
    \put{\small $x^{-\omega_2}$} at 1 -0.4
\endpicture
$$
$$
\beginpicture
\setcoordinatesystem units <0.65cm,0.65cm>         
\setplotarea x from -11 to 6, y from -5 to 5    
    \plot 1.5 -4.34 4.5 0.87 /
    \plot -0.5 -4.34 3.5 2.59 /
    \plot -2.5 -4.34 2.5 4.34 /	
    \put{$\ttH_{\alpha_2^\vee}$} at 2.8 4.6
    \plot -3.5 -2.6 0.5 4.34 /		
    \plot -4.5 -0.87 -1.5 4.34 /
    \plot -1.5 -4.34 -4.5 0.87 /
    \plot 0.5 -4.34 -3.5 2.59 /	
    \plot 2.5 -4.34 -2.5 4.34 /	
    \put{$\ttH_{\alpha_1^\vee}$} at -2.9 4.7
    \plot 3.5 -2.6 -0.5 4.34 /
    \plot 4.5  -0.87 1.5 4.34 /
    \plot -3 3.46  3 3.46 /
    \plot -4 1.73 4 1.73 /			
    \plot -5 0  5 0 /					
    \put{$\ttH_{\varphi^\vee}$} at 5.7 0
    \plot -4 -1.73  4 -1.73 /
    \plot -3 -3.46  3 -3.46 /
    \put{$\scriptstyle +$} at 2.2 4.2 \put{$\scriptstyle -$} at 2.7 4
    \put{$\scriptstyle +$} at -2.2 4.2 \put{$\scriptstyle -$} at -2.7 4
    \put{$\scriptstyle +$} at 4.9 0.25 \put{$\scriptstyle -$} at 4.9 -0.2
    \put{\hbox{Sheet} $1$}[l] at -11 0
    \put{$\bullet$} at 0 0
    \put{\small $1$} at 0 1.1
    \put{\small $s_1$} at -1 0.5 \put{\small $s_2$} at 1 0.5
    \put{\small $s_1s_2$} at -1 -0.5 \put{\small $s_2s_1$} at 1 -0.5
    \put{\small $s_\varphi$} at 0 -1.2
    \put{\small $s_0^\vee$} at 0 2.3
    \put{\small $x^{\alpha_1}$} at 3 2.82
    \put{\small $x^{\alpha_2}$} at -3 2.82
    \put{\small $x^{-\alpha_1}$} at 3 -0.5
    \put{\small $x^{-\alpha_2}$} at -3 -0.5
    \put{\small $x^{-\varphi}$} at 0 -2.3
\endpicture
$$
\end{figure}
%

\begin{example}{\em Using Theorem~\ref{thm.Xtau-X}, we calculate the expansion of $E_{-\alpha_2}$ in the basis of monomials.  The minimal length representative is $m_{-\alpha_2} = s_2s_1s_0^\vee$, and the eight walks of type $\vec{m}_{-\alpha_2}=(2,1,0)$ beginning in the fundamental alcove are:
\begin{center}\begin{tabular}{ccc}
\beginpicture
\setcoordinatesystem units <0.45cm,0.45cm>         
\setplotarea x from -5 to 5, y from -8 to 5   
    \plot 1.5 -4.34 4.5 0.87 /
    \plot -0.5 -4.34 3.5 2.59 /
    \plot -2.5 -4.34 2.5 4.34 /	
    \plot -3.5 -2.6 0.5 4.34 /		
    \plot -4.5 -0.87 -1.5 4.34 /
    \plot -1.5 -4.34 -4.5 0.87 /
    \plot 0.5 -4.34 -3.5 2.59 /	
    \plot 2.5 -4.34 -2.5 4.34 /	
    \plot 3.5 -2.6 -0.5 4.34 /
    \plot 4.5  -0.87 1.5 4.34 /
    \plot -3 3.46  3 3.46 /
    \plot -4 1.73 4 1.73 /			
    \plot -5 0  5 0 /					
    \plot -4 -1.73  4 -1.73 /
    \plot -3 -3.46  3 -3.46 /
    \put{$\bullet$} at 0 0 
    \put{\footnotesize$\displaystyle X^{-\alpha_2}T_1$} at 0 -6
    \setplotsymbol(.)
    \arrow <4pt> [.2,.67] from 0 1.15 to 1 0.57
    \arrow <4pt> [.2,.67] from 1 0.57 to 1 -0.57		
    \arrow <4pt> [.2,.67] from 1 -0.57 to 2 -1.15
    \endpicture
& 
\beginpicture
\setcoordinatesystem units <0.45cm,0.45cm>         
\setplotarea x from -5 to 5, y from -8 to 5    
    \plot 1.5 -4.34 4.5 0.87 /
    \plot -0.5 -4.34 3.5 2.59 /
    \plot -2.5 -4.34 2.5 4.34 /	
    \plot -3.5 -2.6 0.5 4.34 /		
    \plot -4.5 -0.87 -1.5 4.34 /
    \plot -1.5 -4.34 -4.5 0.87 /
    \plot 0.5 -4.34 -3.5 2.59 /	
    \plot 2.5 -4.34 -2.5 4.34 /	
    \plot 3.5 -2.6 -0.5 4.34 /
    \plot 4.5  -0.87 1.5 4.34 /
    \plot -3 3.46  3 3.46 /
    \plot -4 1.73 4 1.73 /			
    \plot -5 0  5 0 /					
    \plot -4 -1.73  4 -1.73 /
    \plot -3 -3.46  3 -3.46 /
    \put{$\bullet$} at 0 0 
    \put{\footnotesize$\displaystyle X^0T_2T_1 
    	\frac{t^{-1/2}-t^{1/2}}{1-Y^{\varphi^\vee-d}}$} at 0 -6
    \setplotsymbol(.)
    \arrow <4pt> [.2,.67] from 0 1.15 to 1 0.57
    \arrow <4pt> [.2,.67] from 1 0.57 to 1 -0.57	
    \plot 1 -0.58 1.5 -0.87 / \plot 1.5 -0.87 1.4 -1.05 / 
    \arrow <4pt> [.2,.67] from 1.4 -1.05 to 0.9 -0.76
    \endpicture
& 
\beginpicture
\setcoordinatesystem units <0.45cm,0.45cm>         
\setplotarea x from -5 to 5, y from -8 to 5    
    \plot 1.5 -4.34 4.5 0.87 /
    \plot -0.5 -4.34 3.5 2.59 /
    \plot -2.5 -4.34 2.5 4.34 /	
    \plot -3.5 -2.6 0.5 4.34 /		
    \plot -4.5 -0.87 -1.5 4.34 /
    \plot -1.5 -4.34 -4.5 0.87 /
    \plot 0.5 -4.34 -3.5 2.59 /	
    \plot 2.5 -4.34 -2.5 4.34 /	
    \plot 3.5 -2.6 -0.5 4.34 /
    \plot 4.5  -0.87 1.5 4.34 /
    \plot -3 3.46  3 3.46 /
    \plot -4 1.73 4 1.73 /			
    \plot -5 0  5 0 /					
    \plot -4 -1.73  4 -1.73 /
    \plot -3 -3.46  3 -3.46 /
    \put{$\bullet$} at 0 0 
    \put{\footnotesize$\displaystyle X^{\alpha_1}T_1T_2  
    		\frac{t^{-1/2}-t^{1/2}}{1-Y^{\alpha_2^\vee-d}}$} at 0 -6
    \setplotsymbol(.)
    \arrow <4pt> [.2,.67] from 0 1.15 to 0.9 0.57
    \plot 0.9 0.57 0.9 0 / \plot 1.1 0 0.9 0 / 
    \arrow <4pt> [.2,.67] from 1.1 0 to 1.1 0.57
    \arrow <4pt> [.2,.67] from 1.1 0.57 to 2 1.15
    \endpicture
\\
\beginpicture
\setcoordinatesystem units <0.45cm,0.45cm>         
\setplotarea x from -5 to 5, y from -8 to 5    
    \plot 1.5 -4.34 4.5 0.87 /
    \plot -0.5 -4.34 3.5 2.59 /
    \plot -2.5 -4.34 2.5 4.34 /	
    \plot -3.5 -2.6 0.5 4.34 /		
    \plot -4.5 -0.87 -1.5 4.34 /
    \plot -1.5 -4.34 -4.5 0.87 /
    \plot 0.5 -4.34 -3.5 2.59 /	
    \plot 2.5 -4.34 -2.5 4.34 /	
    \plot 3.5 -2.6 -0.5 4.34 /
    \plot 4.5  -0.87 1.5 4.34 /
    \plot -3 3.46  3 3.46 /
    \plot -4 1.73 4 1.73 /			
    \plot -5 0  5 0 /					
    \plot -4 -1.73  4 -1.73 /
    \plot -3 -3.46  3 -3.46 /
    \put{$\bullet$} at 0 0 
    \put{\scriptsize$\displaystyle 
    X^0T_2 \frac{(t^{-1/2}-t^{1/2})Y^{\varphi^\vee-d}}
    {1-Y^{\varphi^\vee-d}} 
    \frac{t^{-1/2}-t^{1/2}}{1-Y^{\alpha_2^\vee-d}} $} at 0 -6
    \setplotsymbol(.)
    \arrow <4pt> [.2,.67] from 0 1.15 to 0.9 0.57
    \plot 0.9 0.57 0.9 0 / \plot 1.1 0 0.9 0 / \arrow <4pt> [.2,.67] from 1.1 0 to 1.1 0.57
    \plot 1.1 0.57 1.5 0.86 / \plot 1.4 1.04 1.5 0.86 / \arrow <4pt> [.2,.67] from 1.4 1.04 to 0.95 0.75 
\endpicture
&
\beginpicture
\setcoordinatesystem units <0.45cm,0.45cm>         
\setplotarea x from -5 to 5, y from -8 to 5    
    \plot 1.5 -4.34 4.5 0.87 /
    \plot -0.5 -4.34 3.5 2.59 /
    \plot -2.5 -4.34 2.5 4.34 /	
    \plot -3.5 -2.6 0.5 4.34 /		
    \plot -4.5 -0.87 -1.5 4.34 /
    \plot -1.5 -4.34 -4.5 0.87 /
    \plot 0.5 -4.34 -3.5 2.59 /	
    \plot 2.5 -4.34 -2.5 4.34 /	
    \plot 3.5 -2.6 -0.5 4.34 /
    \plot 4.5  -0.87 1.5 4.34 /
    \plot -3 3.46  3 3.46 /
    \plot -4 1.73 4 1.73 /			
    \plot -5 0  5 0 /					
    \plot -4 -1.73  4 -1.73 /
    \plot -3 -3.46  3 -3.46 /
    \put{$\bullet$} at 0 0 
    \put{\footnotesize$\displaystyle X^{\alpha_2}T_2T_1 \frac{t^{-1/2}-t^{1/2}}{1-Y^{\varphi^\vee-2d}}$} at 0 -6
    \setplotsymbol(.)
    \plot 0 1.15 0.5 0.86 / \plot 0.5 0.86 0.4 0.68 /
    \arrow <4pt> [.2,.67] from 0.4 0.68 to -0.1 0.97
    \arrow <4pt> [.2,.67] from -0.1 0.97 to -1 0.57
    \arrow <4pt> [.2,.67] from -1 0.57 to -2 1.15
    \endpicture
&
\beginpicture
\setcoordinatesystem units <0.45cm,0.45cm>         
\setplotarea x from -5 to 5, y from -8 to 5    
    \plot 1.5 -4.34 4.5 0.87 /
    \plot -0.5 -4.34 3.5 2.59 /
    \plot -2.5 -4.34 2.5 4.34 /	
    \plot -3.5 -2.6 0.5 4.34 /		
    \plot -4.5 -0.87 -1.5 4.34 /
    \plot -1.5 -4.34 -4.5 0.87 /
    \plot 0.5 -4.34 -3.5 2.59 /	
    \plot 2.5 -4.34 -2.5 4.34 /	
    \plot 3.5 -2.6 -0.5 4.34 /
    \plot 4.5  -0.87 1.5 4.34 /
    \plot -3 3.46  3 3.46 /
    \plot -4 1.73 4 1.73 /			
    \plot -5 0  5 0 /					
    \plot -4 -1.73  4 -1.73 /
    \plot -3 -3.46  3 -3.46 /
    \put{$\bullet$} at 0 0 
    \put{\!\!\!\!\!\scriptsize$\displaystyle X^0T_1 
    \frac{t^{-1/2}-t^{1/2}}{1-Y^{\varphi^\vee-2d}}  
    \frac{(t^{-1/2}-t^{1/2})Y^{\varphi^\vee-d}}
    {1-Y^{\varphi^\vee-d}}$} at 0 -6
    \setplotsymbol(.)
    \plot 0 1.15 0.5 0.86 / \plot 0.5 0.86 0.4 0.68 /
    \arrow <4pt> [.2,.67] from 0.4 0.68 to -0.1 0.97
    \arrow <4pt> [.2,.67] from -0.1 0.97 to -0.9 0.57
    \plot -0.9 0.57 -1.5 0.86 / \plot -1.5 0.86 -1.6 0.68 / 
    \arrow <4pt> [.2,.67] from -1.6 0.68 to -1 0.38
    \endpicture
\end{tabular}\end{center}
$$
\beginpicture
\setcoordinatesystem units <0.45cm,0.45cm>         
\setplotarea x from -5 to 5, y from -8 to 5    
    \plot 1.5 -4.34 4.5 0.87 /
    \plot -0.5 -4.34 3.5 2.59 /
    \plot -2.5 -4.34 2.5 4.34 /	
    \plot -3.5 -2.6 0.5 4.34 /		
    \plot -4.5 -0.87 -1.5 4.34 /
    \plot -1.5 -4.34 -4.5 0.87 /
    \plot 0.5 -4.34 -3.5 2.59 /	
    \plot 2.5 -4.34 -2.5 4.34 /	
    \plot 3.5 -2.6 -0.5 4.34 /
    \plot 4.5  -0.87 1.5 4.34 /
    \plot -3 3.46  3 3.46 /
    \plot -4 1.73 4 1.73 /			
    \plot -5 0  5 0 /					
    \plot -4 -1.73  4 -1.73 /
    \plot -3 -3.46  3 -3.46 /
    \put{$\bullet$} at 0 0 
    \put{\footnotesize$\displaystyle
    X^\varphi T_1T_2T_1 \frac{t^{-1/2}-t^{1/2}}{1-Y^{\varphi^\vee-2d}}  
    \frac{t^{-1/2}-t^{1/2}}{1-Y^{\alpha_2^\vee-d}}$} at 0 -6
    \setplotsymbol(.)
    \plot 0 0.8 0.4 0.6 / \plot 0.4 0.6 0.5 0.78 /
    \arrow <4pt> [.2,.67] from 0.5 0.78 to 0 1.05 
    \plot 0 1.05 -0.4 0.83  / \plot -0.4 0.83 -0.5 1.01 / 
    \arrow <4pt> [.2,.67] from -0.5 1.01 to 0 1.3
    \arrow <4pt> [.2,.67] from 0 1.3 to 0 2.4
    \endpicture
\qquad
\beginpicture
\setcoordinatesystem units <0.45cm,0.45cm>         
\setplotarea x from -5 to 5, y from -8 to 5    
    \plot 1.5 -4.34 4.5 0.87 /
    \plot -0.5 -4.34 3.5 2.59 /
    \plot -2.5 -4.34 2.5 4.34 /	
    \plot -3.5 -2.6 0.5 4.34 /		
    \plot -4.5 -0.87 -1.5 4.34 /
    \plot -1.5 -4.34 -4.5 0.87 /
    \plot 0.5 -4.34 -3.5 2.59 /	
    \plot 2.5 -4.34 -2.5 4.34 /	
    \plot 3.5 -2.6 -0.5 4.34 /
    \plot 4.5  -0.87 1.5 4.34 /
    \plot -3 3.46  3 3.46 /
    \plot -4 1.73 4 1.73 /			
    \plot -5 0  5 0 /					
    \plot -4 -1.73  4 -1.73 /
    \plot -3 -3.46  3 -3.46 /
    \put{$\bullet$} at 0 0 
    \put{\footnotesize$\displaystyle 
    X^0 \frac{t^{-1/2}-t^{1/2}}{1-Y^{\varphi^\vee-2d}}  
    \frac{t^{-1/2}-t^{1/2}}{1-Y^{\alpha_2^\vee-d}}
    \frac{(t^{-1/2}-t^{1/2})Y^{\varphi^\vee-d}}
    {1-Y^{\varphi^\vee-d}}$} at 0 -6
    \setplotsymbol(.)
    \plot 0 0.8 0.4 0.6 / \plot 0.4 0.6 0.5 0.78 /
    \arrow <4pt> [.2,.67] from 0.5 0.78 to 0 1 
    \plot 0 1 -0.4 0.7  / \plot -0.4 0.7 -0.5 0.88 / 
    \arrow <4pt> [.2,.67] from -0.5 0.88 to -0.1 1.15
    \plot -0.1 1.15 -0.1 1.72 / \plot -0.1 1.72 0.1 1.72 / 
    \arrow <4pt> [.2,.67] from 0.1 1.72 to 0.1 1.15
    \endpicture
$$
Since $Y^{\varphi^\vee-d}\mathbf{1}= qt^2\mathbf{1}$, $Y^{\alpha_2^\vee-d}\mathbf{1}= qt\mathbf{1}$, and $Y^{\varphi^\vee-2d}\mathbf{1} = q^2t^2\mathbf{1}$, then in the polynomial representation,
\begin{align*}
t^{-\frac12}&E_{-\alpha_2} \mathbf{1}
= X^{-\alpha_2}\mathbf{1} + X^{\alpha_1}\frac{1-t}{1-qt}\mathbf{1} 
	+ X^{\alpha_2}\frac{1-t}{1-q^2t^2}\mathbf{1}
	+ X^\varphi \frac{1-t}{1-q^2t^2}\frac{1-t}{1-qt}\mathbf{1} \\
&+\! X^0\!
	\left(\frac{1-t}{1-qt^2} +\! \frac{1-t}{1-qt^2}\frac{1-t}{1-qt}qt 
	+\! \frac{1-t}{1-q^2t^2}\frac{1-t}{1-qt^2}qt
	+\! \frac{1-t}{1-q^2t^2}\frac{1-t}{1-qt}\frac{1-t}{1-qt^2}q\right)
	\!\mathbf{1}.
\end{align*}
\hfill$\diamond$
}\end{example}

\begin{example}\label{eg.XtoE}{\em
Using Corollary~\ref{cor.X-tau}, we calculate the expansion of $X^{-\alpha_2}$ in the basis of nonsymmetric Macdonald polynomials.  We consider the eight walks of type $\vec{m}_{-\alpha_2}^{-1}=(0,1,2)$ beginning in the fundamental alcove.  Only five of these are contained in the dominant chamber. 
\begin{center}\begin{tabular}{ccc}
\beginpicture
\setcoordinatesystem units <0.45cm,0.45cm>         
\setplotarea x from -5 to 5, y from -7 to 5   
    \plot 1.5 -4.34 4.5 0.87 /
    \plot -0.5 -4.34 3.5 2.59 /
    \plot -2.5 -4.34 2.5 4.34 /	
    \plot -3.5 -2.6 0.5 4.34 /		
    \plot -4.5 -0.87 -1.5 4.34 /
    \plot -1.5 -4.34 -4.5 0.87 /
    \plot 0.5 -4.34 -3.5 2.59 /	
    \plot 2.5 -4.34 -2.5 4.34 /	
    \plot 3.5 -2.6 -0.5 4.34 /
    \plot 4.5  -0.87 1.5 4.34 /
    \plot -3 3.46  3 3.46 /
    \plot -4 1.73 4 1.73 /			
    \plot -5 0  5 0 /					
    \plot -4 -1.73  4 -1.73 /
    \plot -3 -3.46  3 -3.46 /
    \put{$\bullet$} at 0 0
    \put{\footnotesize$\displaystyle
    t^{-1/2}\tau_{m_{-\alpha_2}}^\vee \mathbf{1}
    $} at 0 -6
    \setplotsymbol(.)
    \arrow <4pt> [.2,.67] from 0 1.15 to 0 2.3
    \arrow <4pt> [.2,.67] from 0 2.3 to -1 2.88
    \arrow <4pt> [.2,.67] from -1 2.88 to -1 4.03
    \endpicture
& 
\beginpicture
\setcoordinatesystem units <0.45cm,0.45cm>         
\setplotarea x from -5 to 5, y from -7 to 5    
    \plot 1.5 -4.34 4.5 0.87 /
    \plot -0.5 -4.34 3.5 2.59 /
    \plot -2.5 -4.34 2.5 4.34 /	
    \plot -3.5 -2.6 0.5 4.34 /		
    \plot -4.5 -0.87 -1.5 4.34 /
    \plot -1.5 -4.34 -4.5 0.87 /
    \plot 0.5 -4.34 -3.5 2.59 /	
    \plot 2.5 -4.34 -2.5 4.34 /	
    \plot 3.5 -2.6 -0.5 4.34 /
    \plot 4.5  -0.87 1.5 4.34 /
    \plot -3 3.46  3 3.46 /
    \plot -4 1.73 4 1.73 /			
    \plot -5 0  5 0 /					
    \plot -4 -1.73  4 -1.73 /
    \plot -3 -3.46  3 -3.46 /
    \put{$\bullet$} at 0 0 
    \put{\footnotesize$\displaystyle 
    -t^{-1/2}\tau_{m_{\alpha_2}}^\vee 
    \frac{t^{-1/2}-t^{1/2}}{1-Y^{\varphi^\vee-2d}} \mathbf{1}$} at 0 -6
    \setplotsymbol(.)
    \arrow <4pt> [.2,.67] from 0 1.15 to 0 2.3
    \arrow <4pt> [.2,.67] from 0 2.3 to -0.9 2.88
    \plot -0.9 2.88 -0.9 3.46 / \plot -1.1 3.45 -0.9 3.45 / 
    \arrow <4pt> [.2,.67] from -1.1 3.46 to -1.1 2.88
 \endpicture
& 
\beginpicture
\setcoordinatesystem units <0.45cm,0.45cm>         
\setplotarea x from -5 to 5, y from -7 to 5    
    \plot 1.5 -4.34 4.5 0.87 /
    \plot -0.5 -4.34 3.5 2.59 /
    \plot -2.5 -4.34 2.5 4.34 /	
    \plot -3.5 -2.6 0.5 4.34 /		
    \plot -4.5 -0.87 -1.5 4.34 /
    \plot -1.5 -4.34 -4.5 0.87 /
    \plot 0.5 -4.34 -3.5 2.59 /	
    \plot 2.5 -4.34 -2.5 4.34 /	
    \plot 3.5 -2.6 -0.5 4.34 /
    \plot 4.5  -0.87 1.5 4.34 /
    \plot -3 3.46  3 3.46 /
    \plot -4 1.73 4 1.73 /			
    \plot -5 0  5 0 /					
    \plot -4 -1.73  4 -1.73 /
    \plot -3 -3.46  3 -3.46 /
    \put{$\bullet$} at 0 0 
    \put{\footnotesize$\displaystyle 
    -t^{-1/2}\tau_{m_{\alpha_1}}^\vee 
    \frac{t^{-1/2}-t^{1/2}}{1-Y^{\alpha_2^\vee-d}} \mathbf{1}$} at 0 -6
    \setplotsymbol(.)
    \arrow <4pt> [.2,.67] from 0 1.15 to 0 2.2
    \plot  -0.6 2.42 0 2.12 / \plot -0.5 2.6 -0.6 2.42 / 
    \arrow <4pt> [.2,.67] from -0.5 2.6 to 0.1 2.31 
    \arrow <4pt> [.2,.67] from 0.1 2.31 to 1 2.88
    \endpicture
\end{tabular}
\end{center}
$$
\beginpicture
\setcoordinatesystem units <0.45cm,0.45cm>         
\setplotarea x from -5 to 5, y from -7 to 5    
    \plot 1.5 -4.34 4.5 0.87 /
    \plot -0.5 -4.34 3.5 2.59 /
    \plot -2.5 -4.34 2.5 4.34 /	
    \plot -3.5 -2.6 0.5 4.34 /		
    \plot -4.5 -0.87 -1.5 4.34 /
    \plot -1.5 -4.34 -4.5 0.87 /
    \plot 0.5 -4.34 -3.5 2.59 /	
    \plot 2.5 -4.34 -2.5 4.34 /	
    \plot 3.5 -2.6 -0.5 4.34 /
    \plot 4.5  -0.87 1.5 4.34 /
    \plot -3 3.46  3 3.46 /
    \plot -4 1.73 4 1.73 /			
    \plot -5 0  5 0 /					
    \plot -4 -1.73  4 -1.73 /
    \plot -3 -3.46  3 -3.46 /
    \put{$\bullet$} at 0 0 
    \put{\footnotesize$\displaystyle
    t^{-1/2}\tau_{m_\varphi}^\vee 
    \frac{t^{-1/2}-t^{1/2}}{1-Y^{\alpha_2^\vee-d}} 
    	\frac{t^{-1/2}-t^{1/2}}{1-Y^{\alpha_1^\vee-d}} \mathbf{1}$} at 0 -6
    \setplotsymbol(.)
    \arrow <4pt> [.2,.67] from 0 1.15 to 0 2.1
    \plot  0 2.1 -0.6 2.4  / \plot -0.6 2.4 -0.5 2.58 / 
    	\arrow <4pt> [.2,.67] from -0.5 2.58 to 0 2.3 
    \plot 0 2.3 0.5 2.59 / \plot 0.4 2.77 0.5 2.59 / 
    \arrow <4pt> [.2,.67] from 0.4 2.77 to 0 2.55
\endpicture
\qquad
\beginpicture
\setcoordinatesystem units <0.45cm,0.45cm>         
\setplotarea x from -5 to 5, y from -7 to 5    
    \plot 1.5 -4.34 4.5 0.87 /
    \plot -0.5 -4.34 3.5 2.59 /
    \plot -2.5 -4.34 2.5 4.34 /	
    \plot -3.5 -2.6 0.5 4.34 /		
    \plot -4.5 -0.87 -1.5 4.34 /
    \plot -1.5 -4.34 -4.5 0.87 /
    \plot 0.5 -4.34 -3.5 2.59 /	
    \plot 2.5 -4.34 -2.5 4.34 /	
    \plot 3.5 -2.6 -0.5 4.34 /
    \plot 4.5  -0.87 1.5 4.34 /
    \plot -3 3.46  3 3.46 /
    \plot -4 1.73 4 1.73 /			
    \plot -5 0  5 0 /					
    \plot -4 -1.73  4 -1.73 /
    \plot -3 -3.46  3 -3.46 /
    \put{$\bullet$} at 0 0 
    \put{\footnotesize$\displaystyle
    -t^{-1/2}
    \frac{t^{-1/2}-t^{1/2}}{1-Y^{\varphi^\vee-d}}\ 
    t^{1/2} \ 
    t^{1/2}
    \mathbf{1}$} at 0 -6
    \setplotsymbol(.)
    \plot -0.3 1.1 -0.3 1.72 / \plot -0.3 1.72 -0.05 1.72 /
    \arrow <4pt> [.2,.67] from -0.05 1.72 to -0.05 1 %
    \plot -0.05 1 -0.4 0.75  / \plot -0.4 0.75 -0.27 0.57 / 
    \arrow <4pt> [.2,.67] from -0.27 0.57 to 0.02 0.82 %
    \plot 0.02 0.82 0.33 0.62 / \plot 0.33 0.62 0.4 0.78 / 
    \arrow <4pt> [.2,.67] from 0.4 0.78 to 0.15 1.05 
    \endpicture
$$

Thus,
\begin{align*}
X^{-\alpha_2} 
&= t^{-\frac12} E_{-\alpha_2}
	- \frac{1-t}{1-q^2t^2} t^{-\frac22} E_{\alpha_2} 
	- \frac{1-t}{1-qt}t^{-\frac22} E_{\alpha_1} 
	+ \frac{1-t}{1-qt}\frac{1-t}{1-qt} t^{-\frac32} E_\varphi
	- \frac{1-t}{1-qt^2} E_0.
\end{align*}
}\end{example}



\end{document}